\documentclass[twocolumn]{svjour3} % smallcondensed, smallextended, twocolumn, draft, final

\usepackage[american]{babel}
\usepackage[utf8]{inputenc}
\usepackage{ae,aecompl}
\usepackage[numbers]{natbib}

\usepackage{bm}
\usepackage{url}
\usepackage{color}
\usepackage{booktabs}
\usepackage{graphicx}
\usepackage{subfigure}
\usepackage{float}
\usepackage{stfloats}
\usepackage{tikz}
\usepackage{amsmath} 
\usepackage{amssymb}
\usepackage{commath}
\usepackage{dsfont}

\usepackage{caption}
\journalname{Nonlinear Dynamics}

\begin{document}

% --------------------------------------------------------------
\title{Probabilistic maps on bistable vibration energy harvesters}
%\subtitle{}
%\titlerunning{}

\author{Jo\~ao Pedro Norenberg         \and
        Americo Cunha Jr.              \and
        Samuel da Silva                \and
        \\Paulo Sergio Varoto
}
%\authorrunning{}

\institute{Jo\~ao Pedro Norenberg \at
              S\~ao Paulo State University, Ilha Solteira, SP, Brazil  \\
              ORCID: 0000-0003-3558-4053\\
              \email{jp.norenberg@unesp.br}
           \and
           Americo Cunha Jr. \at
           Rio de Janeiro State University, Rio de Janeiro, RJ, Brazil\\
           ORCID: 0000-0002-8342-0363\\
           \email{americo.cunha@uerj.br}
           \and
            Samuel da Silva \at
              S\~ao Paulo State University, Ilha Solteira, SP, Brazil  \\
              ORCID: 0000-0001-6430-3746\\
              \email{samuel.silva13@unesp.br}
           \and
           Paulo Sergio Varoto \at
           University of S\~ao Paulo, S\~ao Carlos, SP, Brazil \\
           ORCID: 0000-0002-1240-1720\\
           \email{varoto@sc.usp.br}
}

\date{Received: date / Accepted: date}

\maketitle
% --------------------------------------------------------------

% rev by Americo
% --------------------------------------------------------------
\begin{abstract}

This paper analyzes the impact of parametric uncertainties on the dynamics of bistable energy harvesters, focusing on obtaining statistical information about how each parameter's variability affects the energy harvesting process. To model the parametric uncertainties, we use a probability distribution derived from the maximum entropy principle, while polynomial chaos is employed to propagate uncertainty. Conditional probabilities and probability maps are obtained to investigate the effect of uncertainty on harvesting energy. We consider different models of bistable energy harvesters that account for nonlinear piezoelectric coupling and asymmetries. Our findings suggest a higher probability of increasing harvested power in the intrawell motion regime as the excitation frequency increases. In contrast, increasing the excitation amplitude and piezoelectric coupling are more likely to increase power in the chaotic and interwell motion regimes, respectively. An illustrative example is presented to emphasize the importance of investigating the influence when all parameters vary simultaneously.

%This paper aims to analyze the effect of parametric uncertainties in the dynamics of bistable energy harvesters, seeking to obtain statistical information about how the variability of each parameter affects the energy harvesting process. A probability distribution derived from the maximum entropy principle is used to simulate parametric uncertainties, while polynomial chaos is used to propagate uncertainty. Different models of bistable energy harvesters are analyzed, considering nonlinear piezoelectric coupling and asymmetries. According to the findings, there is a greater probability of increasing power harvested in the intrawell motion regime as excitation frequency is increased. However, increasing the excitation amplitude and the piezoelectric coupling showed a greater chance of increasing the power under chaotic and interwell motion regimes, respectively.
 
\end{abstract}
% --------------------------------------------------------------

% rev by Americo
% --------------------------------------------------------------
\section{Introduction}

The conversion of vibrational energy available in the environment into electricity has been extensively explored for powering small electronic components, such as embedded sensors in the Internet of Things (IoT) applications, microelectromechanical systems (MEMS), and nanoelectromechanical systems (NEMS) \cite{Koka_2014,Mahmud2022,Mallick2017,SEOL_2013}. This technology has proven to be promising, particularly in minimizing the costs of battery replacement and disposal, which might potentially have a negative impact on the environment.

Numerous efforts have been dedicated to nonlinear bistable piezo-magnetic-elastic energy harvesters. In contrast to their linear counterparts, the nonlinear harvesters are powerful to generate electricity at a broadband of frequencies. They were first proposed by Cottone et al. \cite{cottone2009p080601} and Erturk et al. \cite{erturk2009p254102}, and have been explored in a number of publications \cite{Lua_2020,19_dAQAQ,Jia2020,Khovanov_2021,Lopes_2019,Mann_2010,NORENBERG2023108542}.

In addition to the complexity due to nonlinear behavior, evaluating the performance of bistable energy harvesting systems can still be a significant challenge when introduced into an environment of uncertainty. Recently, a global sensitivity analysis (GSA) was carried out on bistable energy harvesters in \cite{norenbergNoDy_2022} to verify which parameters have a high or low impact on energy harvesting systems under different scenarios of parametric variability. The findings provided crucial information for a deeper understanding of the system, helped to pinpoint the vital parameters that govern changes in dynamic behavior, and served as an indispensable tool for further robust design, optimization, and response prediction of nonlinear harvesters. However, GSA alone does not inform and explain how the most sensitive parameter affects power harvesting for the better or worse. Using a probabilistic uncertainty quantification (UQ) methodology to model and quantify these impacts due to variability is suitable to comprehend how uncertainties affect the system response quantitatively.

This paper uses UQ analysis to investigate bistable energy harvesting systems addressing nonlinear electromechanical coupling and asymmetries. The GSA aids UQ approaches by serving as a crucial preliminary investigation that helps to develop a more straightforward probabilistic model for the system of interest (see \cite{Nagel2020}). Therefore, this work continues the investigation into the effects of uncertainty on bistable energy harvesters that was started in \cite{norenbergNoDy_2022}. The goal of this work is to investigate the effects of uncertainty in bistable energy harvesters, accounting for the nonlinear electromechanical coupling and asymmetries in the harvester. The analysis provides insights into how uncertainties in the system affect power generation and identifies the key parameters that affect the system's performance. The results of this work may help design and optimize bistable energy harvesters, especially in environments of uncertainty.

The literature provides several studies on uncertainty quantification for energy harvesters. Ali et al. \cite{Ali_2010} conducted pioneering research by evaluating the performance of piezoelectric energy harvesters when subject to uncertainties in natural frequency and damping ratio. The authors calculated the power output and verified the results by comparing them with Monte Carlo simulations. They also optimized specific parameters to maximize harvested power in the presence of uncertainties. Ruiz and Meruane \cite{Ruiz2017_uq_gsa} carried out GSA to examine uncertainty propagation in frequency response functions for unimorph and bimorph piezoelectric energy harvesters. Subsequently, Varoto \cite{Varoto_2019} investigated uncertainties in piezoelectric properties, electrical, geometric, and mechanical boundary conditions. The author performed extensive Monte-Carlo simulations for various frequency ranges since the studied device had multiple natural frequencies.

Various methods have been developed to address the challenges of uncertainty quantification in energy harvesting systems while keeping computational costs manageable. Nanda et al. \cite{Nanda_2015} employed the quadrature method with the maximum entropy principle to assess the impact of uncertain parameters on linear and nonlinear energy harvesting systems. Huang et al. \cite{Huang_2020} used the Chebyshev polynomial approximation to study the dynamics of a nonlinear vibration energy harvester with an uncertain parameter. This method transforms the stochastic energy harvester into a high-dimensional equivalent deterministic system via the Chebyshev polynomial approximation. The mean response of the stochastic energy harvester is then analyzed to understand the stochastic response. The authors found that the random factor could lead to multi-period phenomena, periodic bifurcation behavior, output voltage fluctuation, and changes in subharmonics and superharmonics. In the works \cite{Li_2019,Li_2020}, the authors presented an improved interval extension based on the 2nd-order Taylor series as a new method for uncertainty analysis of a monostable nonlinear energy harvester. The proposed method is suitable for quantifying the uncertainty associated with the excitation frequency. They showed that the output voltage of the nonlinear monostable system is more sensitive to frequency than to the excitation force. Godoy and Trindade \cite{Godoy_2012} studied a cantilever plate with bonded piezoelectric patches and a tip mass serving as an energy harvesting device. The authors considered piezoelectric and dielectric constants of the piezoelectric active layers and the electric circuit equivalent inductance as stochastic parameters.

Optimization studies are also employed in these energy harvesting systems to address parametric uncertainties. Li et al. \cite{Li_2020b} developed a robust optimization method for a nonlinear monostable energy harvester that considers uncertainties. This approach defines a range of variations in the mass, capacitance, and electromechanical coupling coefficient. They obtained an optimal design by maximizing the output voltage center point while minimizing its deviation. On the same direction, Cunha~Jr \cite{cunhajr2021p137} presents a robust numerical framework based on the cross-entropy method, which is capable of obtaining optimal designs even in the presence of noise.

Machine learning techniques, specifically Gaussian processes, have been used recently to establish a parametric relationship between uncertain parameters and the harvested power without the need for governing equations. Chatterjee et al. \cite{Chatterjee_2022} verified this method via direct numerical integration combined with Monte Carlo simulations.

Martins et al. \cite{Martins_2022} presented two methodologies for designing cantilever piezoelectric energy harvesters that account for the presence of uncertain parameters. The methodologies use deterministic and robust optimization to identify optimal designs.

These works in the literature typically focus on carrying out the propagation of uncertainties through the system dynamics, considering joint uncertain effects. However, a comprehensive analysis that examines the specific influence of each uncertain effect on harvester performance is currently lacking. This gap is particularly evident when considering a more realistic scenario where all variables are randomly determined. Although valuable for physical comprehension, this limitation highlights the need for more sophisticated statistical tools to gain a deeper understanding of the phenomena involved in these complex systems. In this work, we aim to fill this gap by analyzing the effect of the most sensitive parameters of energy harvesting systems identified in our previous work \cite{norenbergNoDy_2022} in a comprehensive excitation scenario. We construct a polynomial chaos expansion to obtain probabilistic maps and uncertainty propagation of the recovered energy, and we compute statistics to enhance the energy process. Furthermore, we present an illustrative example to emphasize the significance of calculating probability maps that incorporate the joint variation of parameters.

The manuscript is organized as follows. Section 2 presents the dynamical systems of interest, and Section 3 introduces the probabilistic technique approach. Section 4 details the numerical experiments conducted to analyze the uncertainty effect of the harvesters, followed by a discussion of the results. Finally, Section 5 presents the main conclusions of this study.
% --------------------------------------------------------------

% rev by Americo
% --------------------------------------------------------------
\section{Bistable energy harvesters}
\label{nonlinear_syst_sect}

\begin{figure}
\centering
\subfigure[symmetric harvester]{\includegraphics[scale=0.4]{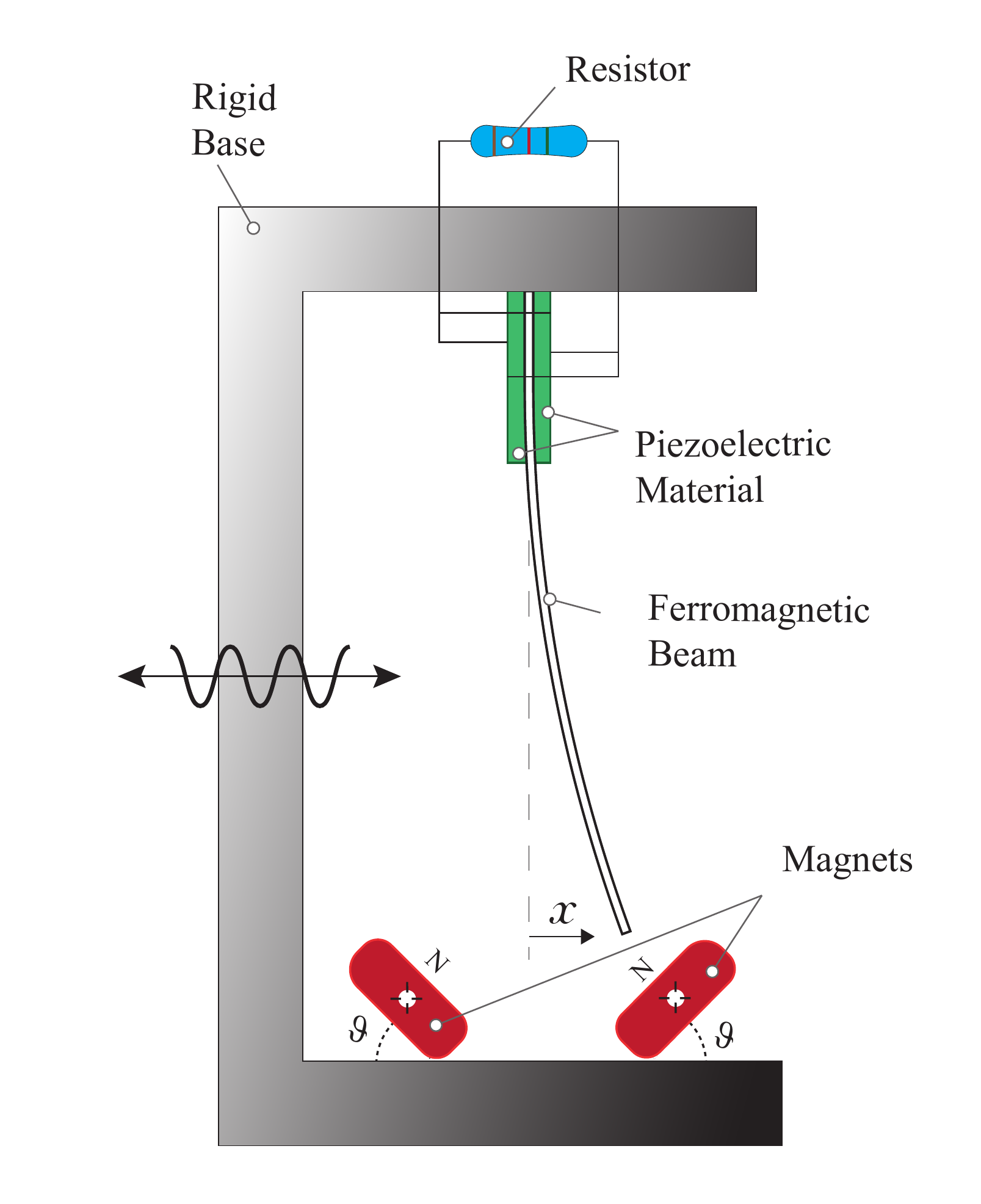}}
\subfigure[asymmetric harvester]{\includegraphics[scale=0.4]{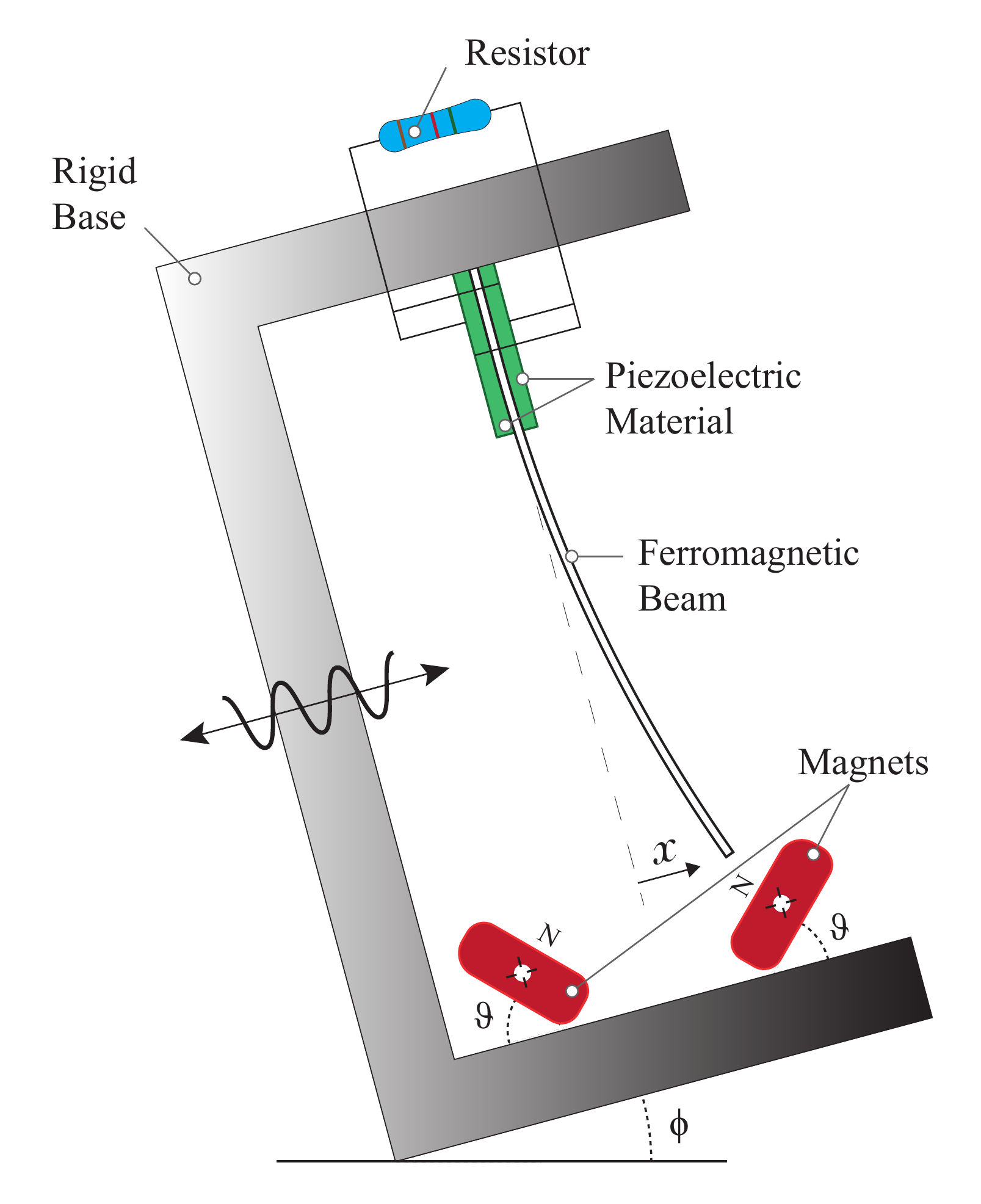}}
\caption{Illustration of the two piezo-magneto-elastic energy harvesting systems: (a) symmetric configuration; (b) asymmetric configuration.}
\label{harvesting_device_fig}
\end{figure}

Figure~\ref{harvesting_device_fig} shows bistable piezo-magneto-elastic energy harvesting systems studied in this work. The systems include symmetric and asymmetric configurations with both linear and nonlinear piezoelectric coupling, which are described in detail in our previous work \cite{norenbergNoDy_2022}. The system depicted in Fig.\ref{harvesting_device_fig}a is a symmetric bistable energy harvester that consists of a vertically oriented clamped-free ferromagnetic elastic beam with piezoelectric layers attached to its highest part and two magnets placed symmetrically on its lower part. An external periodic force excites the rigid base, and the piezoelectric layers convert the kinetic energy into an electrical signal dissipated in the resistor. The system in Fig.\ref{harvesting_device_fig}b is an asymmetric bistable energy harvester, which is placed at inclined surface $\phi$ with non-identical magnets, introducing asymmetries on restoring force. 

According to \cite{12_duToit}, although a piezoelectric energy harvester is commonly modeled with a linear relationship between material strain and the piezoelectric coefficient, neglecting the nonlinear effects of piezoelectricity under high deformation conditions can result in an underestimation of the harvested power. To incorporate a more accurate model, we considered a piezoelectric nonlinear coupling proposed by \cite{13_triplett}, which accounts for a piezoelectric coefficient that depends on the material strain. The governing equations of motion for the bistable energy harvesting model, which takes into account both asymmetries and nonlinear electromechanical coupling as proposed in \cite{norenbergNoDy_2022}, are expressed as
\begin{equation}
    \begin{aligned}
        \ddot{\mathnormal{x}}+2\;\xi\;\dot{\mathnormal{x}}-\frac{1}{2}\;\mathnormal{x}\;(1+2\delta\mathnormal{x}-\mathnormal{x}^2)-(1+\beta\left|x\right|)\chi\; \mathnormal{v} = \\ \mathnormal{f}\;\cos{\left(\Omega\; t\right)}  + \mathnormal{p}\sin{\phi},
    \end{aligned}
	\label{ode_mechanical}
\end{equation}
\begin{equation}
	 \dot{\mathnormal{v}} + \lambda\;\mathnormal{v}+(1+\beta\left|x\right|)\kappa\; \dot x = 0,
	\label{ode_electrical}
\end{equation}
\begin{equation}
	x(0) = x_0, ~ \dot{x}(0) = \dot{x}_0, ~ v(0) = v_0 \, ,
	\label{initial_conditions}
\end{equation}
where $t$ denotes time; $x$ is the modal amplitude of oscillation; $v$ is the voltage in the resistor; $\xi$ is the damping ratio; $f$ is the rigid base oscillation amplitude; $\Omega$ is the external excitation frequency; $\lambda$ is a reciprocal time constant; the piezoelectric coupling terms are represented by $\chi$, in the mechanical equation, and by $\kappa$ in the electrical one; $\delta$ is a coefficient of the quadratic nonlinearity; $\mathnormal{p}$ is the equivalent dimensionless constant of gravity of ferromagnetic beam; $\beta$ is the dimensionless nonlinear coupling term; $x_{0}$ and $\dot{x}_{0}$ are the initial conditions; $v_{0}$ indicates the initial voltage over the resistor. The upper dot is an abbreviation for time derivative. All of these variables are dimensionless. The instantaneous power output at time $t$ is 
\begin{equation}
     P(t) = \lambda \mathnormal{v}(t)^2 \, .
     \label{eq_power}
\end{equation}

\begin{figure}
    \centering
    \subfigure[1-periodic intra-well motion]{\includegraphics[width=0.45\textwidth]{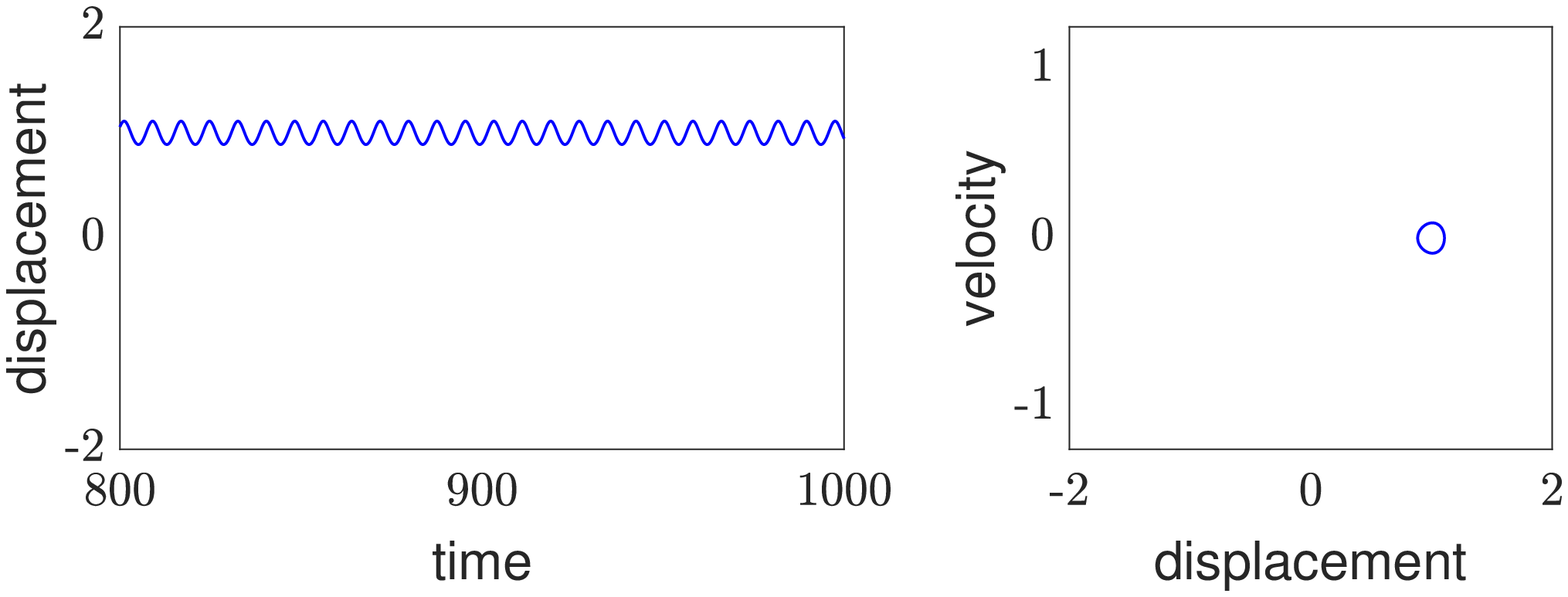}}
    \subfigure[chaotic motion]{\includegraphics[width=0.45\textwidth]{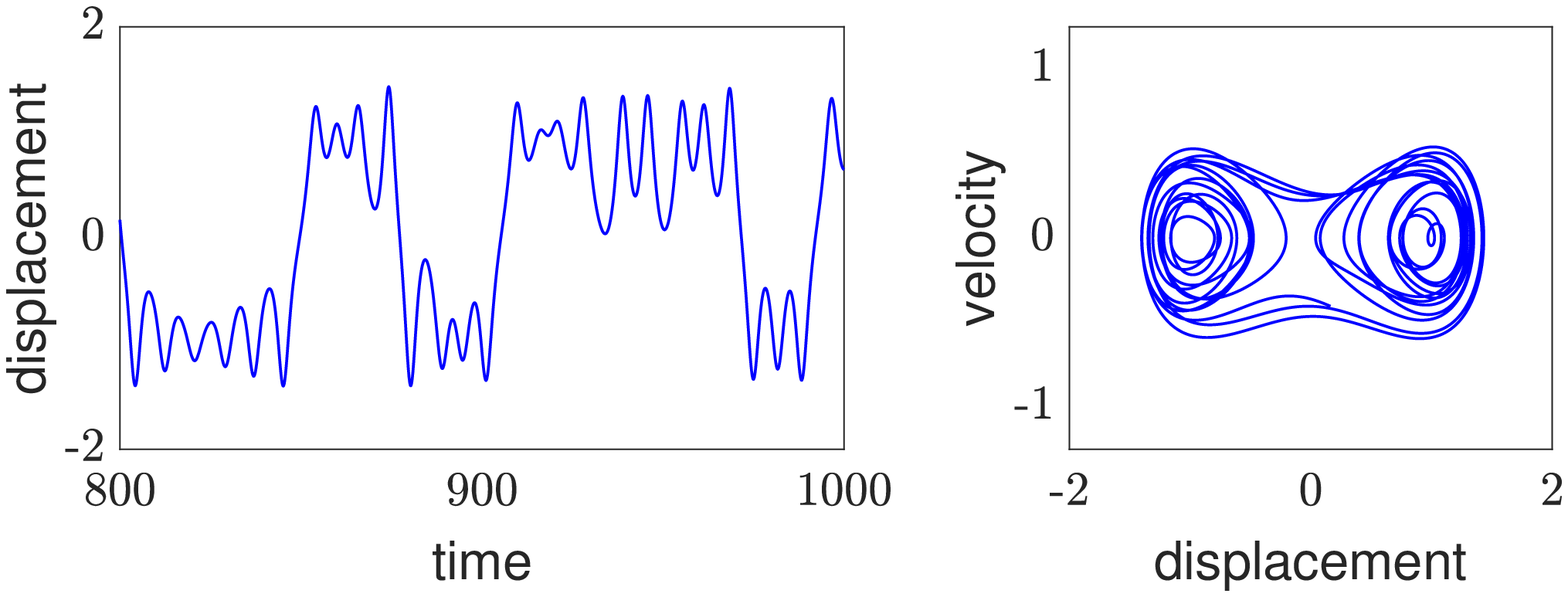}}
    \subfigure[1-periodic inter-well motion]{\includegraphics[width=0.45\textwidth]{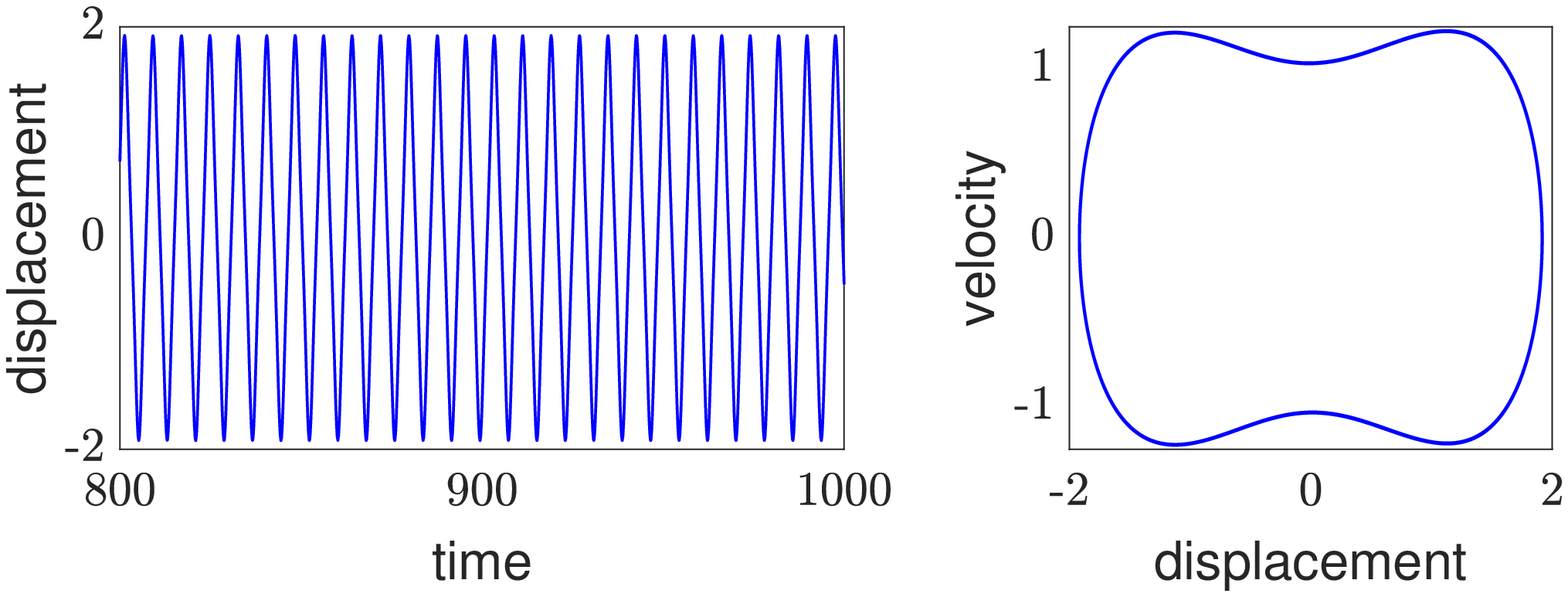}}
    \caption{Typical of dynamic motion by the time-series for the bistable energy harvesting system. The time series in (a) has regular steady-state dynamics at a single-well motion, and in (b) has chaotic steady-state dynamics, while in (c) regular steady-state dynamics at a double-well motion is observed.}
    \label{fig:time}
\end{figure}

This bistable oscillator exhibits three distinct steady-state responses as a result of its nonlinearity. These responses include oscillation with inter-well motion, which can be either chaotic or regular (non-chaotic), as well as intra-well motion with a regular response. They are primarily determined by external excitation conditions. Figure \ref{fig:time} provides a visual representation of these behaviors, showing typical time series and phase portraits for the bistable device.

% --------------------------------------------------------------
% rev by Americo
% --------------------------------------------------------------
\section{Probabilistic approach}

To model the parametric uncertainties of the dynamical systems of interest, a probability space $\left(\Theta; \Sigma; \mathbb{P}\right)$ is considered, where $\Theta$ is a sample space, $\Sigma$ is a $\sigma$-field over $\Theta$, and $\mathbb{P}$ is a probability measure. However, the probability distribution of random parameters cannot be arbitrarily chosen without violating physical principles and creating an inconsistent model. In situations where there is limited information, it can be challenging to determine an unbiased probability distribution function. To address this, the maximum entropy principle provides a formalism that yields the least biased distribution consistent with the available information. According to \cite{Kapur_1992}, this principle offers a rational approach to obtain a suitable joint distribution for the uncertain parameters. The principle aims to choose the least biased distribution, which maximizes entropy while being consistent with the available information about the random parameters \cite{Cunha_uq2017,Soize_2017}.

The entropy of a random variable $X$ is defined as
\begin{equation}
    \mathcal{S}\left( p_X \right) =  -\int_{\mathbb{R}} p_X(x) \, \ln(p_X(x)) \, dx,
    \label{eq:entropy}
\end{equation}
where ${p_X}$ is the probability density function (PDF) of the random variable $X$. Thus, to specify $p_X$, it is necessary to maximize the entropy $\mathcal{S}$, subject to the constraints (known information).

The only known information of the variables for this problem of interest is their supports, Supp$~p_X = \left[a,b\right]$. In other words, we assume that we only know minimum ($a$) and maximum ($b$) values for each parameter based on its physical meaning. Thus, the optimization problem is defined by maximizing Eq.~\eqref{eq:entropy} subjected to $\int_a^b p_X(x) \, dx = 1$.

The Lagrange multipliers strategy is used for minimizing $\mathcal{S}$ subject to the constraint. The Lagrangian function is given by
\begin{equation}
  \begin{aligned}
    \mathcal{L}(p_X,\lambda_0) = -\int_a^b p_X(x) \, \ln(p_X(x)) \, dx \, \\
    - (\lambda_0-1)\left( \int_a^b p_X(x) \, dx - 1\right).
 \end{aligned}
\end{equation}

To obtain the stationary points of $ \mathcal{L}$ in the function of $p_X$ and $\lambda_0$, partial derivatives should be zero. In this way, the extreme conditions are given by
\begin{equation}
    \frac{\partial\mathcal{L}}{\partial p_X}(p_X,\lambda_0) = 0 ~~~ \Rightarrow ~~~ p_X(x) = \mathds{1}_{[a,b]}(x) \, e^{-\lambda_0}, \label{eq:L_px}
\end{equation}
\begin{equation}
    \frac{\partial\mathcal{L}}{\partial \lambda_0}(p_X,\lambda_0) = 0 ~~~ \Rightarrow ~~~ \int_a^b p_X(x) \, dx = 1 \, , 
    \label{eq:L_lab}
\end{equation}
where $\mathds{1}_{[a,b]}(x)$ is a function that returns a value of 1 if it is within the interval $[a,b]$ and 0 otherwise.

Combining Eqs.~\eqref{eq:L_px} and \eqref{eq:L_lab} yields
\begin{equation}
    \int_a^b \mathds{1}_{[a,b]}(x) \, e^{-\lambda_0}~dx = 1~ ~~~ \Rightarrow ~~~ e^{-\lambda_0} = \frac{1}{b-a},
\end{equation}
hence, $X$ has a uniform distribution on $[a, b]$, i.e., 
\begin{equation}
    p_X(x) = \frac{1}{b-a} \, \mathds{1}_{[a,b]}(x).
\end{equation}

Therefore, without known correlation information between the random parameters, the maximum entropy formalism suggests that the parameters are statistically independent. In such a case, a uniform distribution can be used as the marginal distribution for each parameter. Given the available information, this principle ensures that no prior assumptions lead to unbiased distribution.

The mathematical model that predicts the power output for an energy harvesting system can be abstracted as a nonlinear deterministic functional $\mathcal{M}$ that maps the input parameters random vector $\textbf{X}$ into a quantity of interest $\mathcal{Y}$, where $\textbf{X} = (\lambda,\kappa,\mathnormal{f},\Omega,\beta,\delta,\phi)$ represents the input parameters. Based on the sensitivity analysis conducted on our previous paper \cite{norenbergNoDy_2022}, it is worth noting that the parameters $\xi$ and $\chi$ did not influence the power generation. Consequently, they can be treated as having constant values. Thus, the low-dimensional probabilistic model can be represented (explicitly showing the dependence of the random parameters) as follows
\begin{equation}
    \mathcal{Y} = \mathcal{M}(\textbf{X}) .
\end{equation}
This generic notation helps explain the inference methodology used below.

To estimate the probability distribution of the quantity of interest, which in this case is the output power, we need to solve an uncertainty propagation problem \cite{Cunha_uq2017,Soize_2017}. This problem involves determining the distribution of $\mathcal{Y}$ given the probabilistic law of $\textbf{X}$. Several methods can be used to address this problem, including Monte Carlo simulation \cite{cunhajr2014p1355,kroese2011} and polynomial chaos expansion \cite{ghanem2003,xiu2002p619}. Although both methods are effective, the polynomial chaos expansion (PCE) is preferred in this study due to its low computational cost \cite{31_sudret,37_Crestaux,38_Palar}. This method is also accurate and efficient, especially for nonlinear effects in stochastic analysis \cite{34_Oladyshkin,35_SEPAHVAND}.

The PCE is written by
\begin{equation}
    \mathcal{Y} \approx 
    \sum_{\alpha \in \cal{A}} \mathnormal{y}_\alpha\psi_\alpha(\textbf{X}),
\end{equation}
where, $\psi_\alpha$ are multivariate polynomials of $\textbf{X}$, mutually orthonormal with respect to the PDF $p_{\textbf{X}}(\textbf{x})$; and $\mathnormal{y}_\alpha$ are unknown deterministic coefficients. The truncation set $\cal{A} \subset \mathbb{N}^\mathnormal{M}$ is determined from possible multi-indices of multivariate polynomials. The unknown coefficients can be determined using a non-intrusive least-squares regression technique by taking samples from the dynamic system response \cite{31_sudret,37_Crestaux,38_Palar}.

Due to the orthonormality property of the PCE basis, and the fact that $\psi_0 \equiv 1$ and $\mathbb{E}~[\psi_\alpha(\textbf{X})] = 0 ~\forall~\alpha~\neq 0$, the mean value and the variance of the system response can be estimated as
\begin{equation}
    \mathbb{E}\left[\mathcal{Y}\right] \approx \mathnormal{y}_0 \,,
\end{equation}
and
\begin{equation}
    \mathbb{E} \left[\left(\mathcal{Y}- \mathbb{E}\left[ \mathcal{Y} \right]\right)^2\right] \approx \sum_{\stackrel{\alpha\neq0}{\alpha \in \cal A}}\mathnormal{y}_\alpha^2.
\end{equation}

After defining the PCE expansion, ${k}$ independent samples of $\textbf{X}$ are drawn from its distribution analytically and without significative computational cost. Each sample is given as input to the model $\mathcal{M}$, resulting in a set of possible realizations for the quantity of interest
\begin{eqnarray}
\begin{aligned}
    \mathcal{Y}^{(1)} &= \mathcal{M}( \textbf{X}^{(1)}) \\ 
    \mathcal{Y}^{(2)} &= \mathcal{M}( \textbf{X}^{(2)}) \\
    \vdots \;\;\;\; & \;\;\;\;\;\;\;\;\;\;\vdots  \\
    \mathcal{Y}^{(\mathnormal{k})} &= \mathcal{M}( \textbf{X}^{(\mathnormal{k})}) 
\end{aligned}
    % \mathcal{Y}^{(1)},\mathcal{Y}^{(2)},\dots,\mathcal{Y}^{(\mathnormal{k})} = \mathcal{M}( \mathnormal{\Theta}^{(1)},\xi,\chi), \mathcal{M}( \mathnormal{\Theta}^{(2)},\xi,\chi), \dots, \mathcal{M}( \mathnormal{\Theta}^{(\mathnormal{k})},\xi,\chi)
\end{eqnarray}
where $\mathcal{Y}^{(\mathnormal{k})}$ samples are used to estimate statistics of $\mathcal{Y}$ non-parametrically, i.e., without assumptions about the shape of its PDF \cite{Wasserman_2007}. In this work, due to its simplicity and effectiveness, the technique used to estimate the probability density function is the kernel density estimator \cite{Wasserman_2007}. 
% --------------------------------------------------------------

% rev by Americo
% --------------------------------------------------------------
\section{Results and discussion}

In this study, the numerical results are presented for three different bistable models: the symmetric bistable energy harvester with linear piezoelectric coupling, the same system with nonlinear piezoelectric coupling, and the asymmetric bistable energy harvesting system. Their physical parameters are subject to a uniform distribution, as described in the previous section, with a variance coefficient of 20\% around their nominal values. The nominal values for the parameters are assumed to be as follows: $\xi=0.01$, $\chi=0.05$, $\lambda=0.05$, $\kappa=0.5$, $\Omega=0.8$, $\beta=1$,  $\delta = 0.15$ and $\phi=10^\circ$. For the excitation amplitude, various nominal values are used to explore different types of dynamic oscillations, such as the interwell and intrawell motion shown in Fig.~\ref{fig:time}.
Similarly to the other variables, the excitation amplitude is also uniformly varied around its nominal value with a variance coefficient of 20\% in each situation. The same methodology was applied to each model, allowing for an objective comparison of the findings. Finally, we provide an illustrative example to demonstrate the implications of considering the joint variation of parameters and emphasize the importance of calculating probability maps.

\subsection{Symmetric bistable energy harvester with linear piezoelectric coupling}

\begin{figure*}
    \centering
    \subfigure[$\mathnormal{f}=0.041$]{\includegraphics[width=0.32\textwidth]{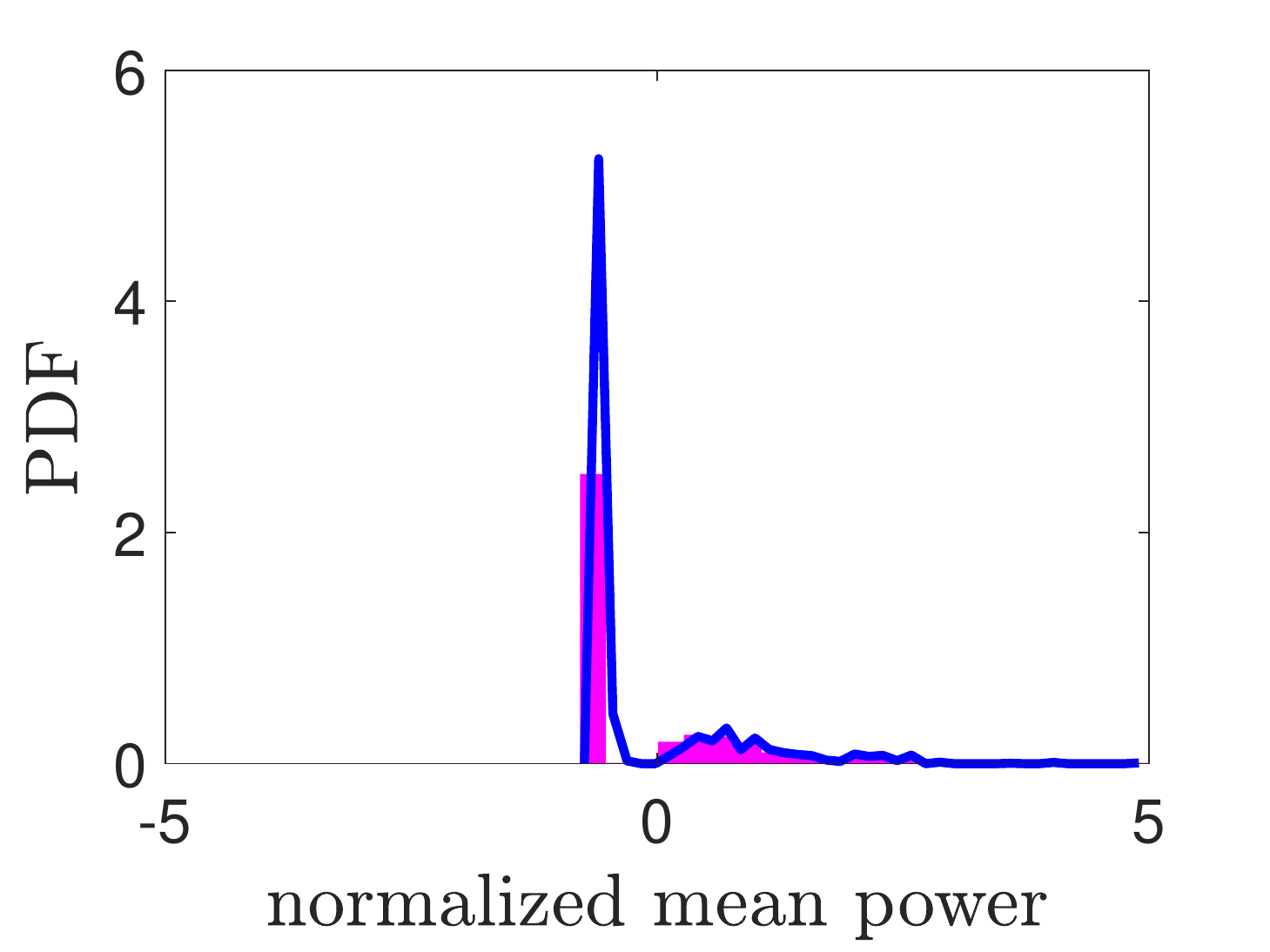}}
    \subfigure[$\mathnormal{f}=0.060$]{\includegraphics[width=0.32\textwidth]{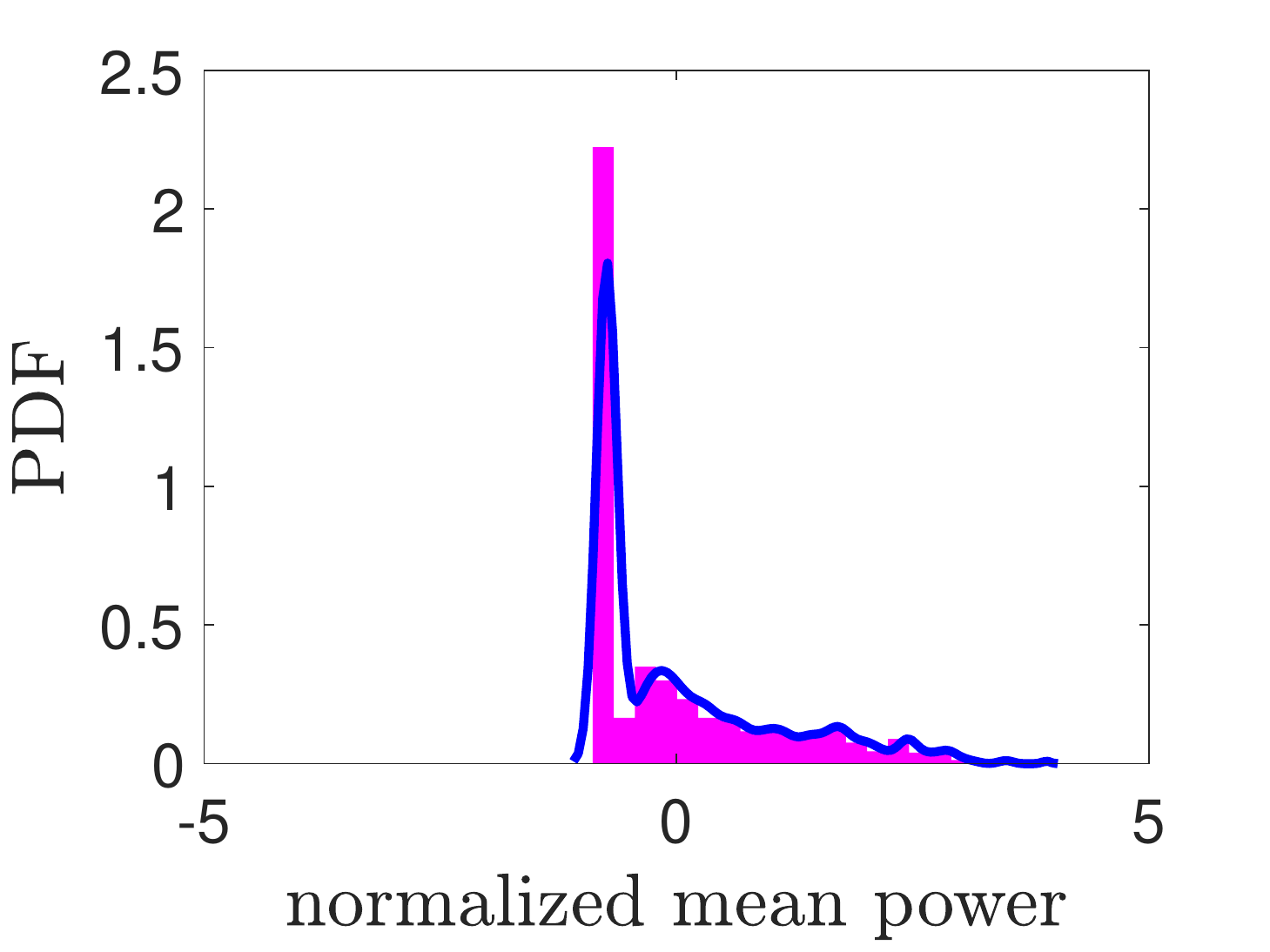}}
    \subfigure[$\mathnormal{f}=0.083$]{\includegraphics[width=0.32\textwidth]{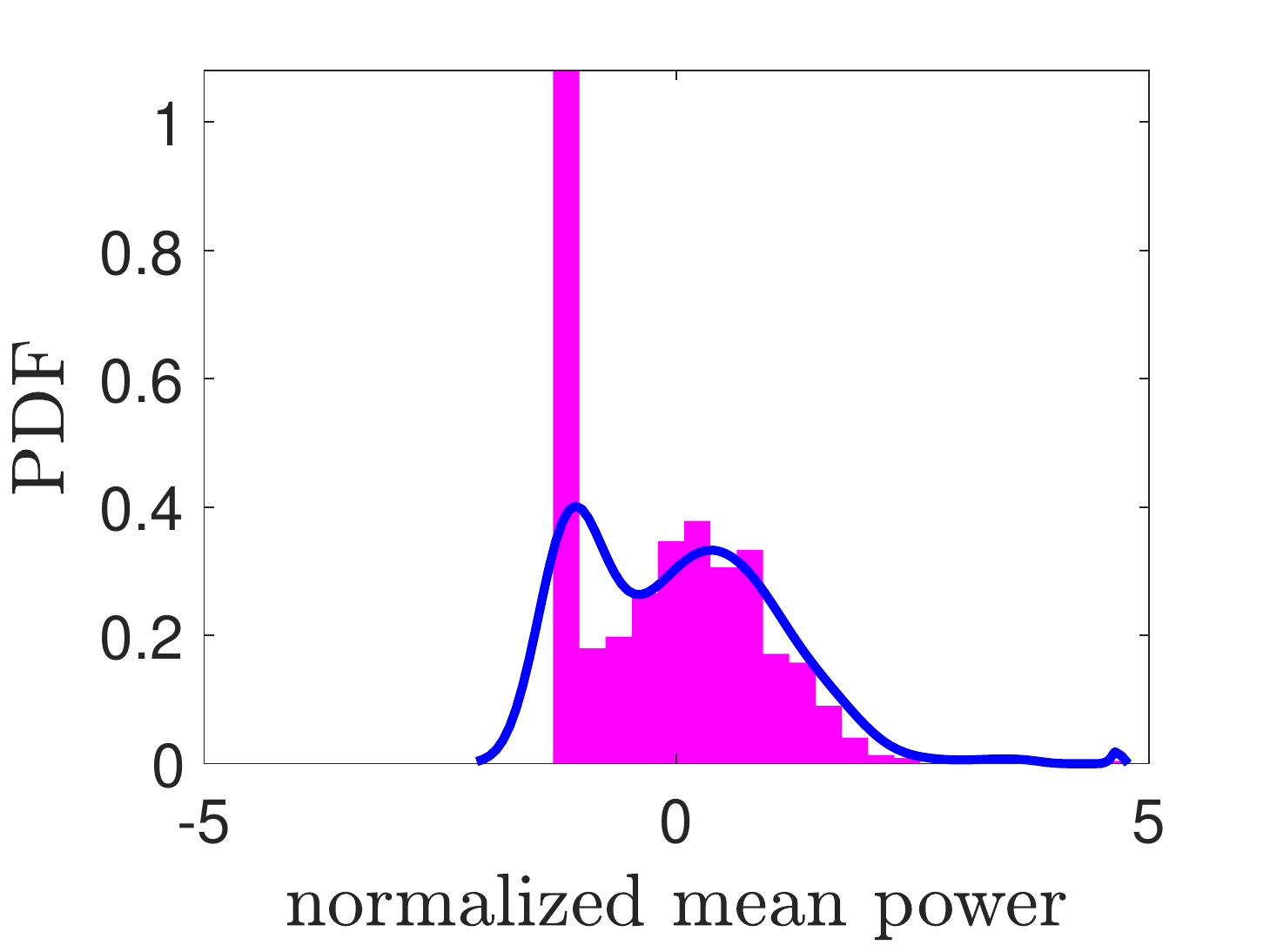}}
    \subfigure[$\mathnormal{f}=0.091$]{\includegraphics[width=0.32\textwidth]{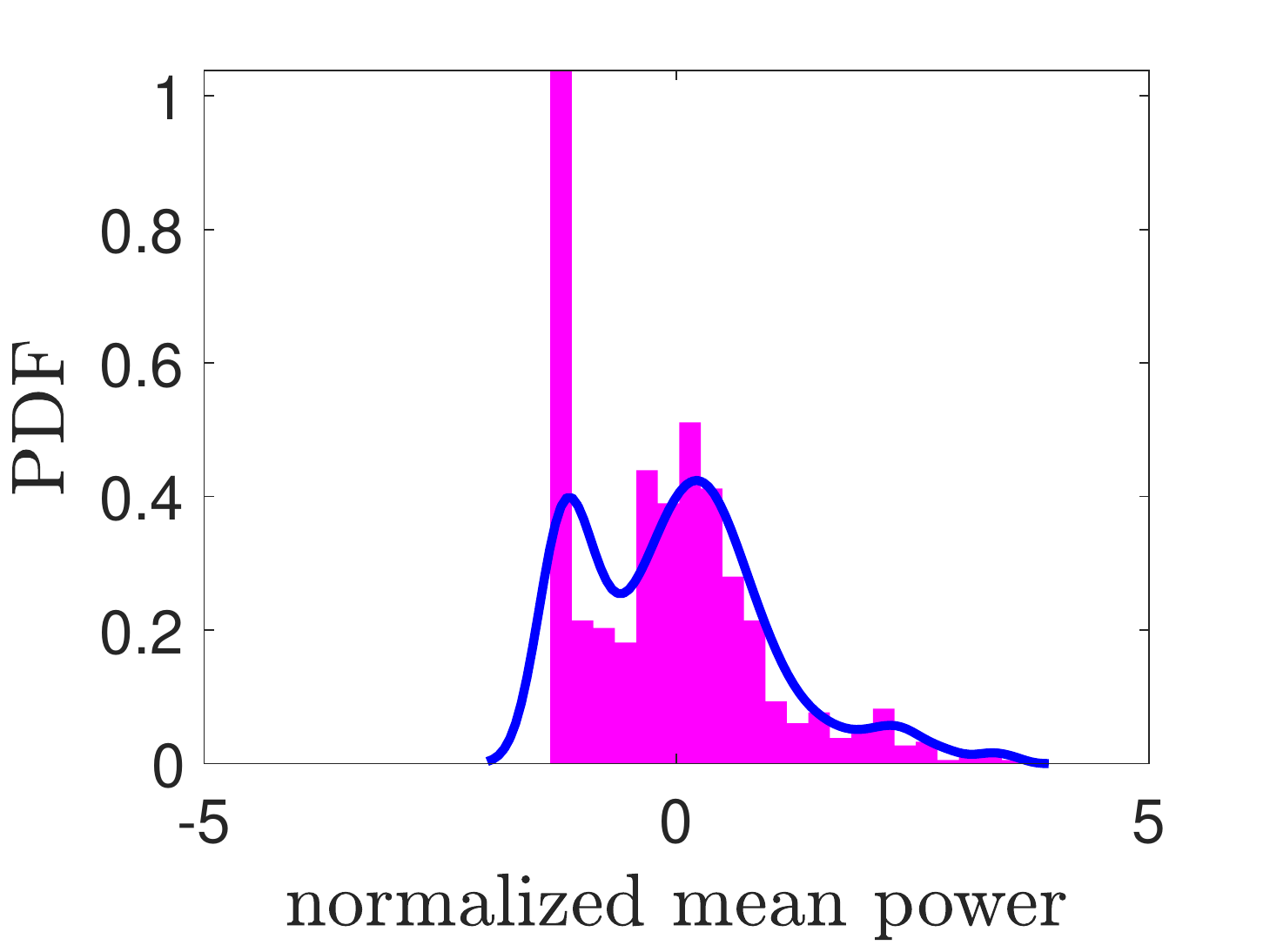}}
    \subfigure[$\mathnormal{f}=0.105$]{\includegraphics[width=0.32\textwidth]{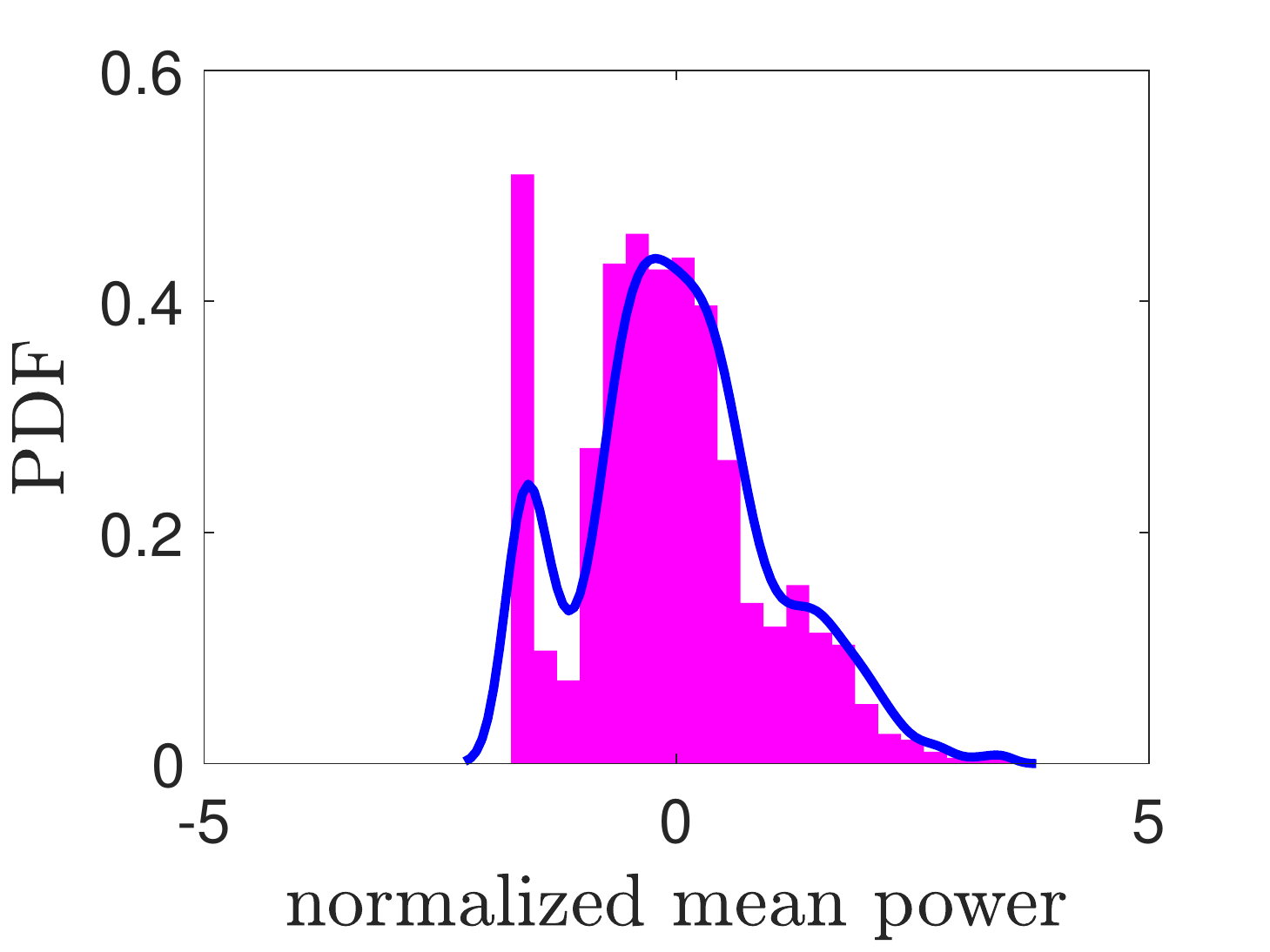}}
    \subfigure[$\mathnormal{f}=0.115$]{\includegraphics[width=0.32\textwidth]{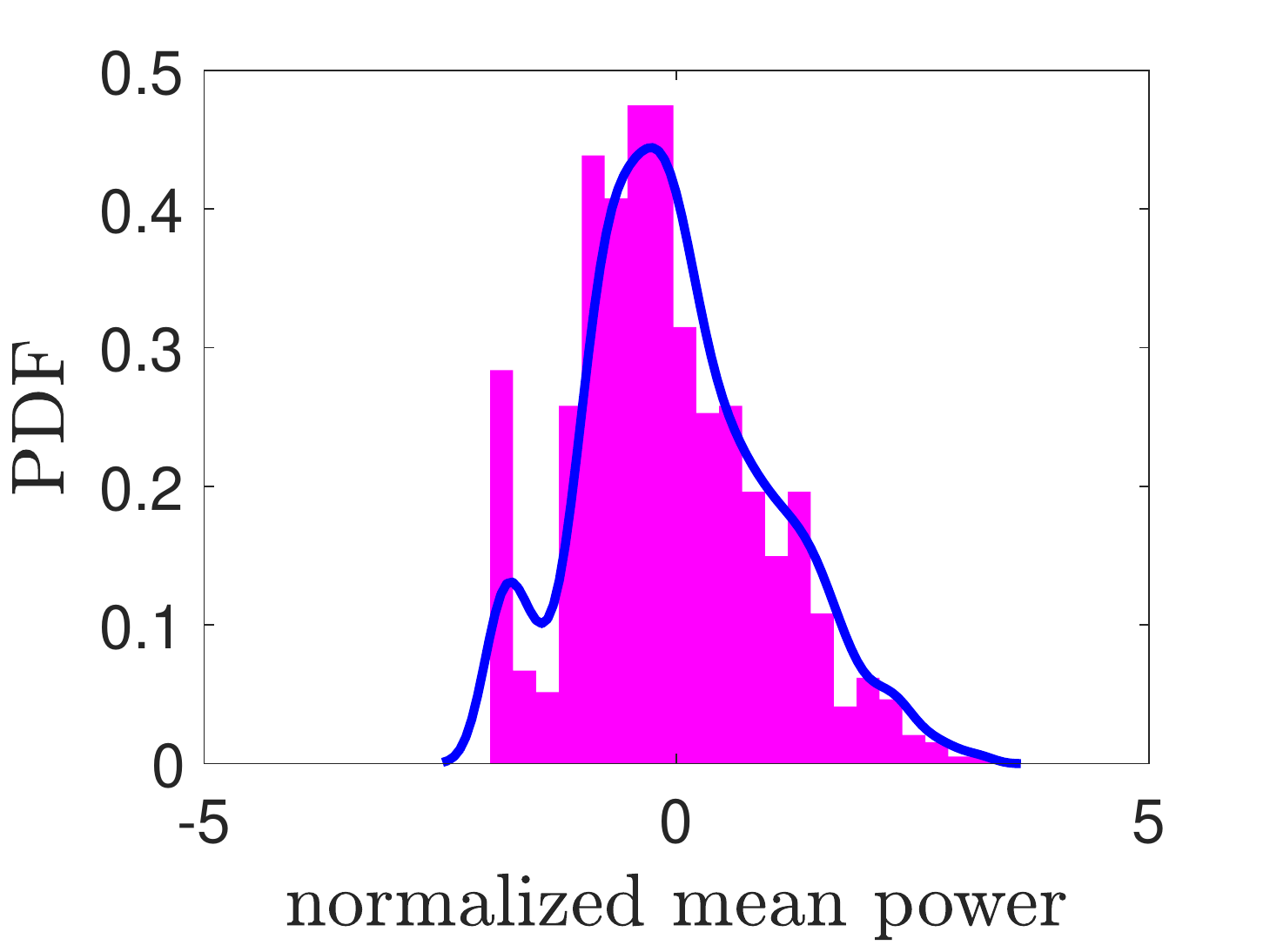}}
    \subfigure[$\mathnormal{f}=0.147$]{\includegraphics[width=0.32\textwidth]{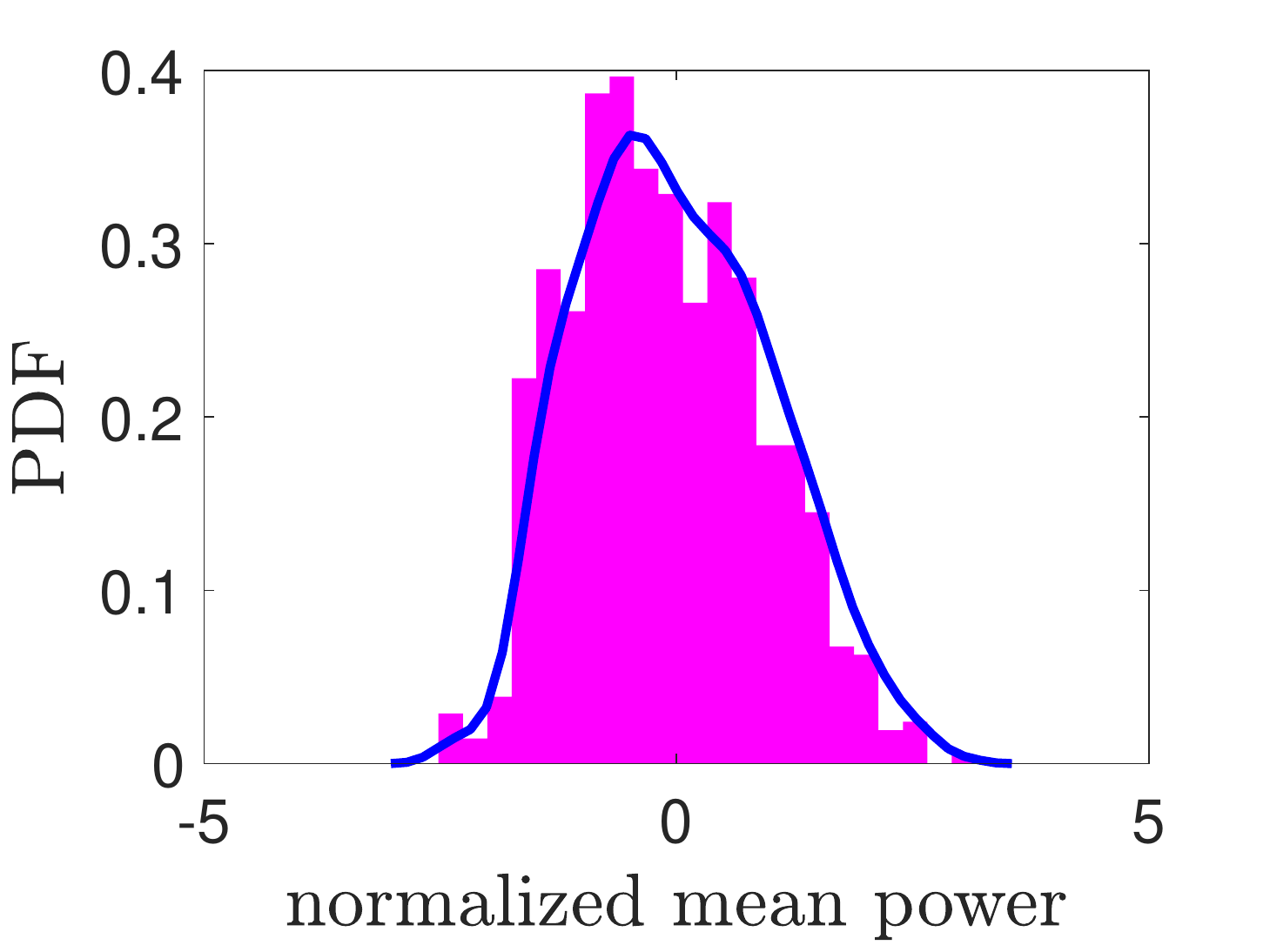}}
    \subfigure[$\mathnormal{f}=0.200$]{\includegraphics[width=0.32\textwidth]{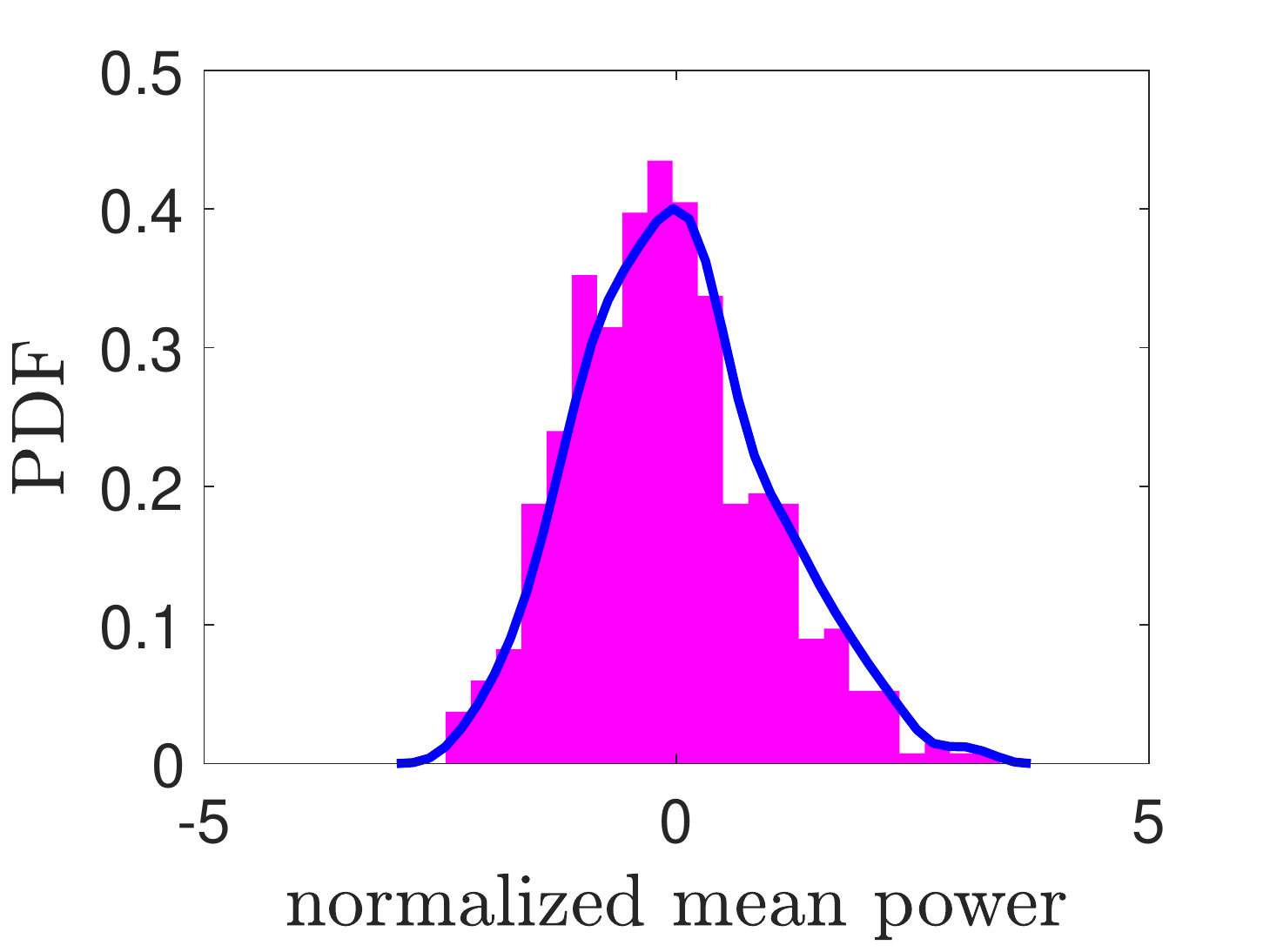}}
    \subfigure[$\mathnormal{f}=0.250$]{\includegraphics[width=0.32\textwidth]{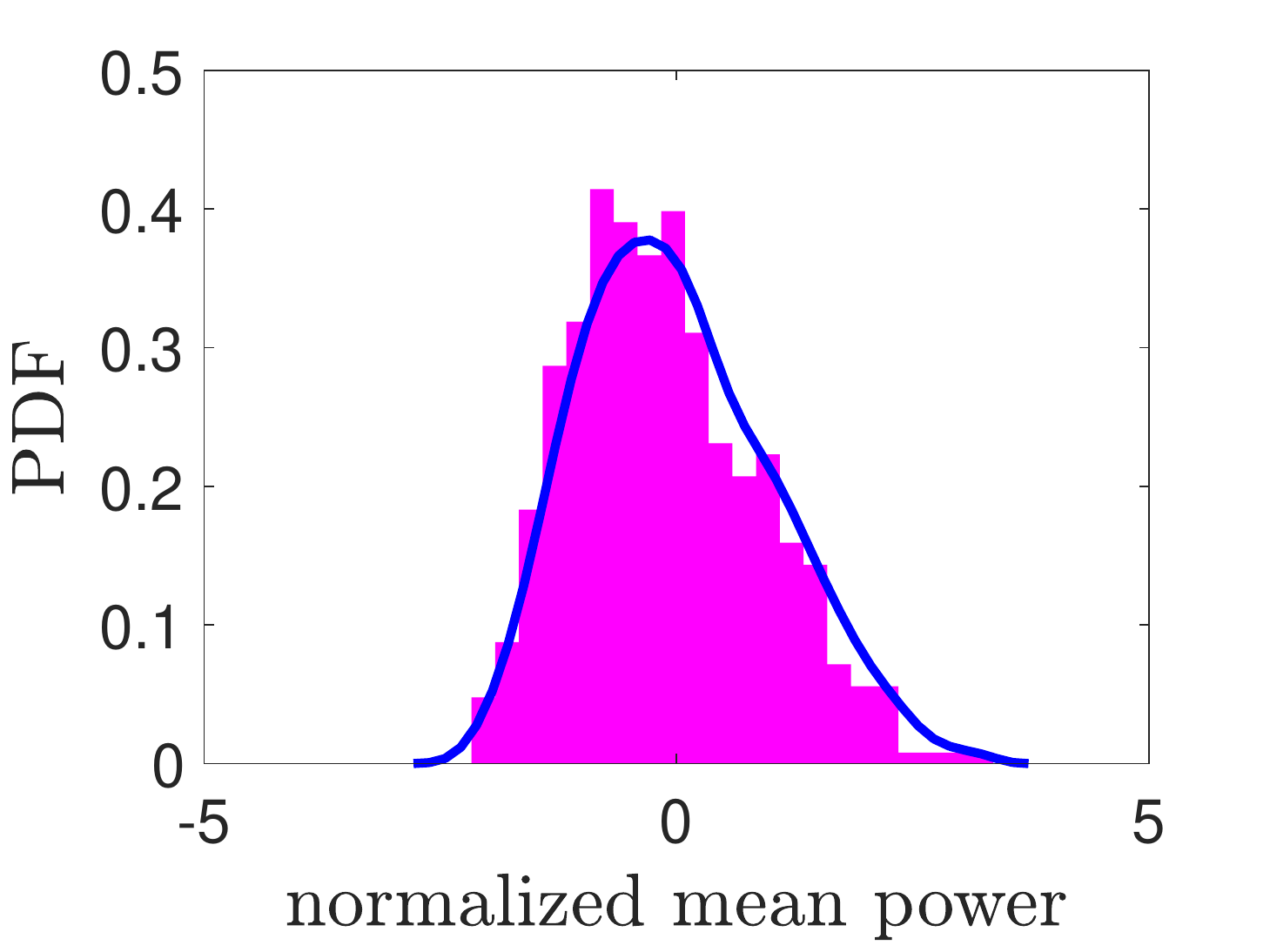}}
    \caption{Probability density function of the normalized mean power for the symmetric energy harvester model with linear piezoelectric coupling under different excitation amplitudes. The kernel density function is represented by the blue line. }
    \label{fig:pdf_pmeh}
\end{figure*}

The findings of the symmetric bistable energy harvester with linear piezoelectric coupling ($\delta = 0$, $\phi = 0$ and $\beta = 0$) are presented in Figure~\ref{fig:pdf_pmeh}, where the histograms and the PDF of the normalized\footnote{Normalization here means zero mean and unit standard deviation.} mean output power are shown for a range of the amplitude of excitation. These results show a bimodal distribution for low values of $\mathnormal{f}$ ($<0.115$), indicating smaller mean power values caused by monostable vibrations. As the amplitude of excitation increases, the second peak (positive values) becomes more prominent until a unimodal distribution is obtained ($\mathnormal{f}>0.115$), where bistable vibrations occur. The nature of the distribution is highly dependent on the system's dynamic behavior and can provide insights into its underlying characteristics, such as monostable and bistable oscillations.

\begin{figure*}
    \centering
    \includegraphics[width=1\textwidth]{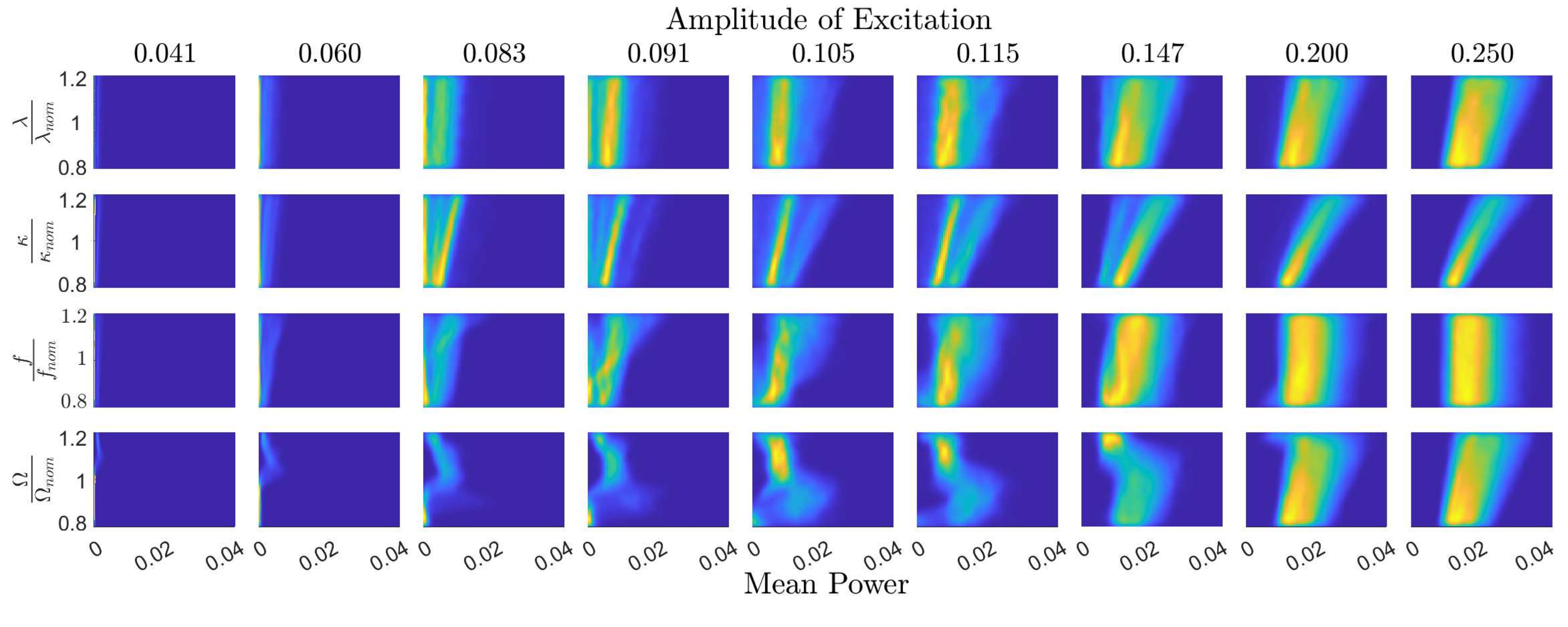}
    \caption{Joint-CDF of mean power conditioned on each parameter of interest ($\lambda$, $\kappa$, $\mathnormal{f}$, $\Omega$) under different values of excitation amplitude for the symmetric model with linear piezoelectric coupling.}
    \label{fig:joint_BEH}
\end{figure*}

Figure~\ref{fig:joint_BEH} displays the contour map of the joint-CDF (joint cumulative distribution function) of the mean power conditioned on each parameter of interest under different nominal excitation conditions. This figure shows the correlation of each parameter with the mean power when all parameters are random. For low amplitudes of excitation ($\mathnormal{f} = 0.041$ and $0.060$), the system exhibits low variability, and the mean power is not significantly affected by parametric changes. Only a slight variation is recorded when $\Omega$ is larger than its nominal value. As the amplitude of excitation increases, the system undergoes significant changes, mainly for $\Omega$ and $\mathnormal{f}$, where discontinuities are observed. Additionally, higher values of $\kappa$ result in a power increase, while the same does not hold for $\lambda$. When the amplitude of excitation is high ($\mathnormal{f} = 0.200$ and $0.250$), the effect of $\kappa$ becomes more pronounced, leading to a considerable increase in power, indicating a positive correlation. The same is observed for $\lambda$, but not for $\mathnormal{f}$ and $\Omega$. Therefore these results demonstrated the correlation between the variable and the mean power.

\begin{figure*}
    \centering
    \includegraphics[width=0.8\textwidth]{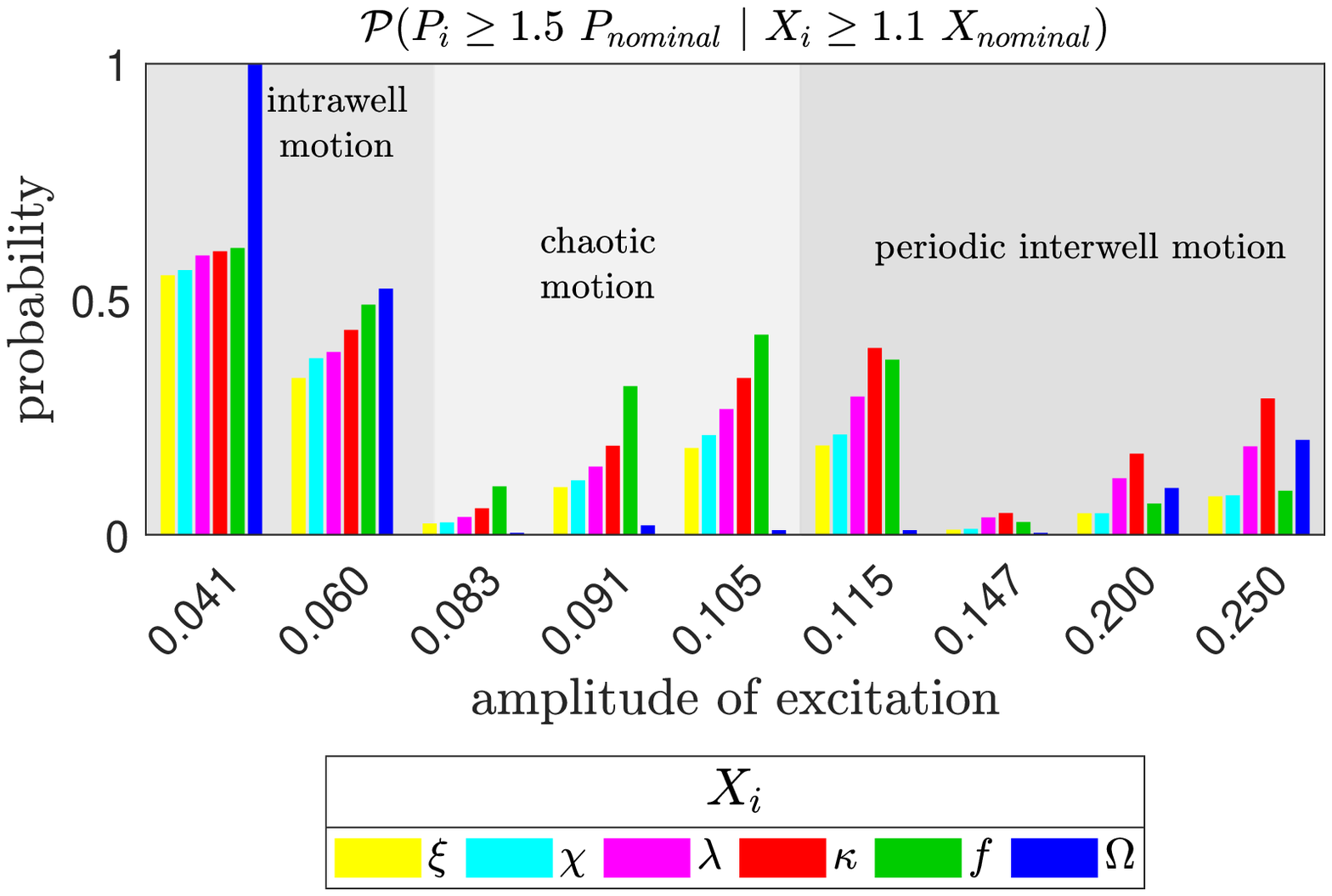}
    \vspace{-1.5cm}
    \caption{Probability of increasing the nominal mean power by 50\% as parameter $X_\mathnormal{i}$ is increased by 10\%, plotted against the excitation amplitude, for the symmetric model with linear piezoelectric coupling.}
    \label{fig:prob_BEH}
\end{figure*}

Figure~\ref{fig:prob_BEH} shows the conditional probability of increasing mean power by 50\% of the nominal power given a 10\% increase in a parameter of interest. This analysis provides insight into the probability of improving energy harvesting as a quantitative improvement. The analysis of the obtained results were categorized based on the dynamic behavior of the system, specifically into three distinct types: intrawell motion, chaotic motion, and periodic interwell motion. This dynamic classification can be seen in more detail in \cite{norenbergNoDy_2022,NORENBERG2023108542}. In the intrawell motion region, increasing $\Omega$ by 10\% of its nominal value is crucial to achieving over 80\% probability of increasing power. In chaotic motion regions, increasing excitation amplitude yields a 40\% probability of increasing mean power, but a high $\Omega$ value reduces the probability of improving energy harvesting. This is because an increased frequency in chaotic regions leads to intrawell behavior that generates less energy. In interwell motion regions, $\kappa$ is the key parameter to generate more energy, as increasing $\kappa$ generates a 20\% chance of improving the harvesting process.

\begin{figure*}
    \centering
    \subfigure[$\mathnormal{f}=0.041$ for $\lambda$]{\includegraphics[width=0.32\textwidth]{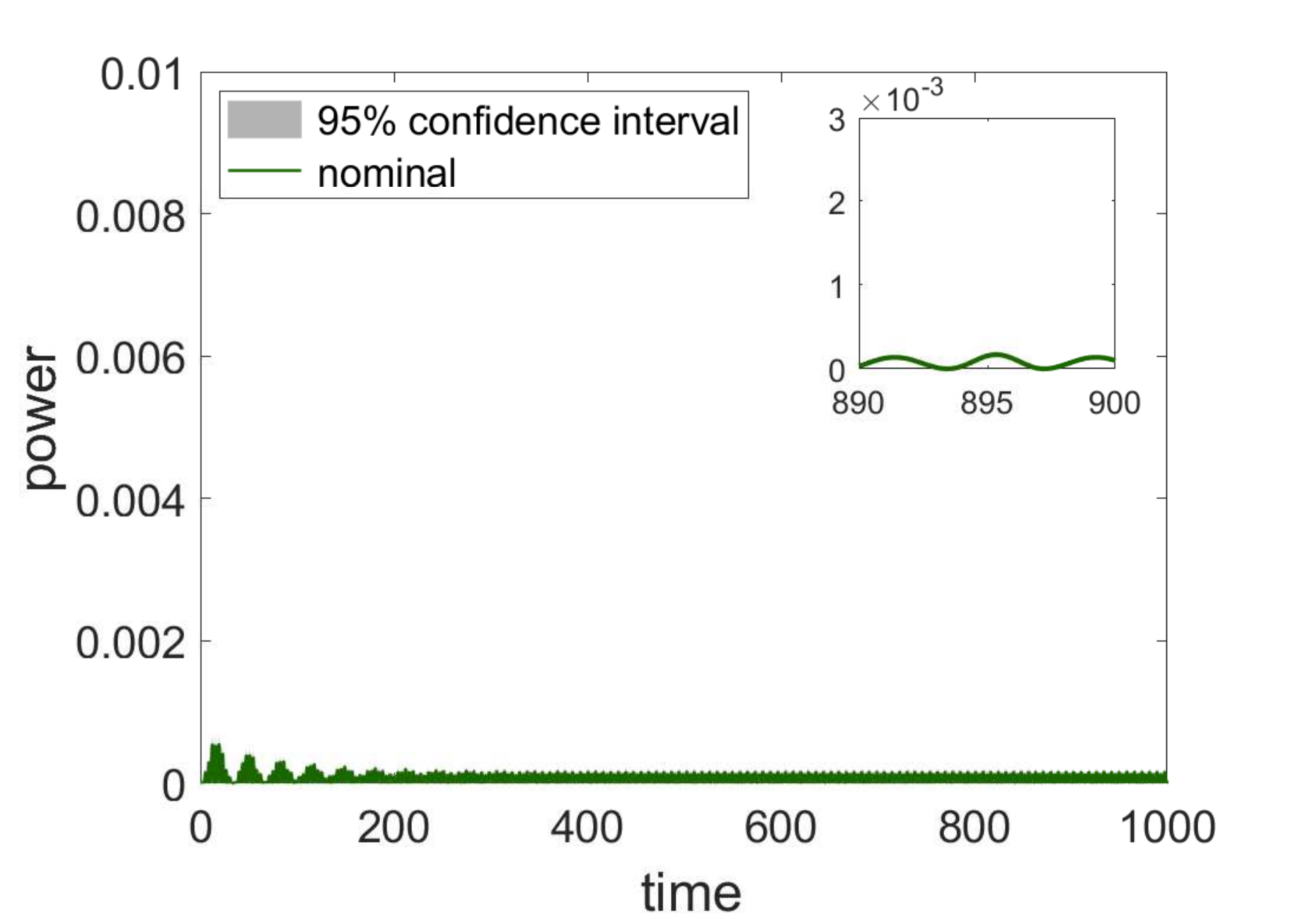}}
    \subfigure[$\mathnormal{f}=0.091$ for $\lambda$]{\includegraphics[width=0.32\textwidth]{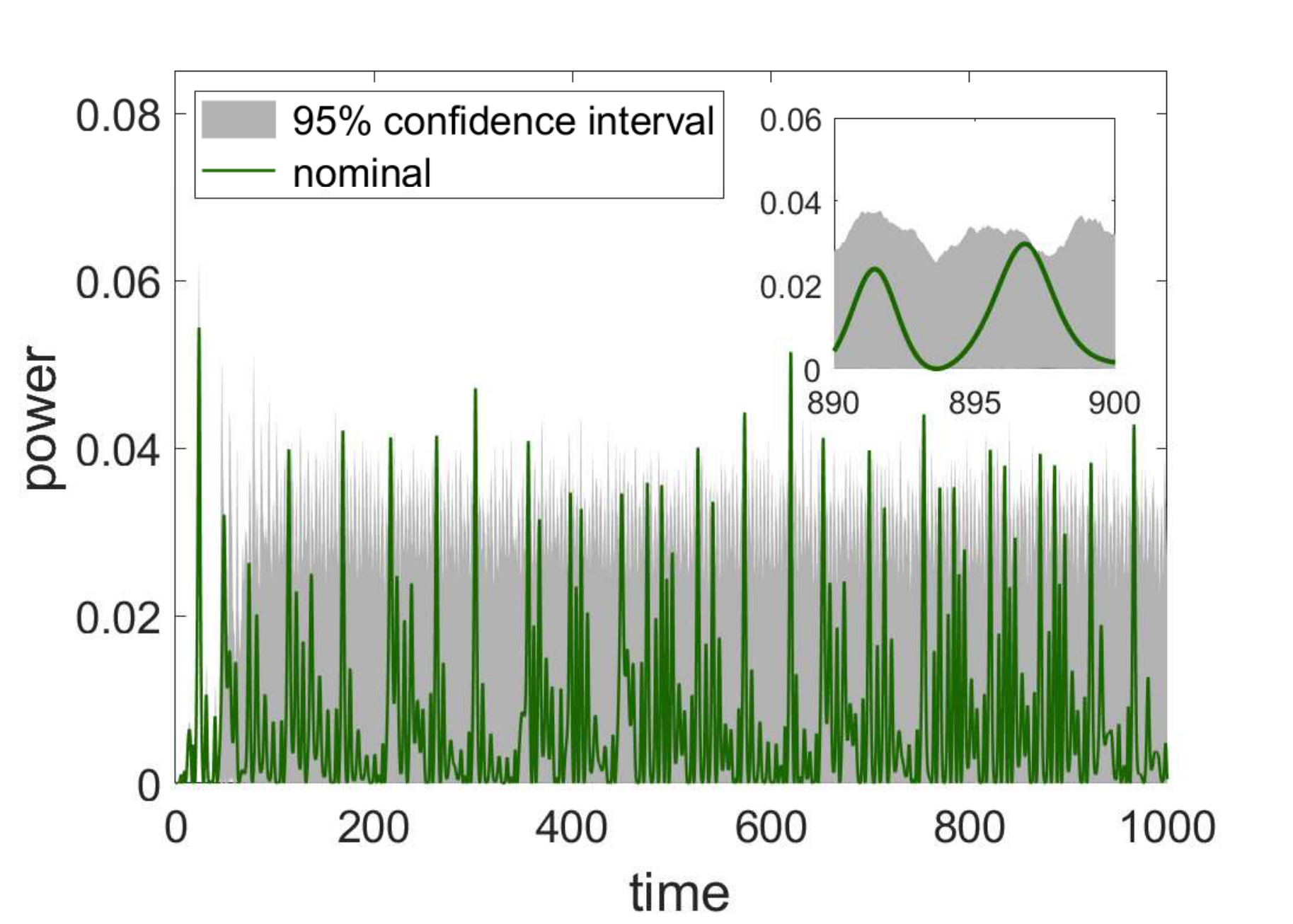}}
    \subfigure[$\mathnormal{f}=0.250$ for $\lambda$]{\includegraphics[width=0.32\textwidth]{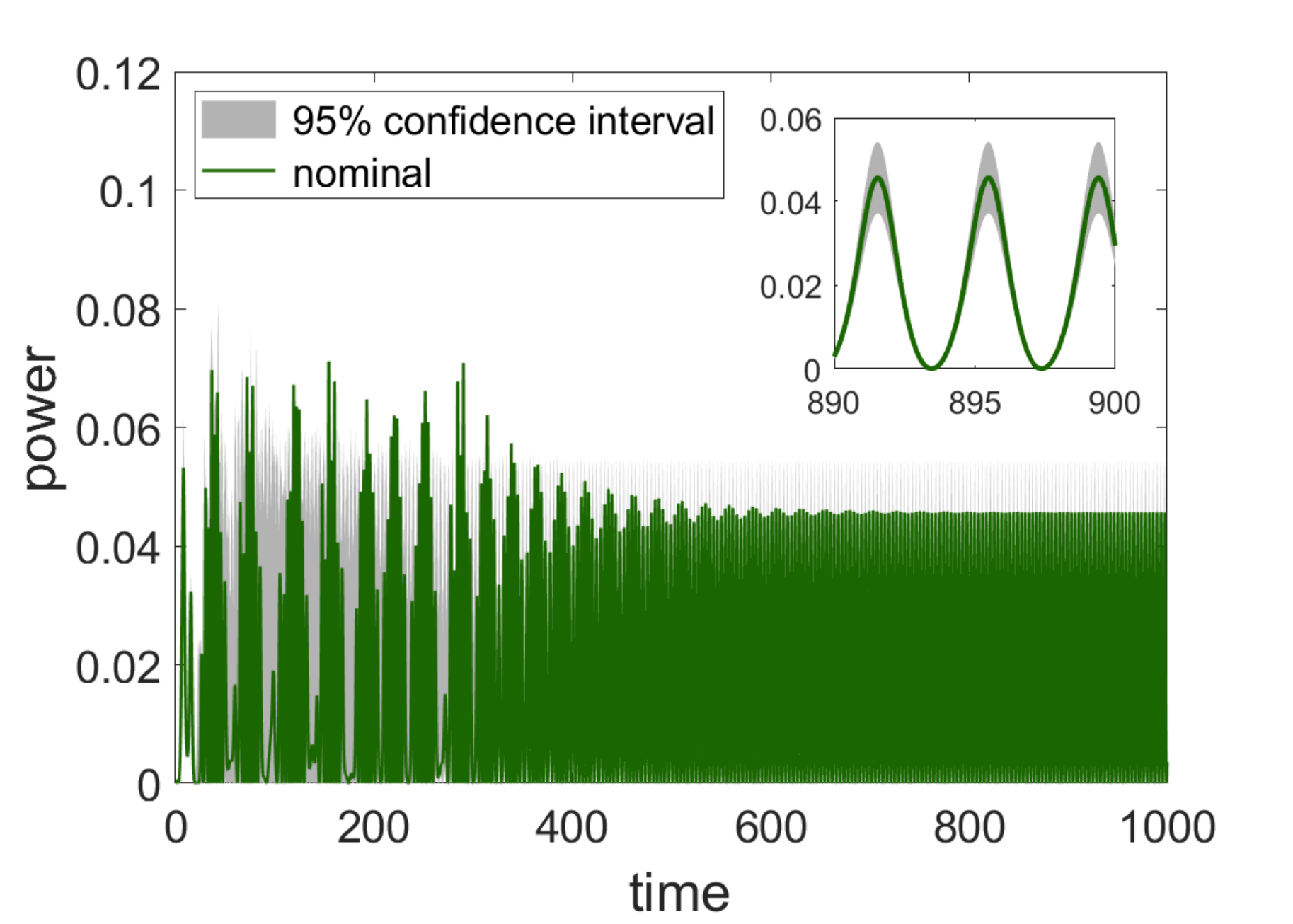}}
    \subfigure[$\mathnormal{f}=0.041$ for $\kappa$]{\includegraphics[width=0.32\textwidth]{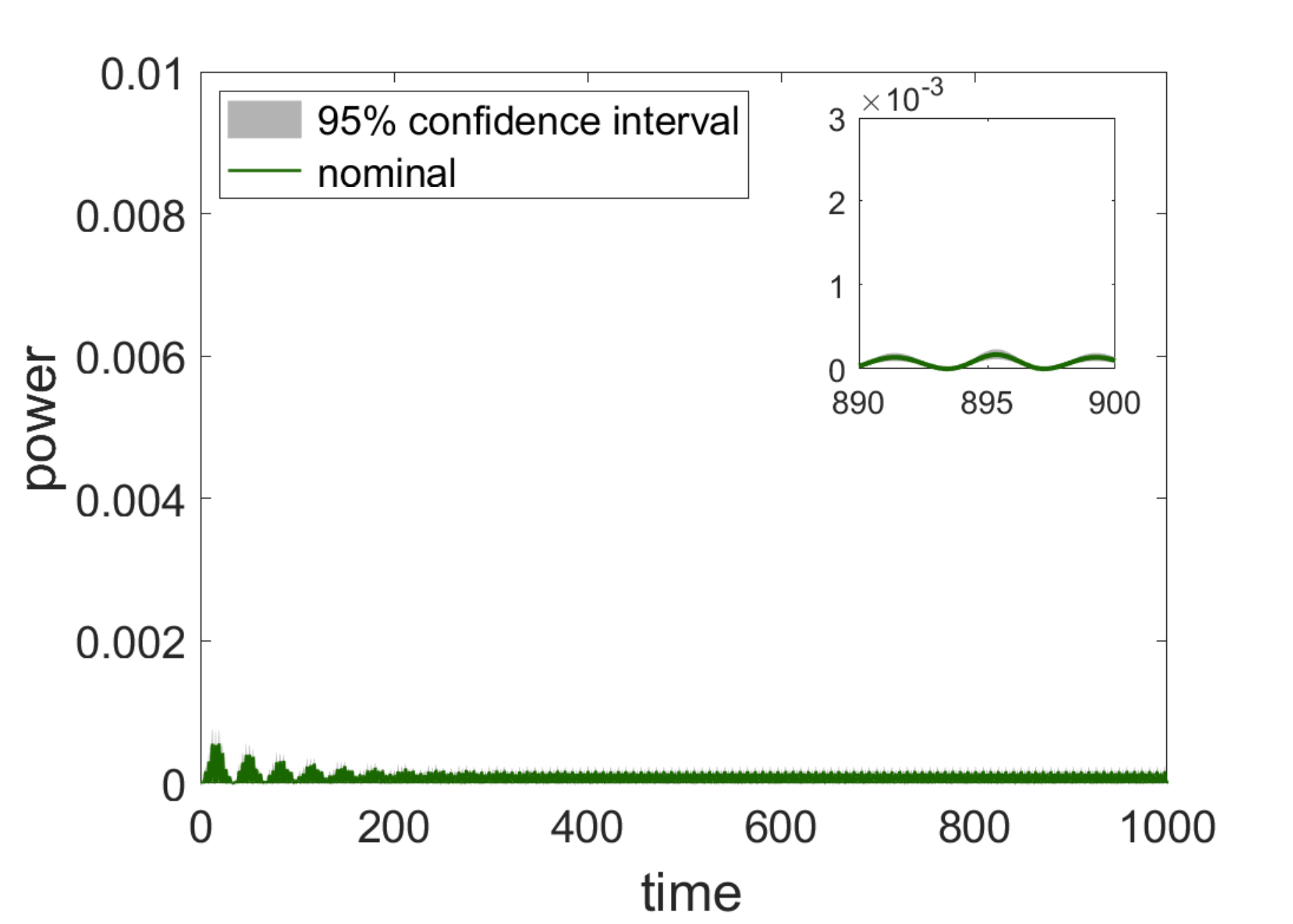}}
    \subfigure[$\mathnormal{f}=0.091$ for $\kappa$]{\includegraphics[width=0.32\textwidth]{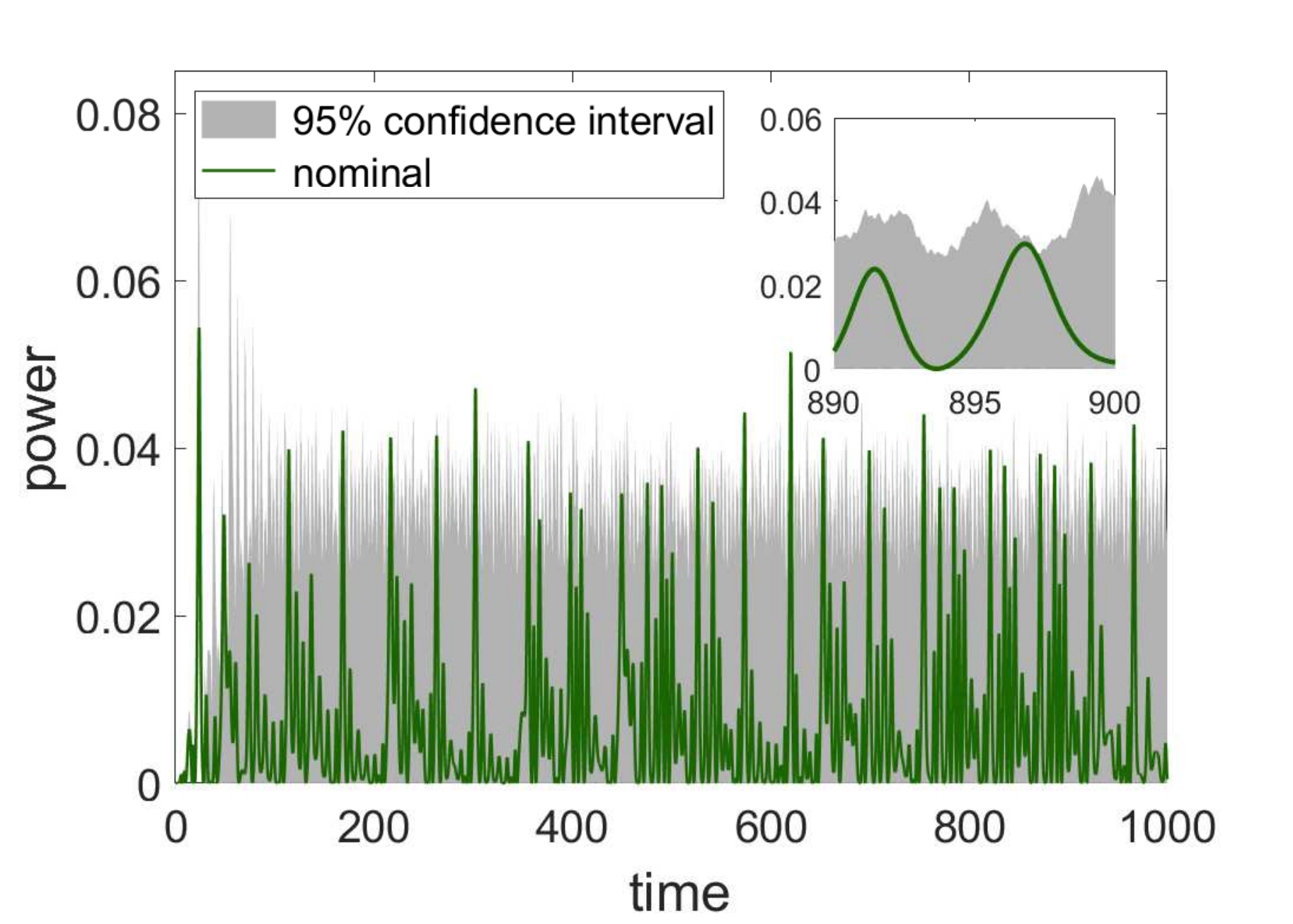}}
    \subfigure[$\mathnormal{f}=0.250$ for $\kappa$]{\includegraphics[width=0.32\textwidth]{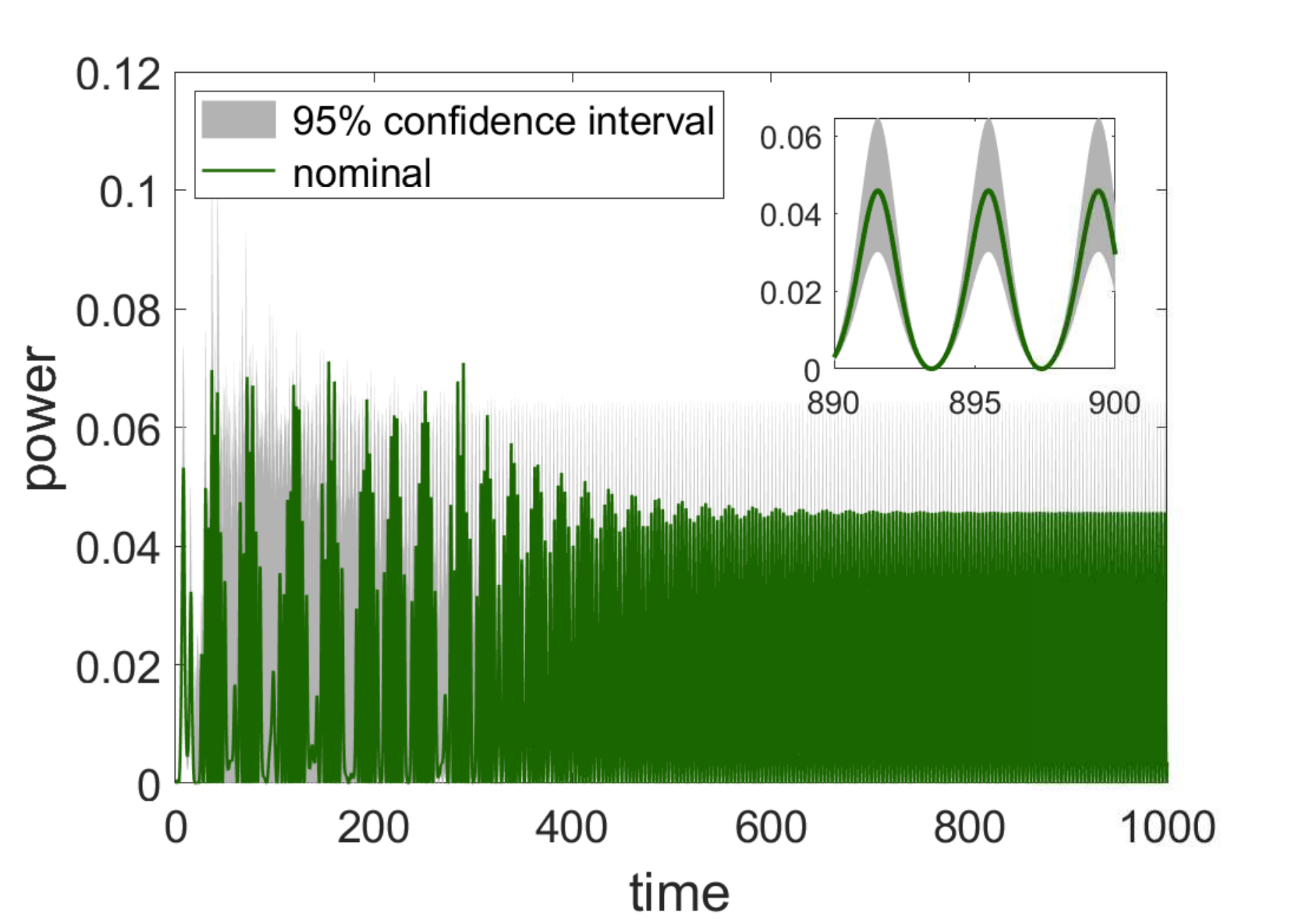}}
    \subfigure[$\mathnormal{f}=0.041$ for $f$]{\includegraphics[width=0.32\textwidth]{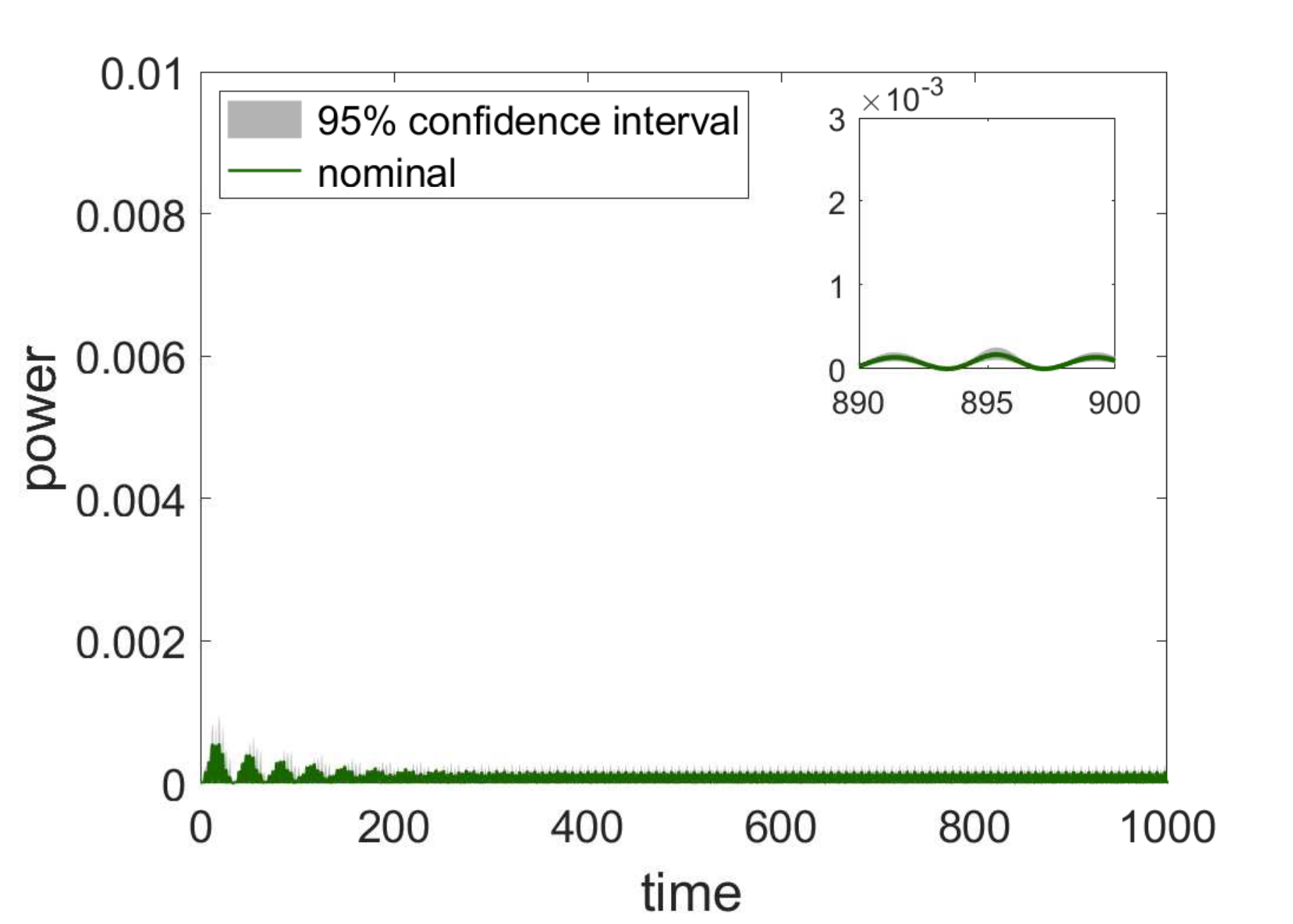}}
    \subfigure[$\mathnormal{f}=0.091$ for $f$]{\includegraphics[width=0.32\textwidth]{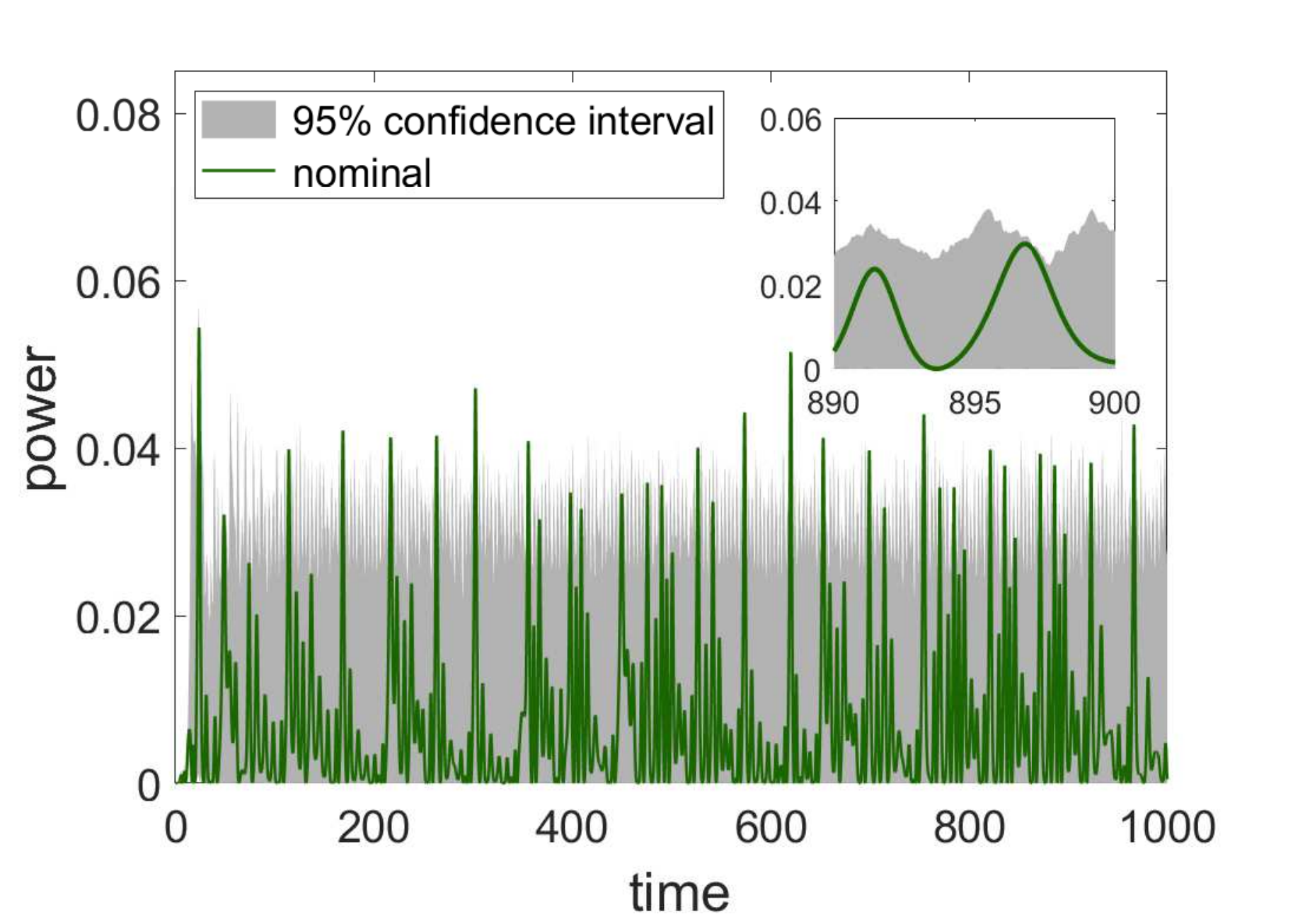}}
    \subfigure[$\mathnormal{f}=0.250$ for $f$]{\includegraphics[width=0.32\textwidth]{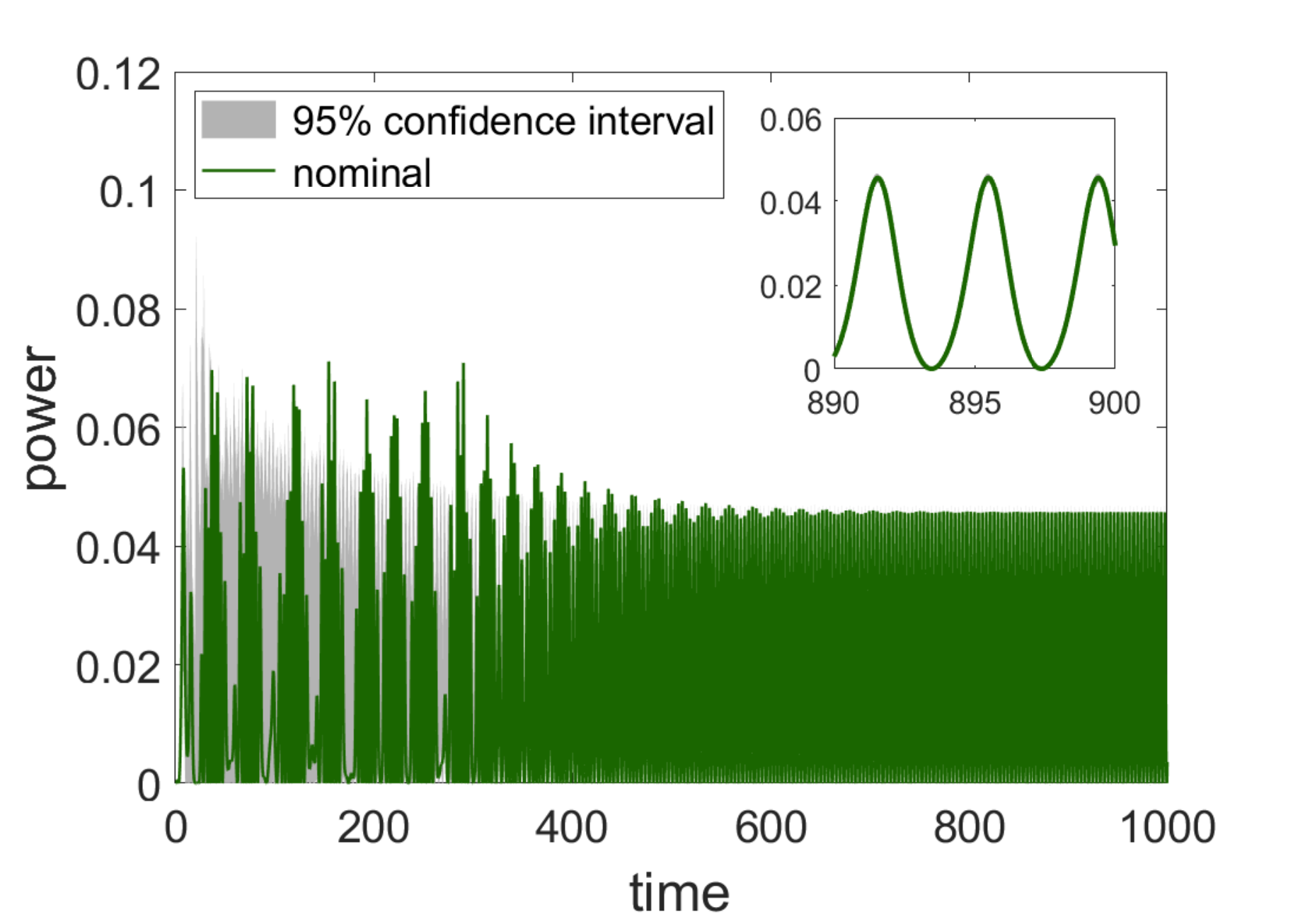}}
    \subfigure[$\mathnormal{f}=0.041$ for $\Omega$]{\includegraphics[width=0.32\textwidth]{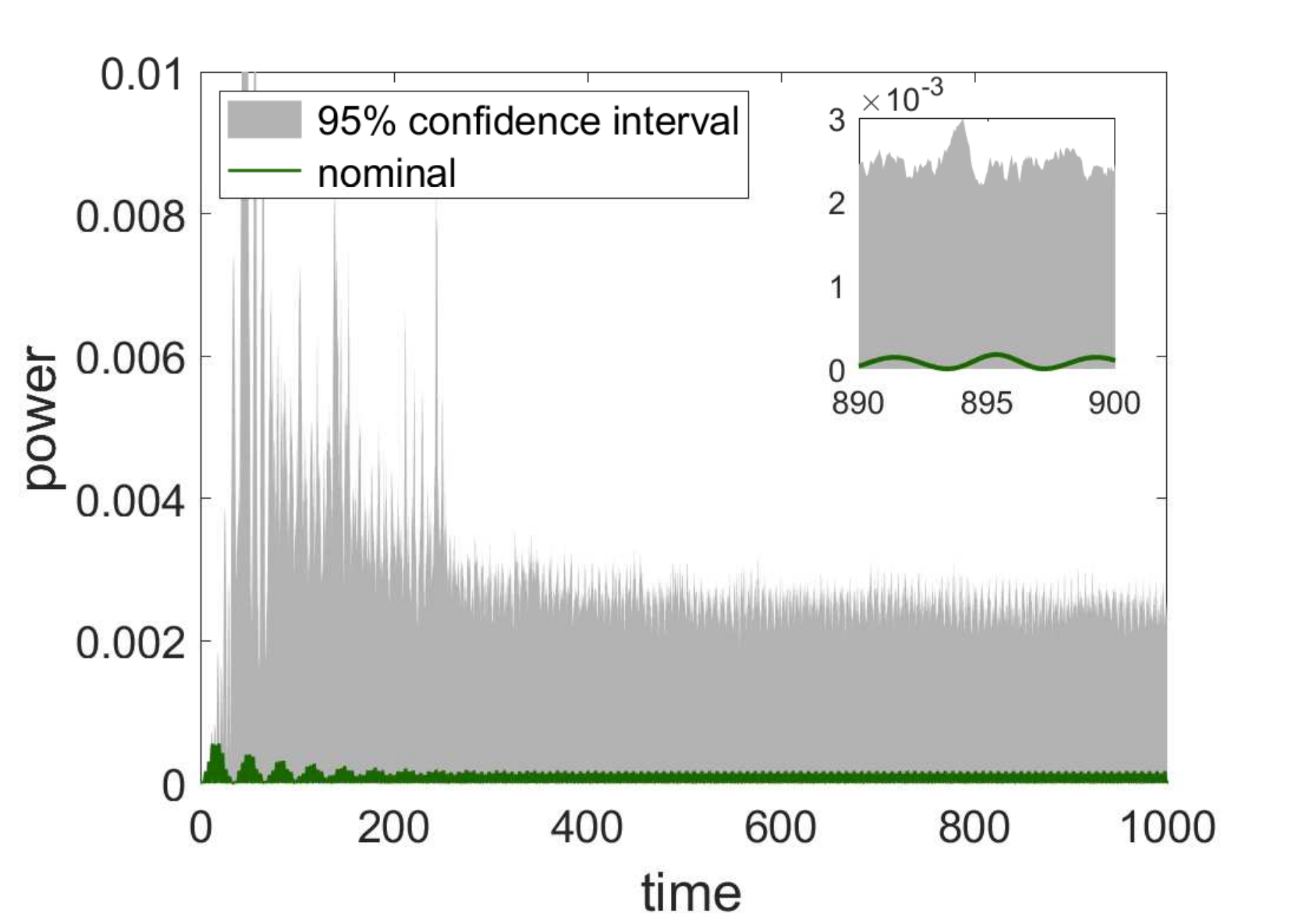}}
    \subfigure[$\mathnormal{f}=0.091$ for $\Omega$]{\includegraphics[width=0.32\textwidth]{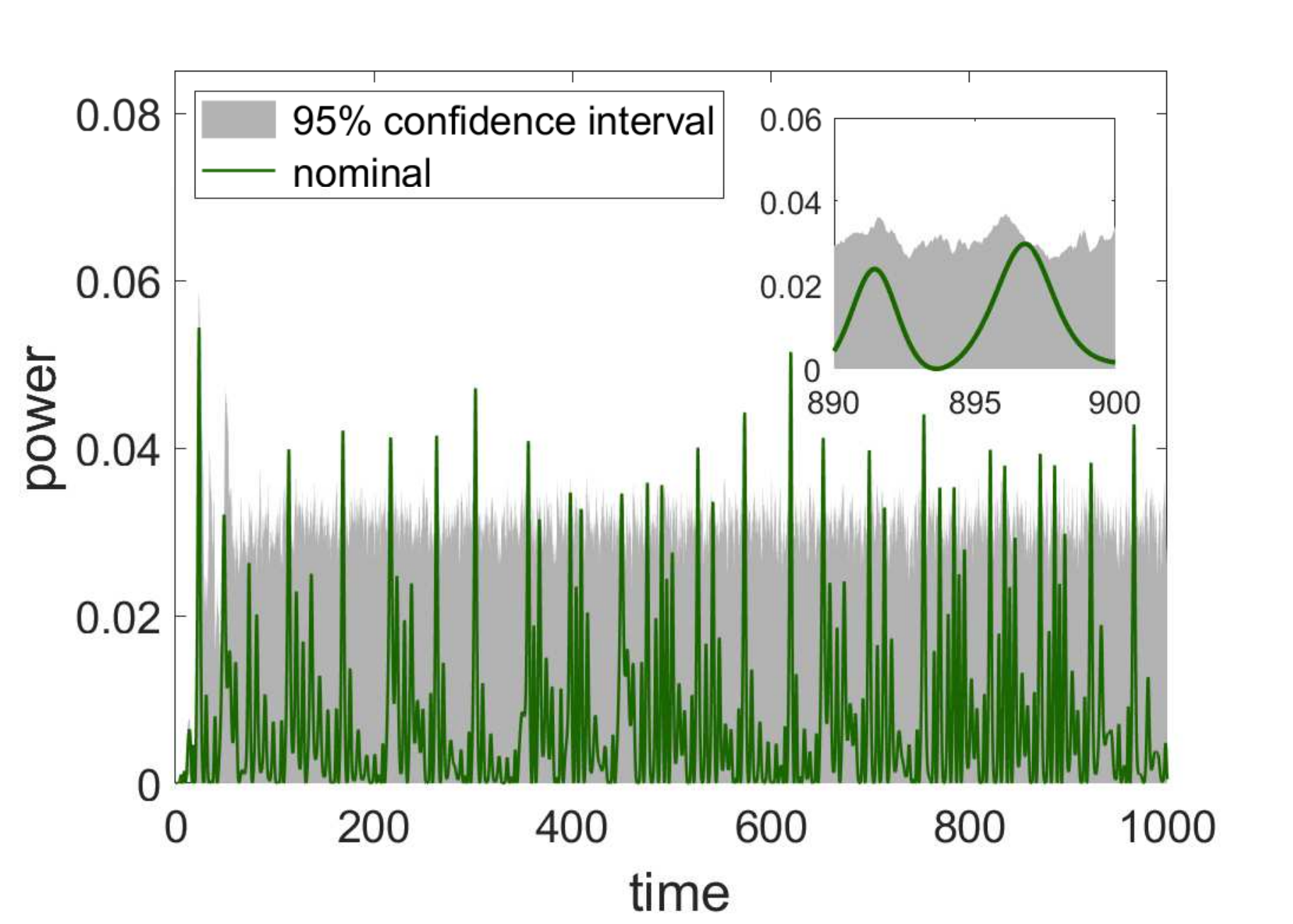}}
    \subfigure[$\mathnormal{f}=0.250$ for $\Omega$]{\includegraphics[width=0.32\textwidth]{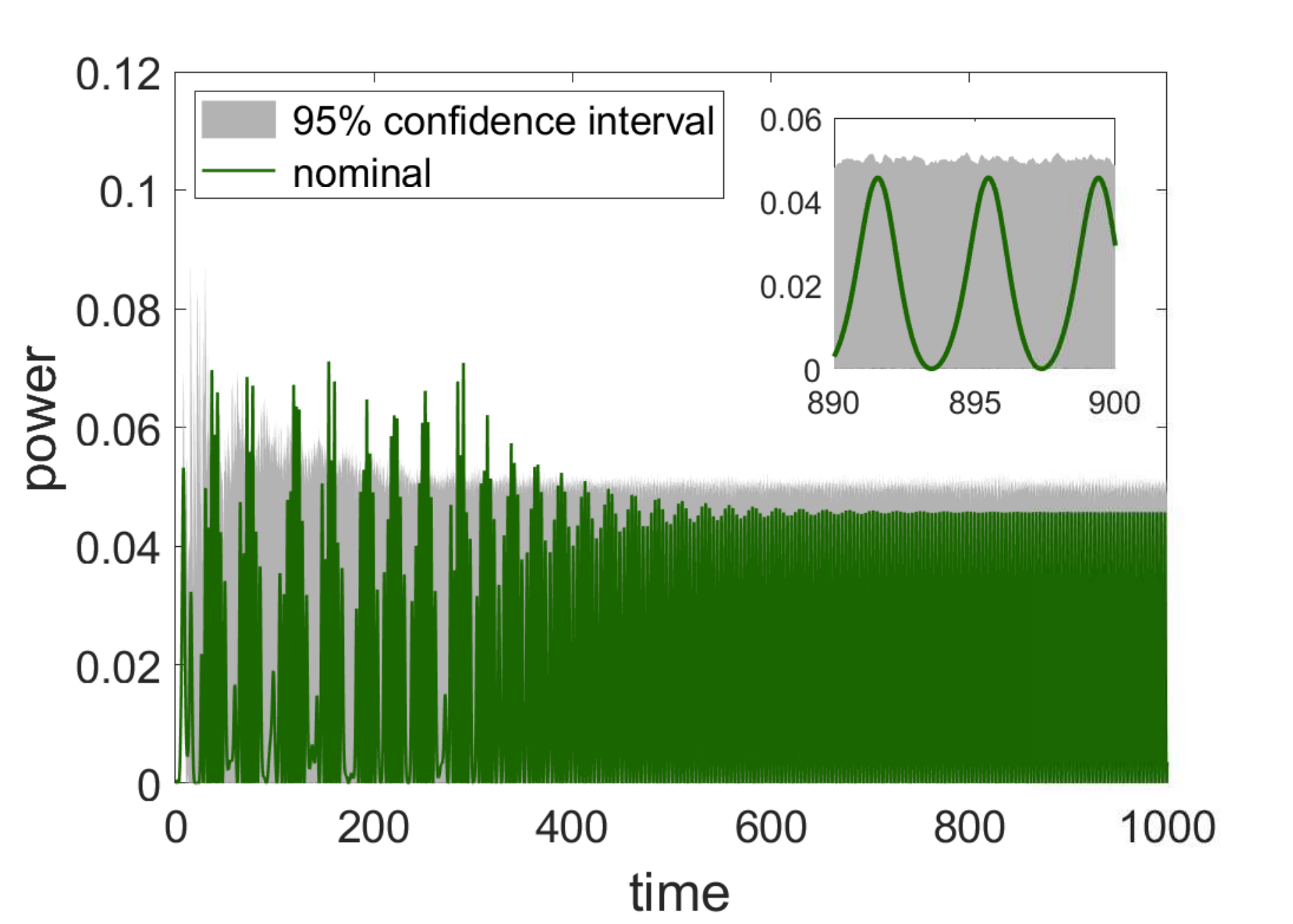}}
    \caption{Propagated uncertainty in the output power time series of the symmetric model with linear piezoelectric coupling is shown under individual parameters: $\lambda$ (first row), $\kappa$ (second row), $\mathnormal{f}$ (third row), and $\Omega$ (fourth row). The columns are divided according to the different motion states of the system: intrawell (left), chaos (middle), and interwell (right).}
    \label{fig:up_BEH}
\end{figure*}

In order to visualize the time response, Fig.~\ref{fig:up_BEH} depicts the uncertainty propagation of the output power over time, with a 95\% confidence interval. The nominal output power time series is also displayed for reference. To demonstrate the diverse behaviors, we present three cases of dynamic behavior: intrawell (left column), chaotic (middle column), and interwell (right column). Each case focuses on the influence of a single random parameter, aiming to isolate its effect on the electrical power. 

For the intrawell motion, it is observed that $\lambda$, $\kappa$, and $\mathnormal{f}$ do not significantly affect the mean power, resulting in a narrow confidence interval. However, for $\Omega$, the interval is substantially wider, and higher values of harvested energy are observed. For the chaotic motion, all parameters alter the confidence interval due to the sensitivity of chaotic behavior to small variability. Finally, for the interwell motion, variations in $\kappa$ and $\Omega$ affect the power. While the effect of $\Omega$ amplifies the envelope in the horizontal direction, it has minimal impact on the amplitude. In contrast, the effect of $\kappa$ increases the amplitude of the envelope, resulting in an increase or decrease in the generated power.

%%%%%%%%%%%%%%%%%%%%%%%%%%%%%%%%%%%%%%
\subsection{Symmetric bistable energy harvester with nonlinear piezoelectric coupling }

\begin{figure*}
    \centering
    \subfigure[$\mathnormal{f}=0.041$]{\includegraphics[width=0.32\textwidth]{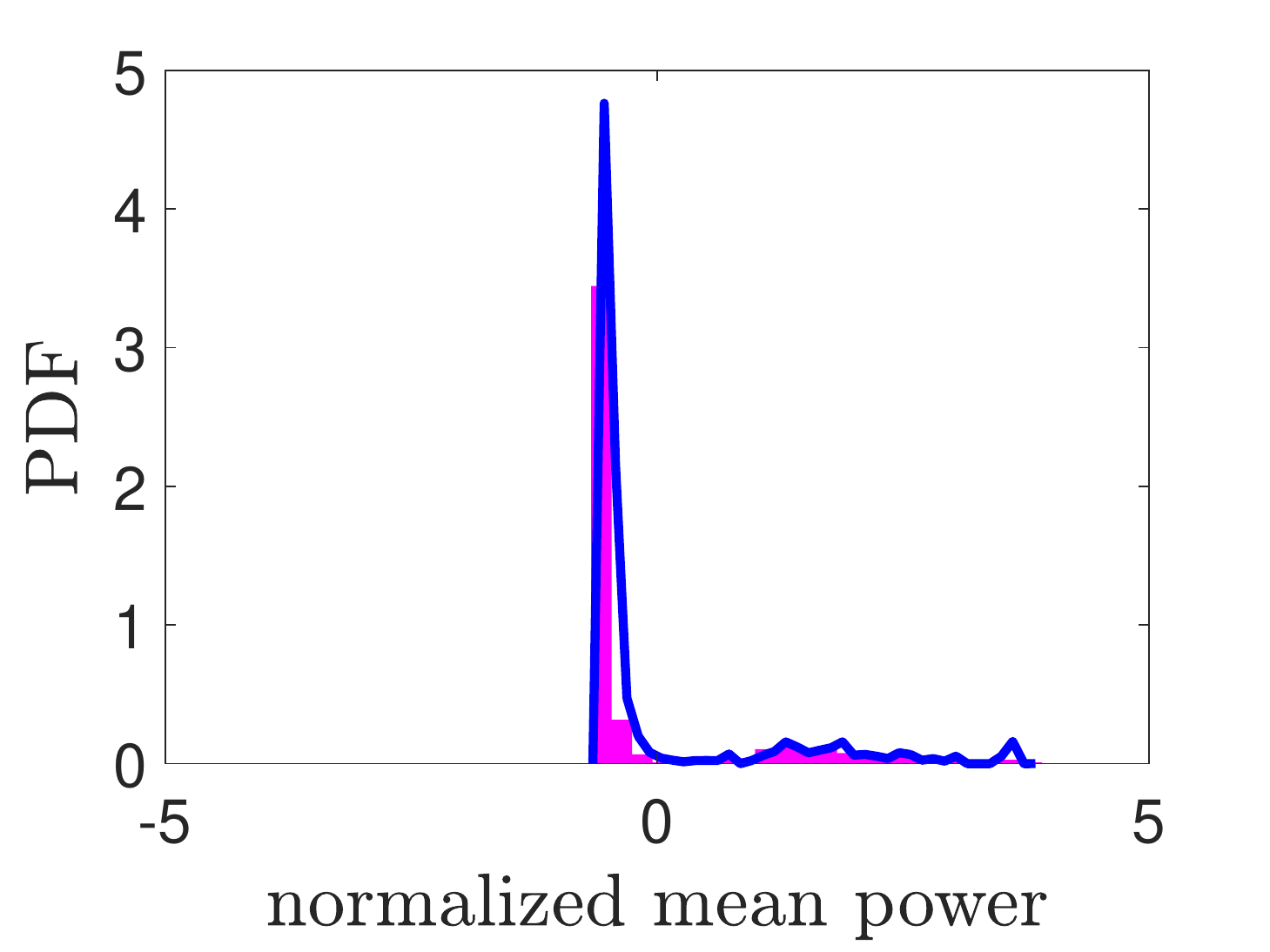}}
    \subfigure[$\mathnormal{f}=0.060$]{\includegraphics[width=0.32\textwidth]{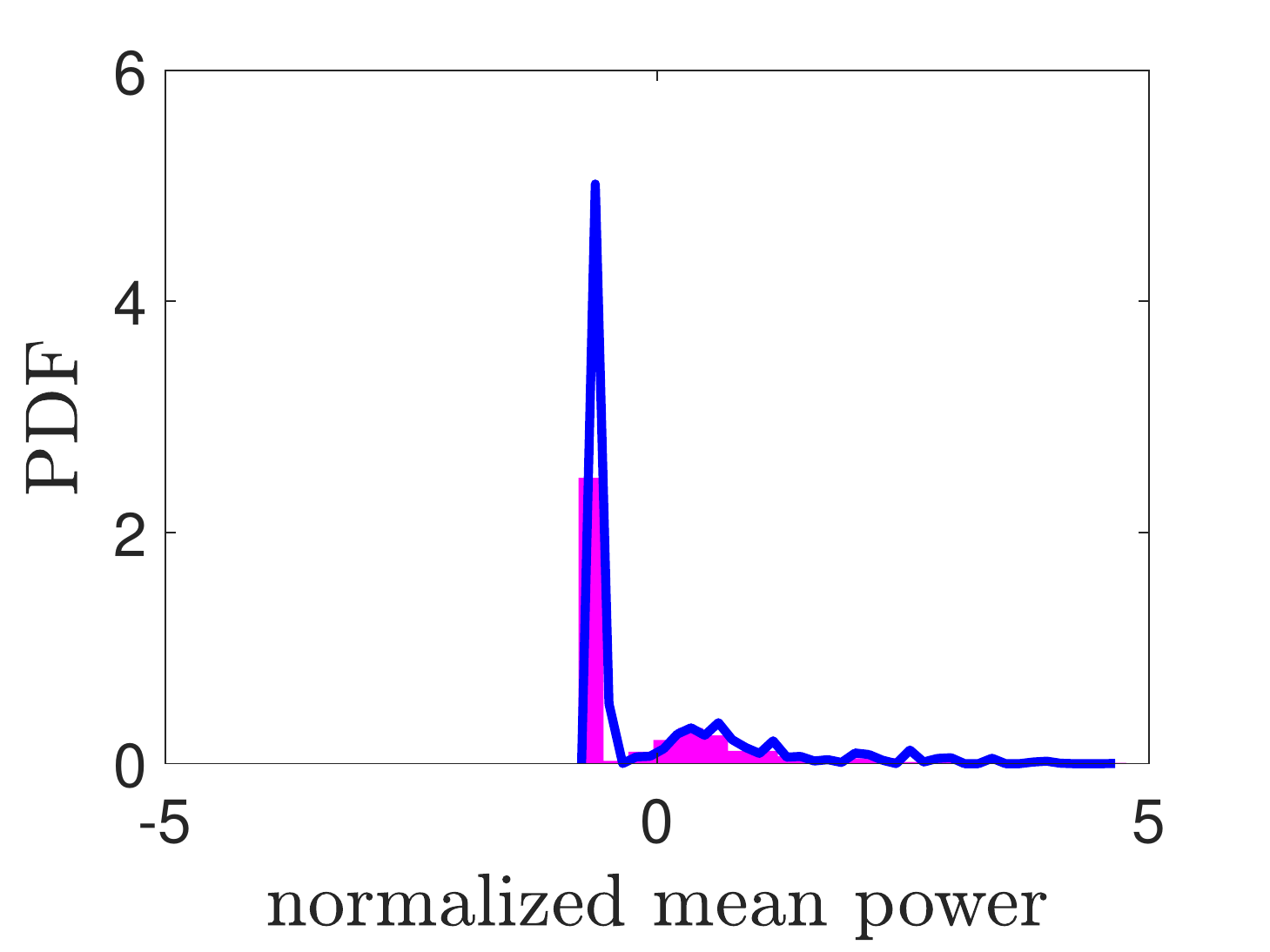}}
    \subfigure[$\mathnormal{f}=0.083$]{\includegraphics[width=0.32\textwidth]{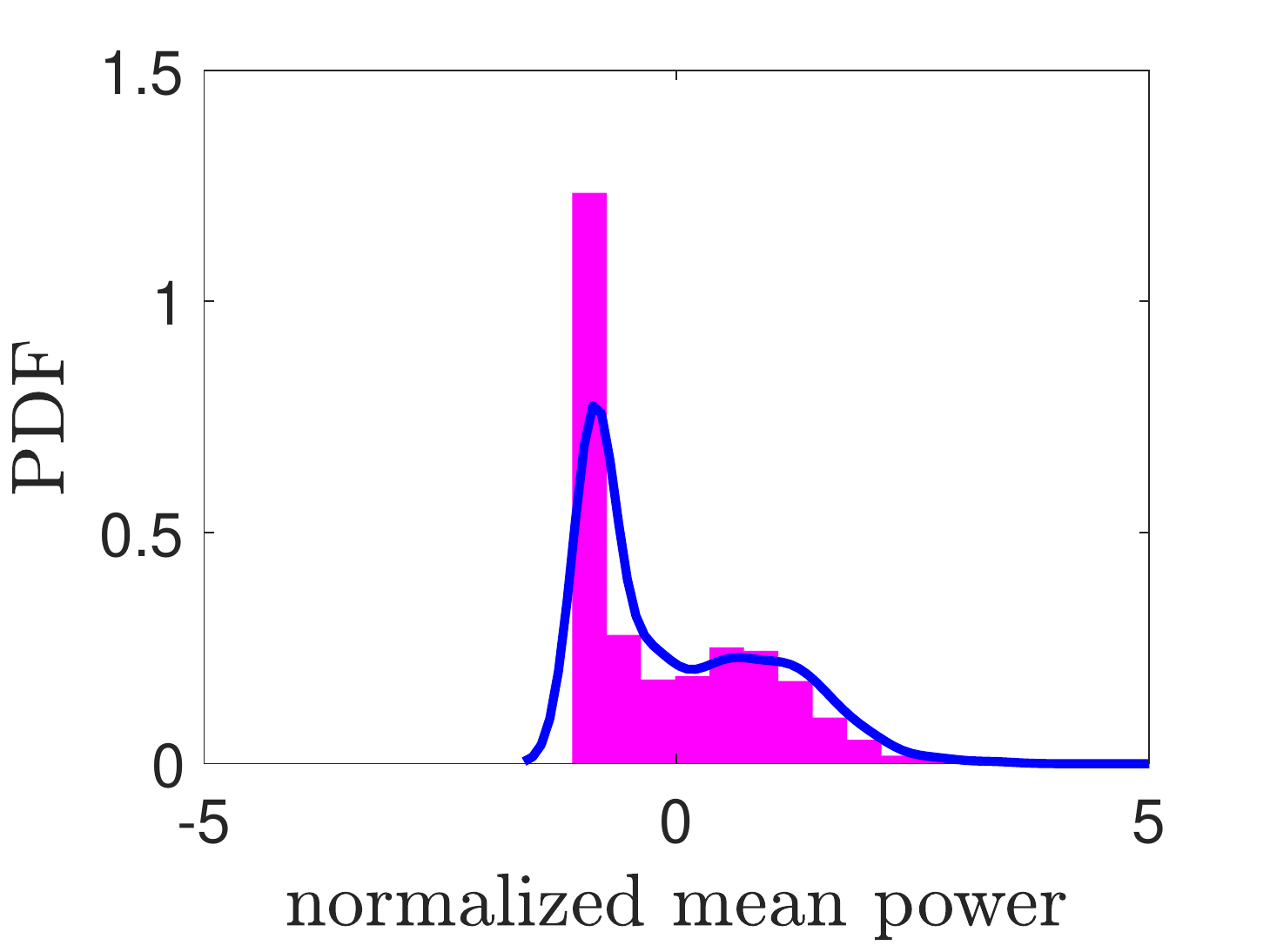}}
    \subfigure[$\mathnormal{f}=0.091$]{\includegraphics[width=0.32\textwidth]{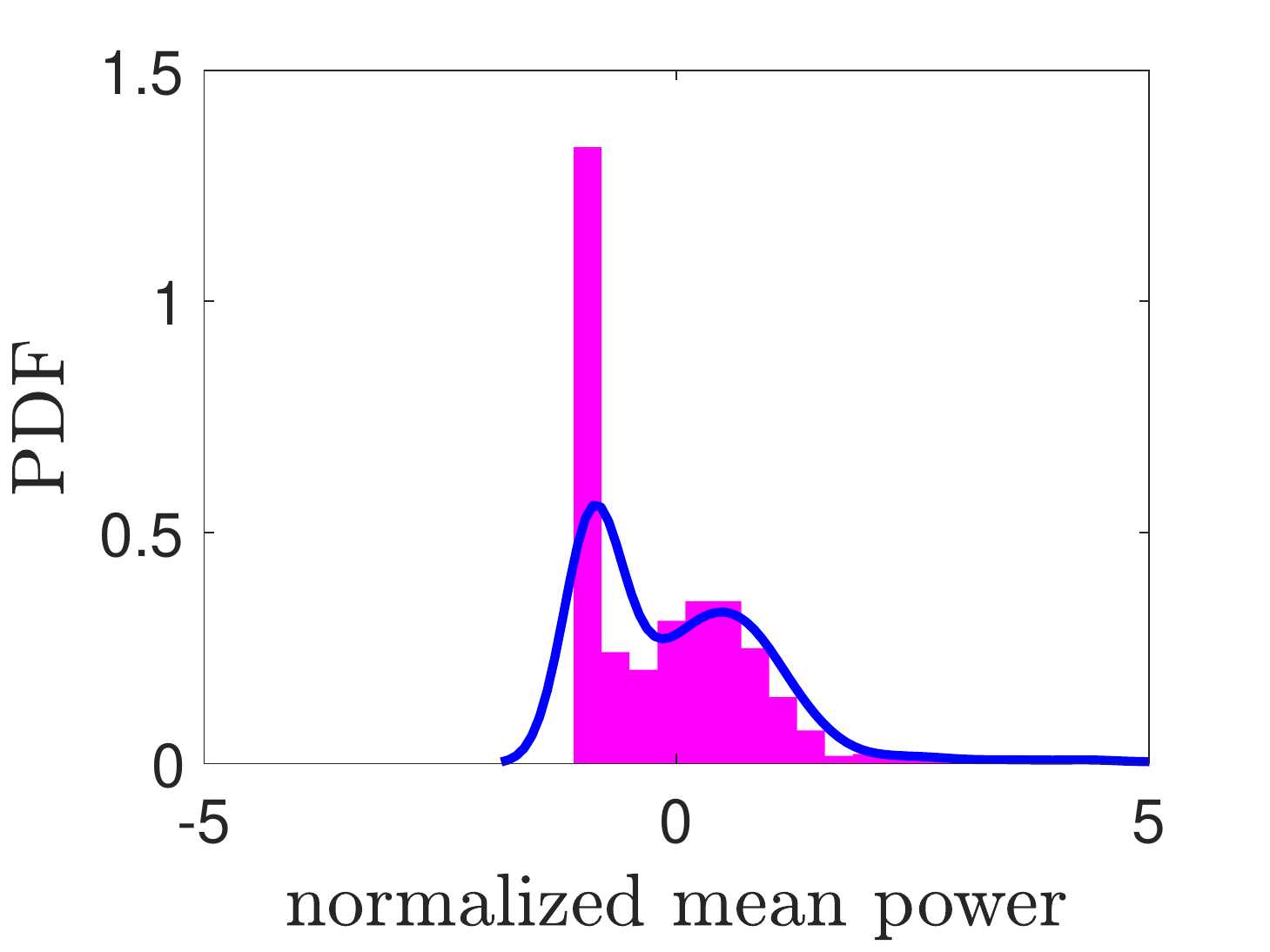}}
    \subfigure[$\mathnormal{f}=0.105$]{\includegraphics[width=0.32\textwidth]{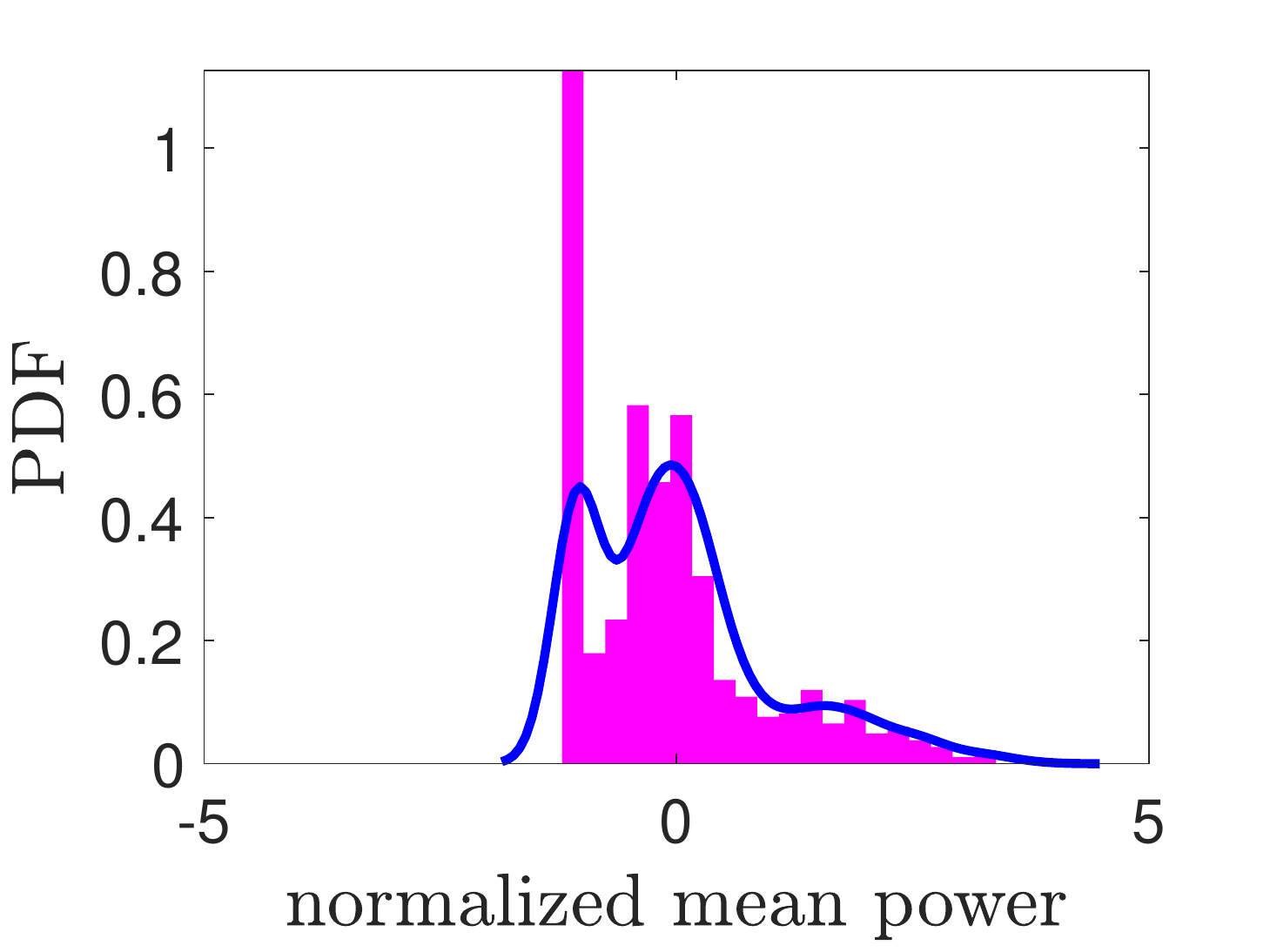}}
    \subfigure[$\mathnormal{f}=0.115$]{\includegraphics[width=0.32\textwidth]{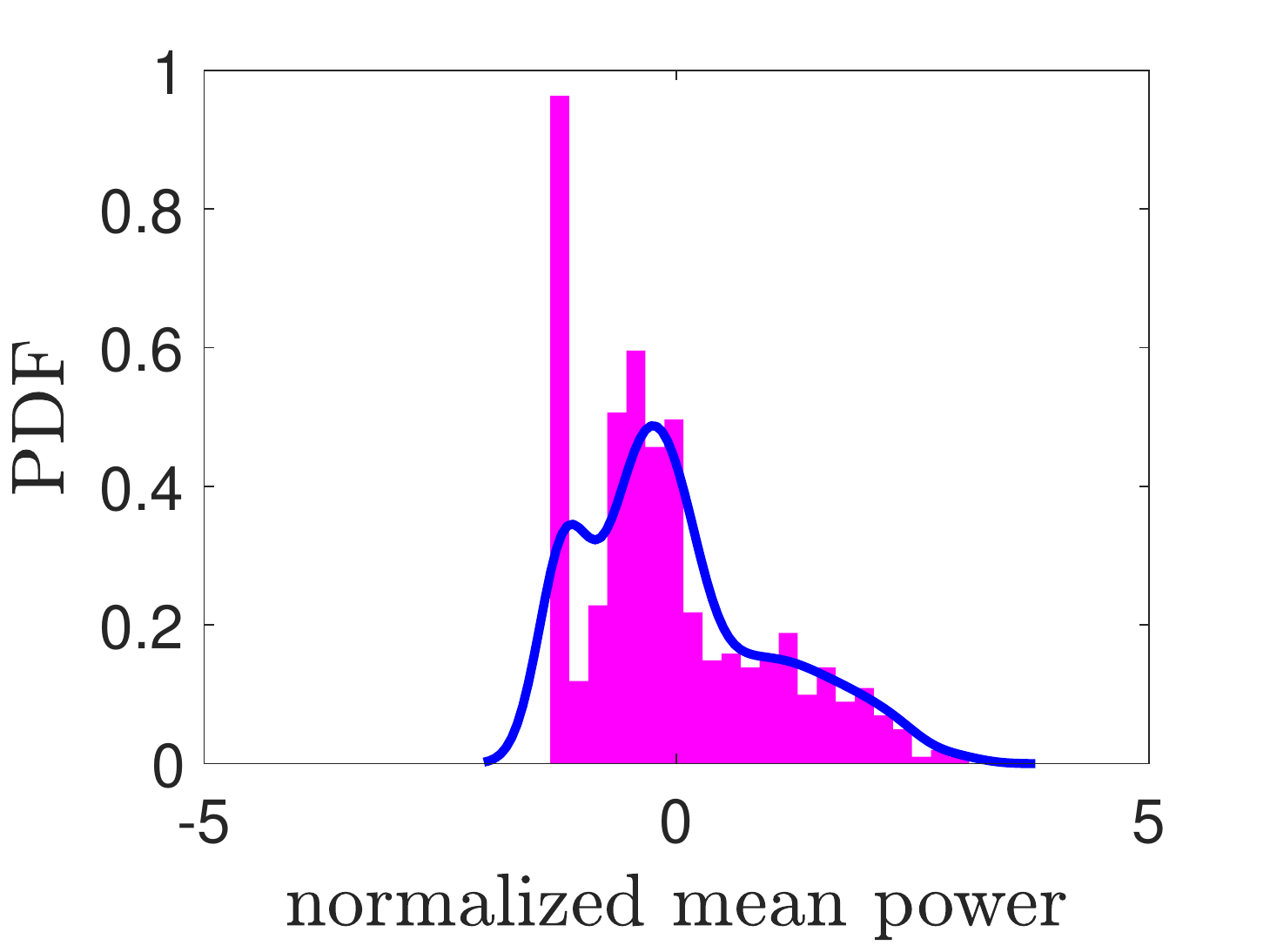}}
    \subfigure[$\mathnormal{f}=0.147$]{\includegraphics[width=0.32\textwidth]{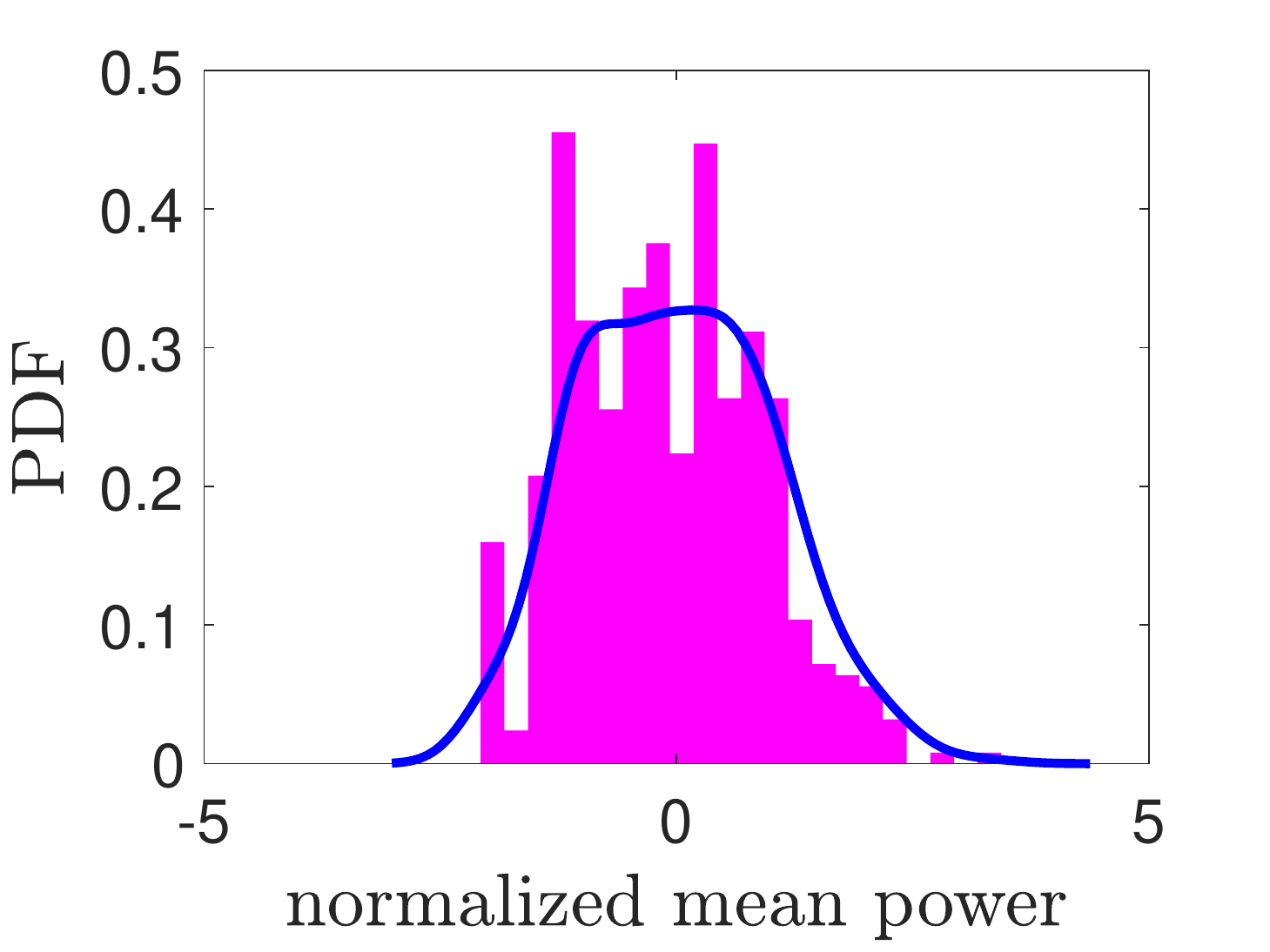}}
    \subfigure[$\mathnormal{f}=0.200$]{\includegraphics[width=0.32\textwidth]{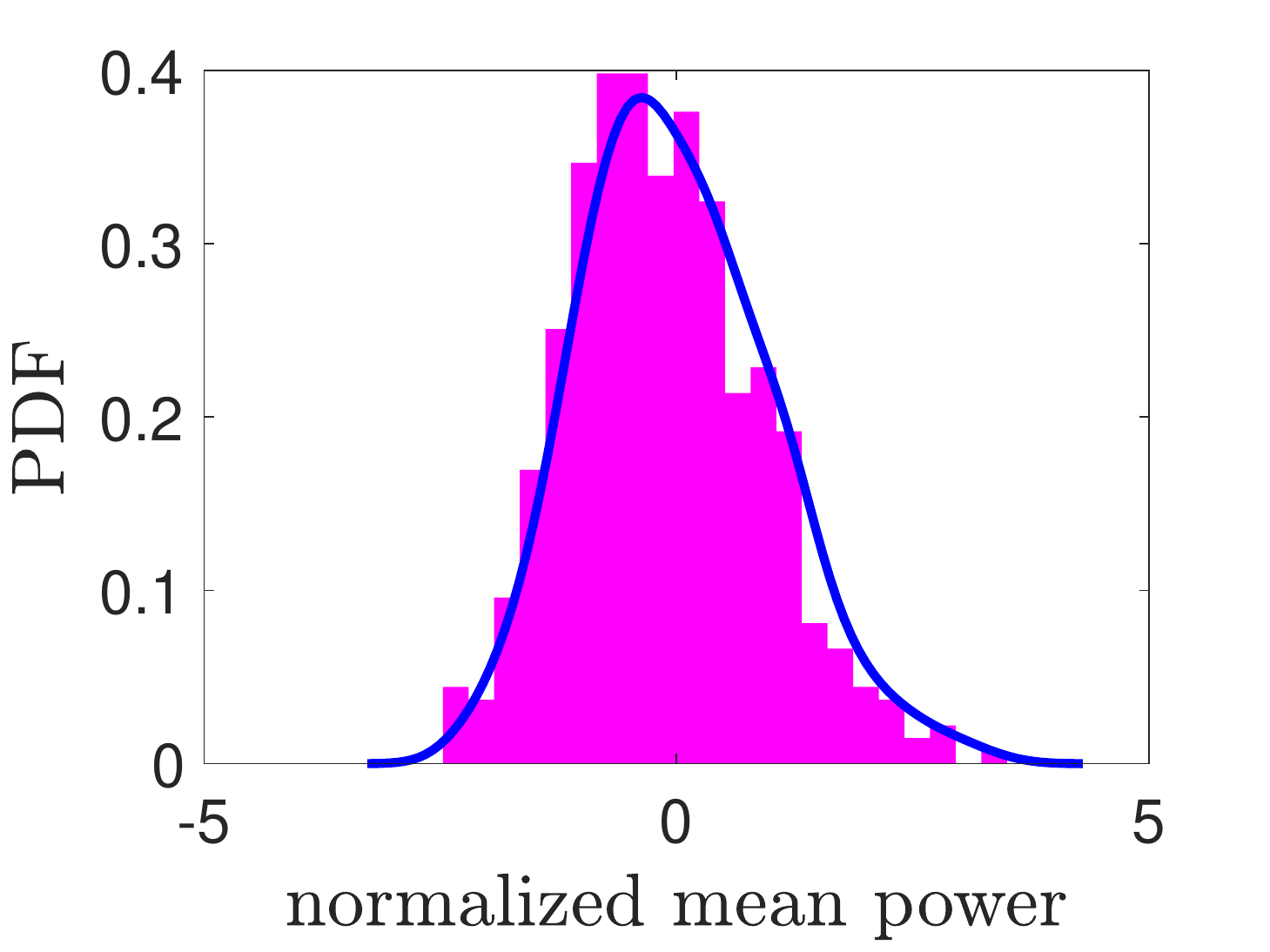}}
    \subfigure[$\mathnormal{f}=0.250$]{\includegraphics[width=0.32\textwidth]{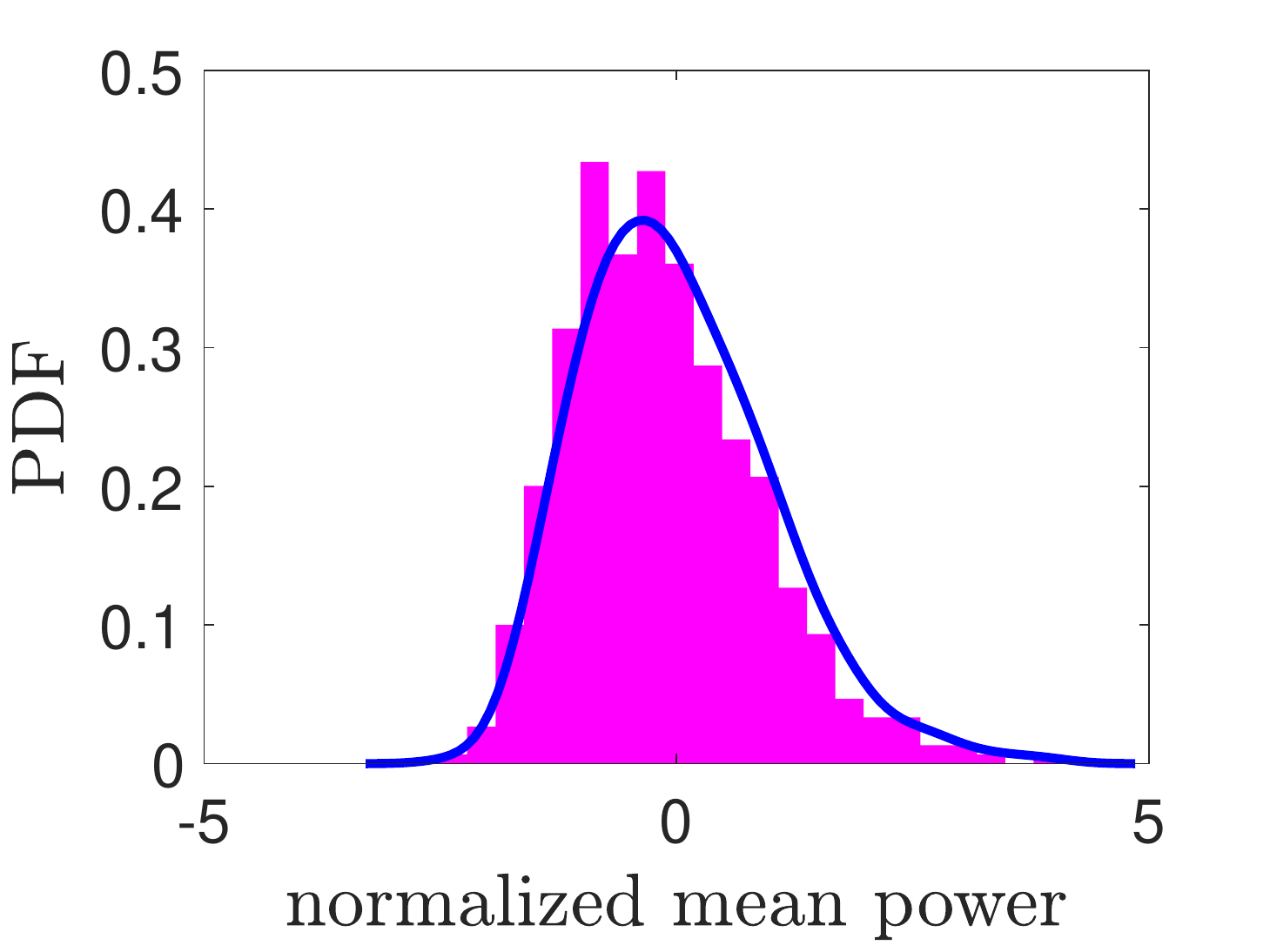}}
    \caption{Probability density function of the normalized mean power for the symmetric energy harvester model with nonlinear piezoelectric coupling under different excitation amplitudes. The kernel density function is represented by the blue line.}
    \label{fig:pdf_pmehn}
\end{figure*}

Figure~\ref{fig:pdf_pmehn} shows the histograms and the PDF of the normalized mean power for a range of excitation amplitudes in a bistable energy harvester with nonlinear coupling ($\delta = 0$, $\phi = 0$ and $\beta \neq 0$). A bimodal distribution is observed for low values of $\mathnormal{f}$ ($<0.115$). As the amplitude of excitation increases, the second peak becomes more prominent, resulting in a unimodal distribution. The observed behavior is similar to that of the previously analyzed bistable model.

\begin{figure*}
    \centering
    \includegraphics[width=1\textwidth]{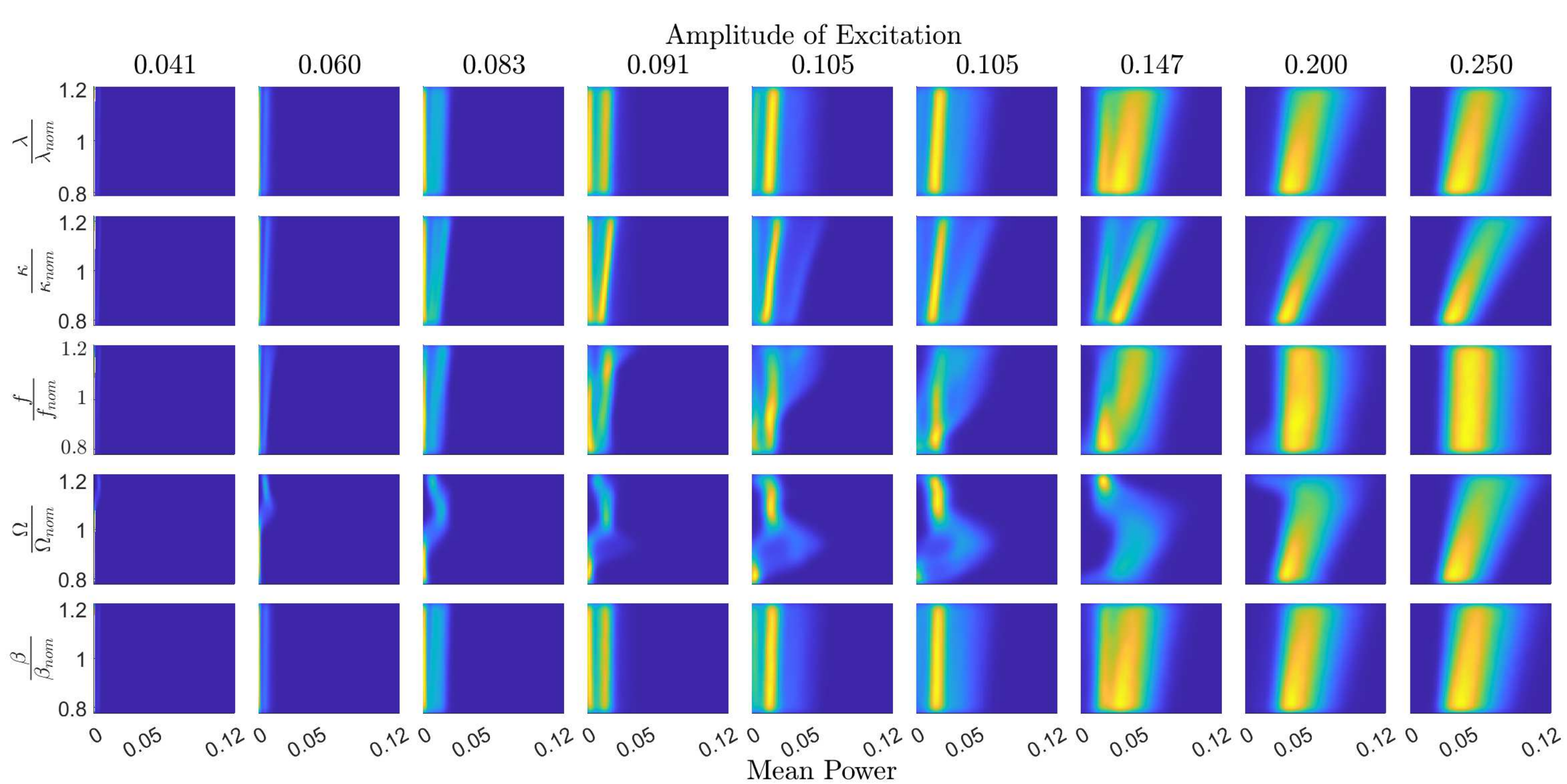}
    \caption{Joint-CDF of mean power conditioned on each parameter of interest ($\lambda$, $\kappa$, $\mathnormal{f}$, $\Omega$) under different values of excitation amplitude for the symmetric model with nonlinear piezoelectric coupling.}
    \label{fig:joint_BEHn}
\end{figure*}

In Figure~\ref{fig:joint_BEHn}, we present the joint-CDF of the mean power conditioned on each parameter of interest under various nominal excitation conditions, for the bistable energy harvester with nonlinear coupling. The impact of the parameters on the mean power remains similar to the case of linear coupling across the explored region. In particular, for low and mid-range values of the excitation amplitude ($\mathnormal{f}<0.200$), the additional term ($\beta$) does not introduce any significant changes. However, for high excitation amplitudes ($\mathnormal{f} =$ 0.200 and 0.250), a slight positive correlation is observed, where increasing $\beta$ leads to a higher mean power.

%Figure~\ref{fig:prob_BEHn} presents the conditional probability of increasing the mean power by 50\% of the nominal power given the parameter of interest also increased by 10\%. For intrawell motion, $\Omega$ is essential for increasing power. For chaotic motion, the increase of $\mathnormal{f}$ is fundamental to improving the mean harvested power, while the increase of $\Omega$ is not attractive, similar to the linear coupling case. Finally, for interwell motion, $\kappa$ is the most critical parameter for increasing energy generation. The probability of increasing power given a higher $\beta$ follows the trend of the other parameters, presenting little relevance.

\begin{figure*}
    \centering
    \includegraphics[width=0.7\textwidth]{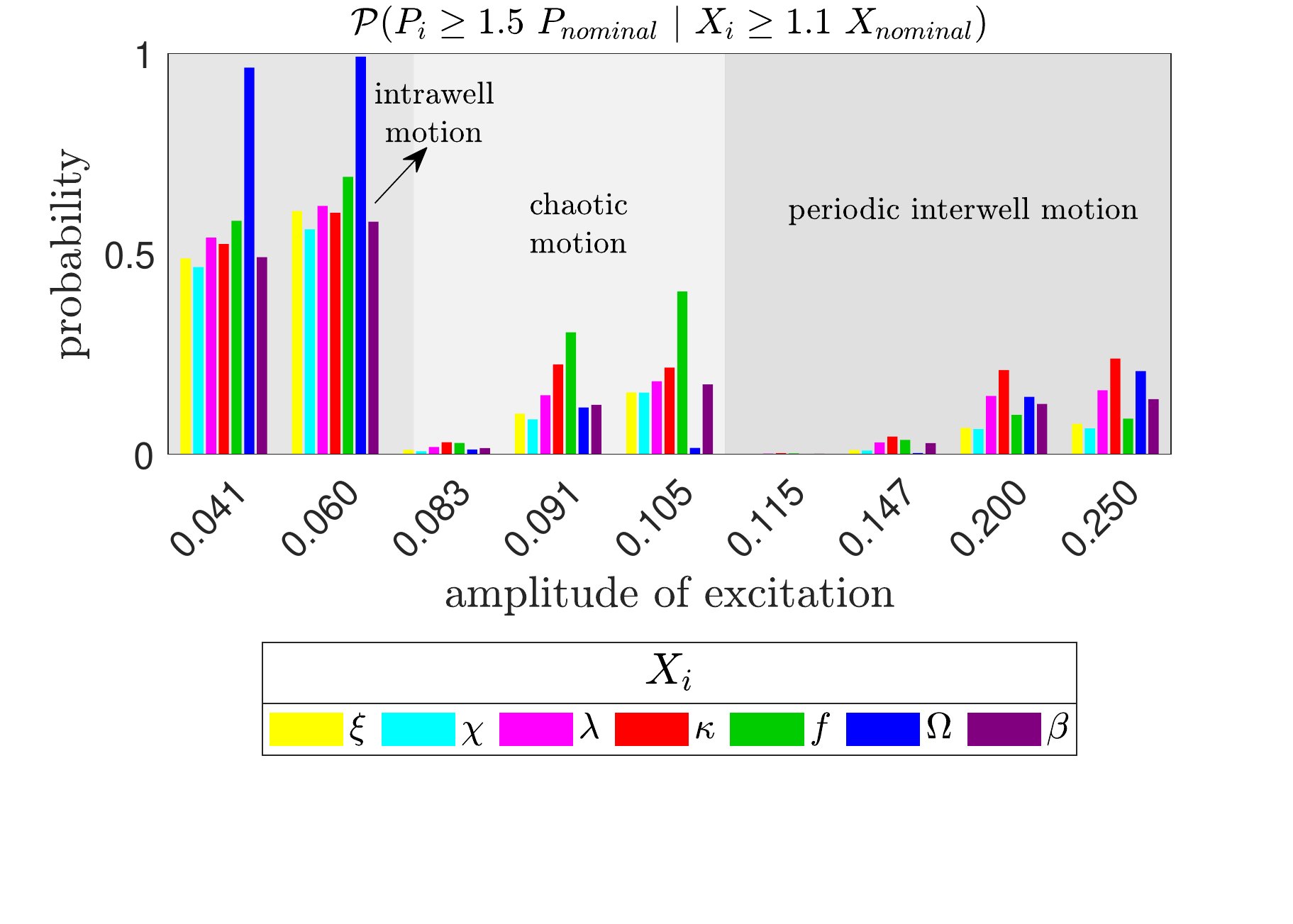}
    \vspace{-1cm}
    \caption{Probability of increasing the nominal mean power by 50\% as parameter $X_\mathnormal{i}$ is increased by 10\%, plotted against the excitation amplitude, for the symmetric model with nonlinear piezoelectric coupling.}
    \label{fig:prob_BEHn}
\end{figure*}

In Figure~\ref{fig:prob_BEHn}, we analyze the conditional probability of increasing the mean power by 50\% of the nominal power, given that one of the parameters of interest increased by 10\%. The results show that the effects of the parameters remain similar to the case of linear coupling in the checked region, with some nuances.

For intrawell motion, the increase of $\Omega$ is essential for improving power generation about the probability of 100\%, while the other parameters ($\lambda$, $\kappa$, and $\mathnormal{f}$) do not have a significant highlight. In chaotic motion, the increase of excitation amplitude ($\mathnormal{f}$) is fundamental to improving the mean harvested power, while the increase of $\Omega$ is not attractive, consistent with the linear coupling case.
In interwell motion, the most critical parameter for increasing energy generation is $\kappa$, with a 20\% chance of improving the harvesting process. On the other hand, the probability of increasing power given a higher $\beta$ follows the trend of the other parameters and presents little relevance. Overall, these results provide insights into the key parameters that affect the energy harvesting process and can guide the design and optimization of bistable energy harvesters.

\begin{figure*}
    \centering
    \subfigure[$\mathnormal{f}=0.041$ for $\lambda$]{\includegraphics[width=0.32\textwidth]{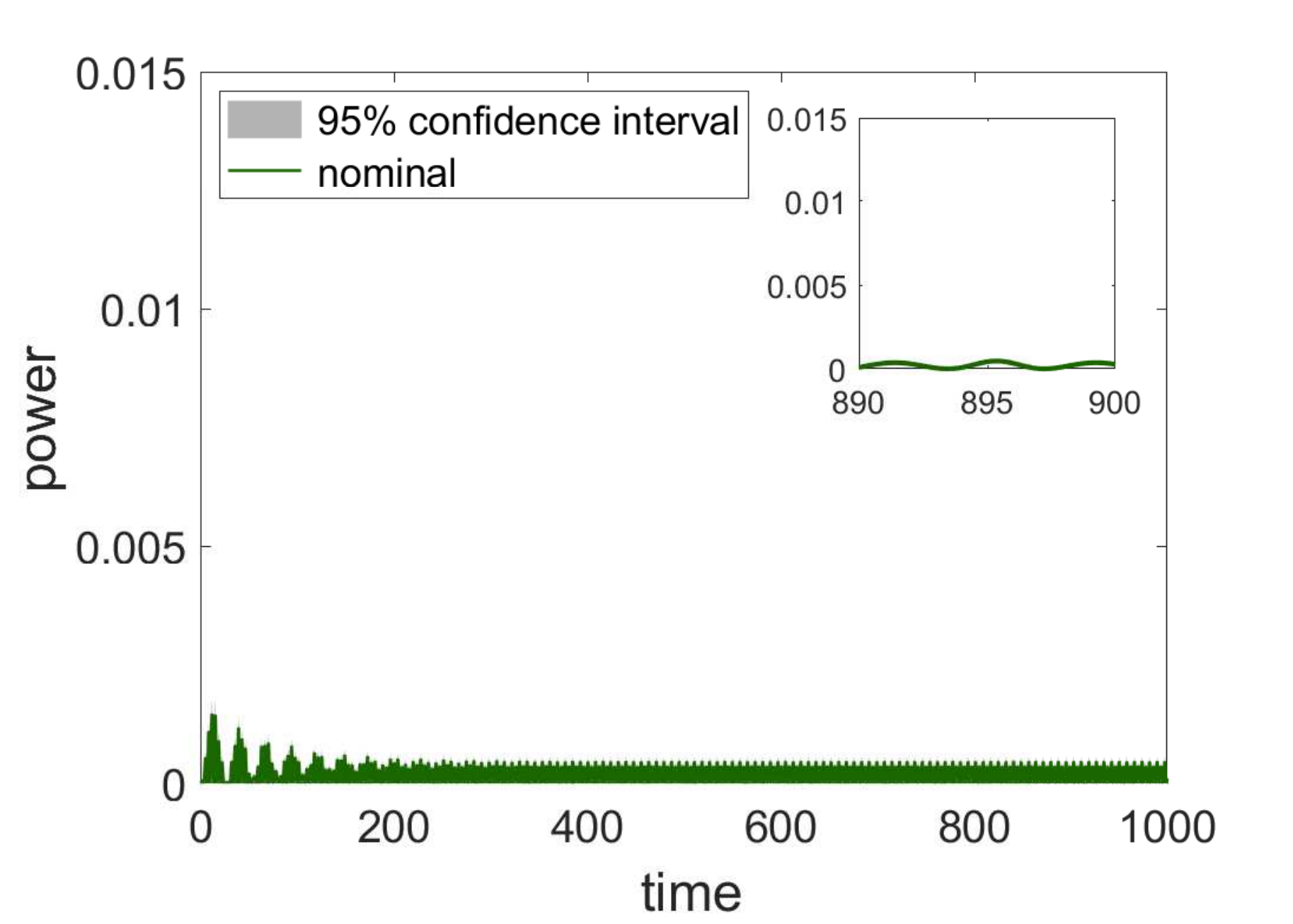}}
    \subfigure[$\mathnormal{f}=0.091$ for $\lambda$]{\includegraphics[width=0.32\textwidth]{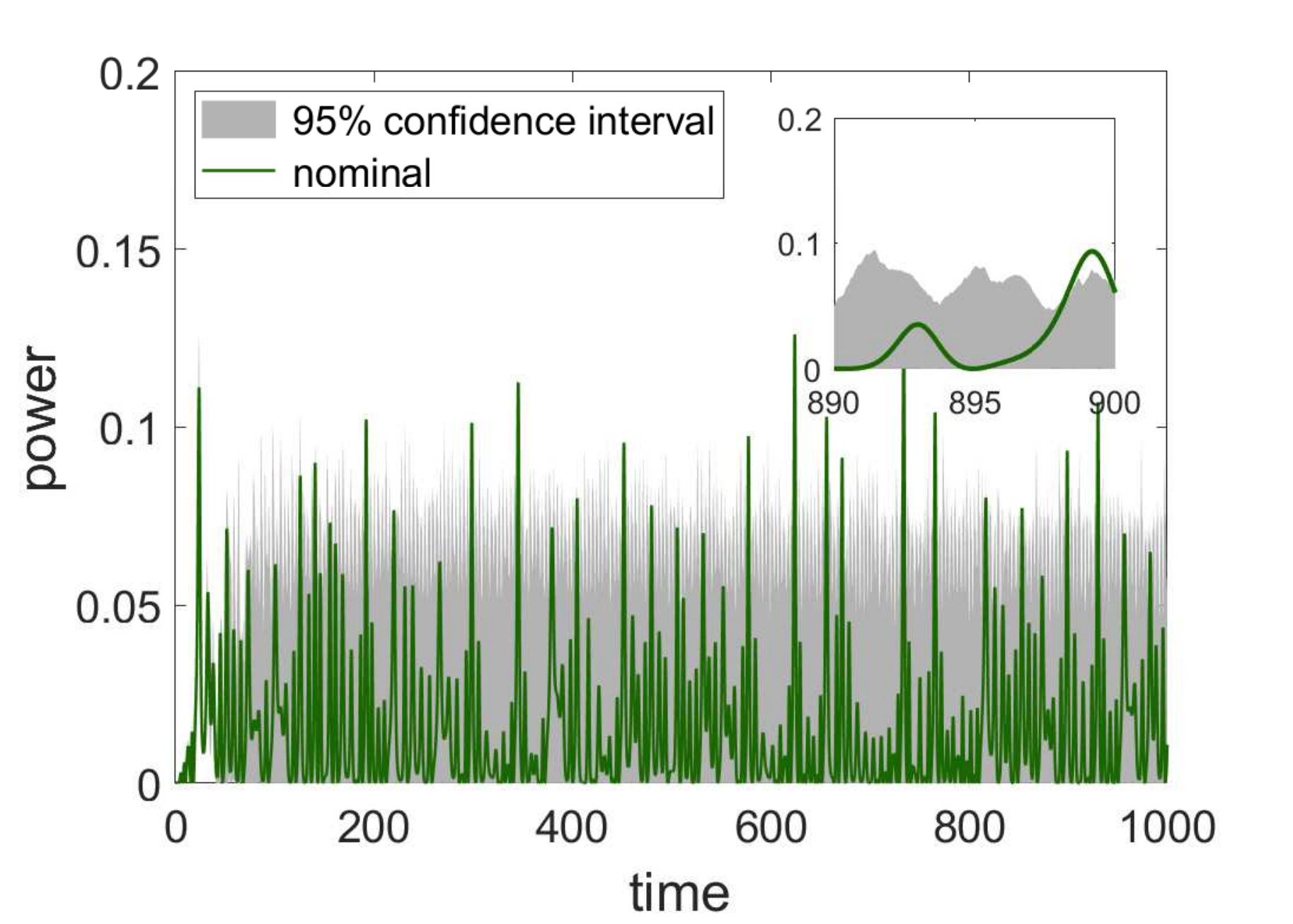}}
    \subfigure[$\mathnormal{f}=0.250$ for $\lambda$]{\includegraphics[width=0.32\textwidth]{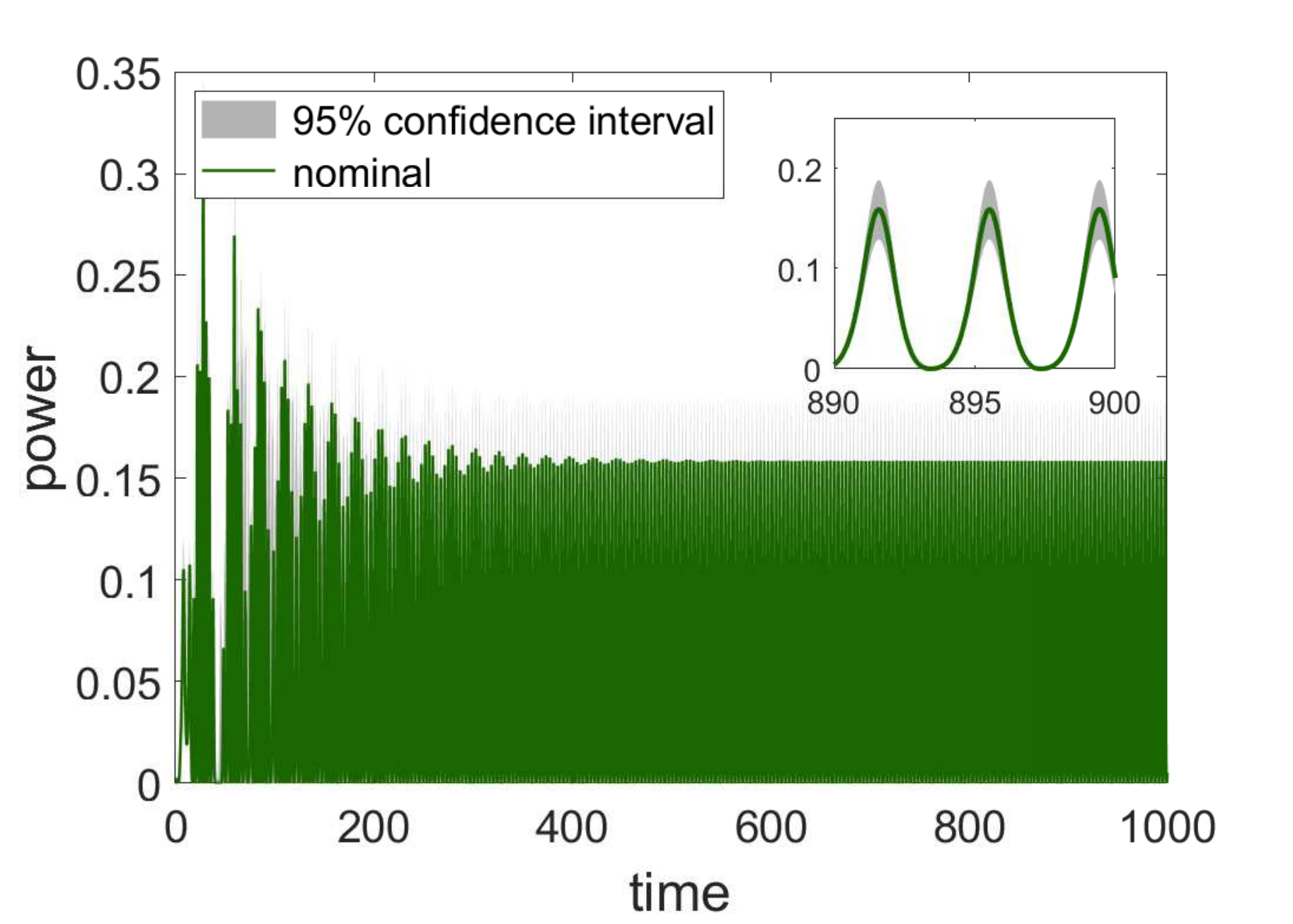}}
    \subfigure[$\mathnormal{f}=0.041$ for $\kappa$]{\includegraphics[width=0.32\textwidth]{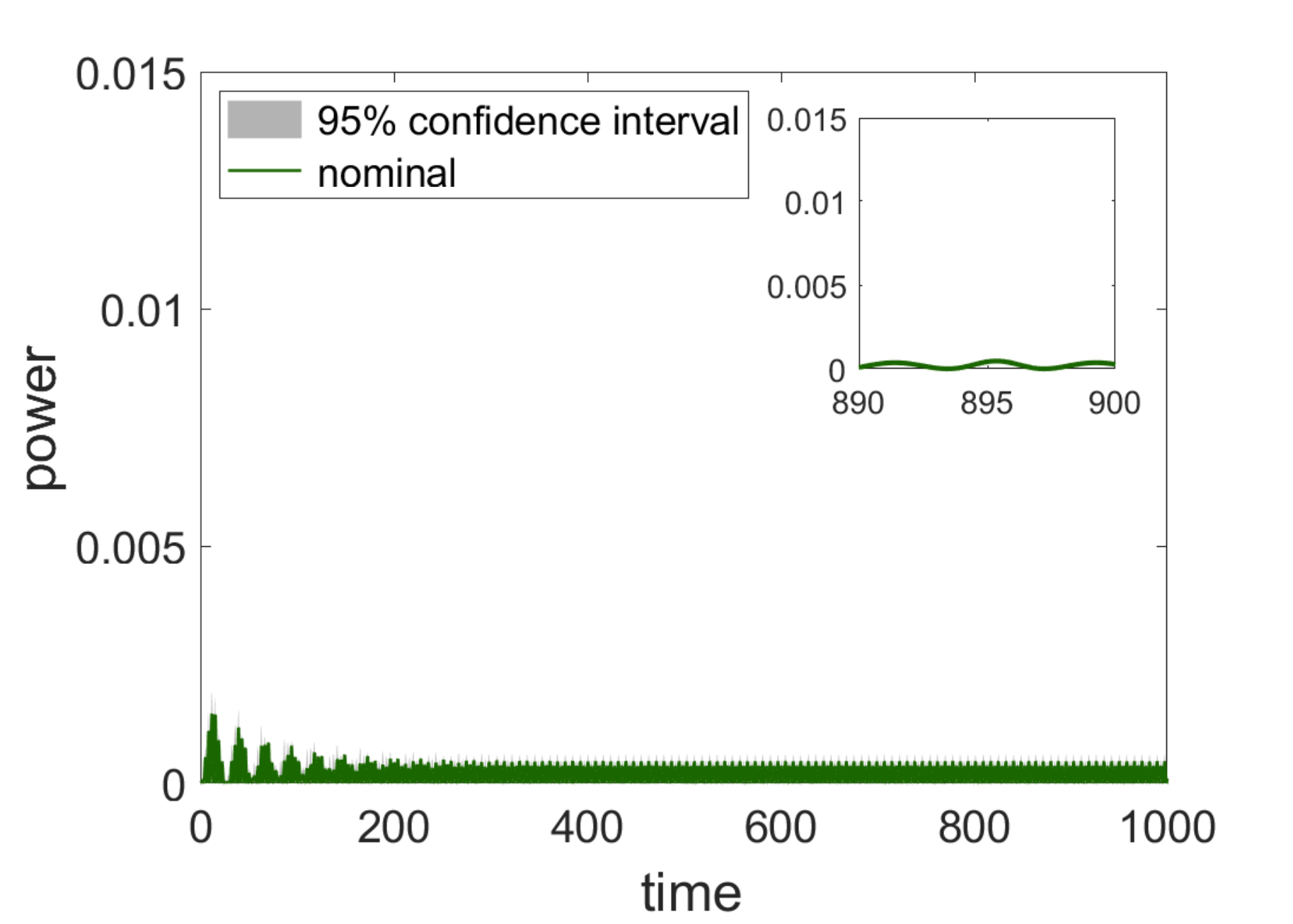}}
    \subfigure[$\mathnormal{f}=0.091$ for $\kappa$]{\includegraphics[width=0.32\textwidth]{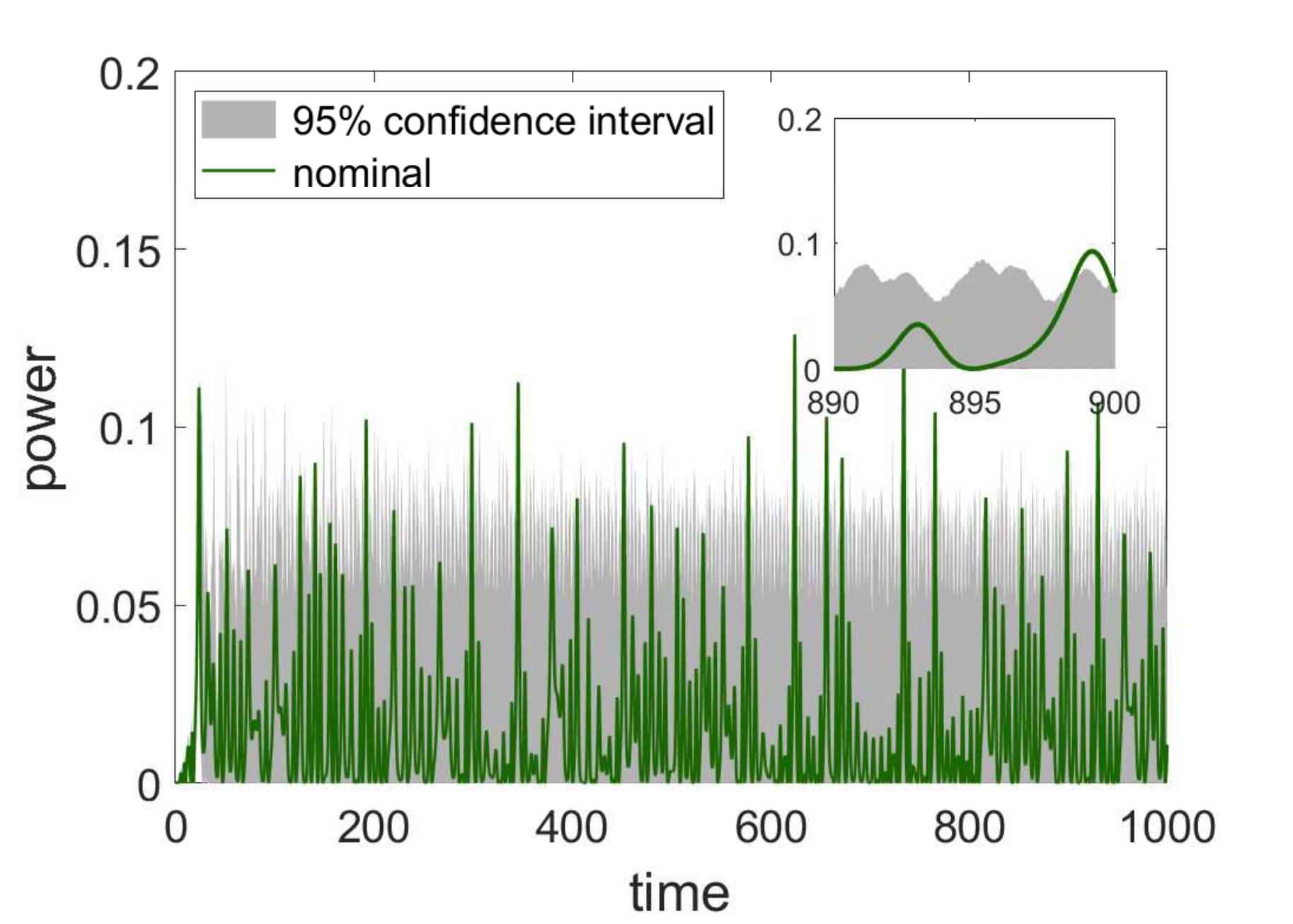}}
    \subfigure[$\mathnormal{f}=0.250$ for $\kappa$]{\includegraphics[width=0.32\textwidth]{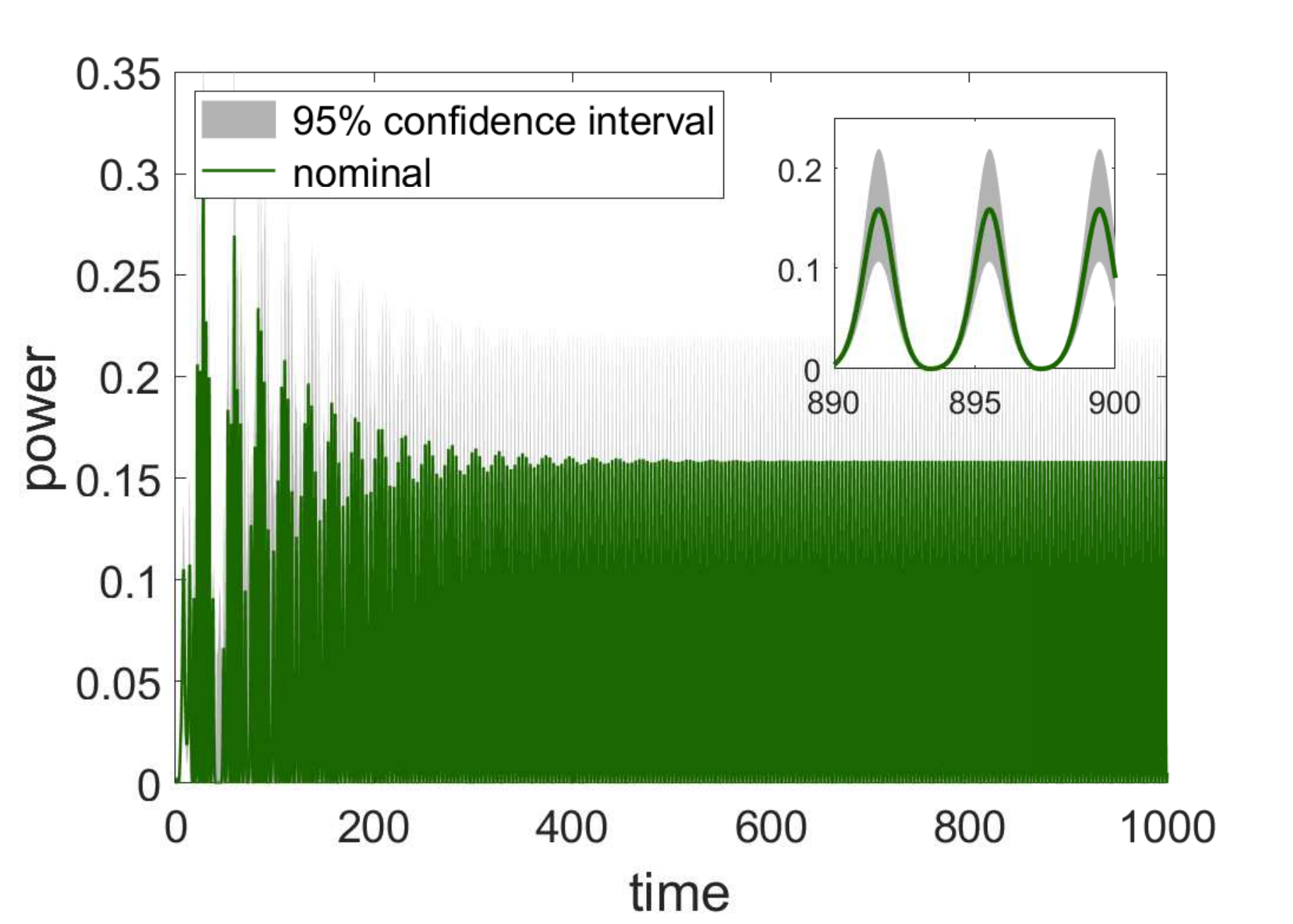}}
    \subfigure[$\mathnormal{f}=0.041$ for $f$]{\includegraphics[width=0.32\textwidth]{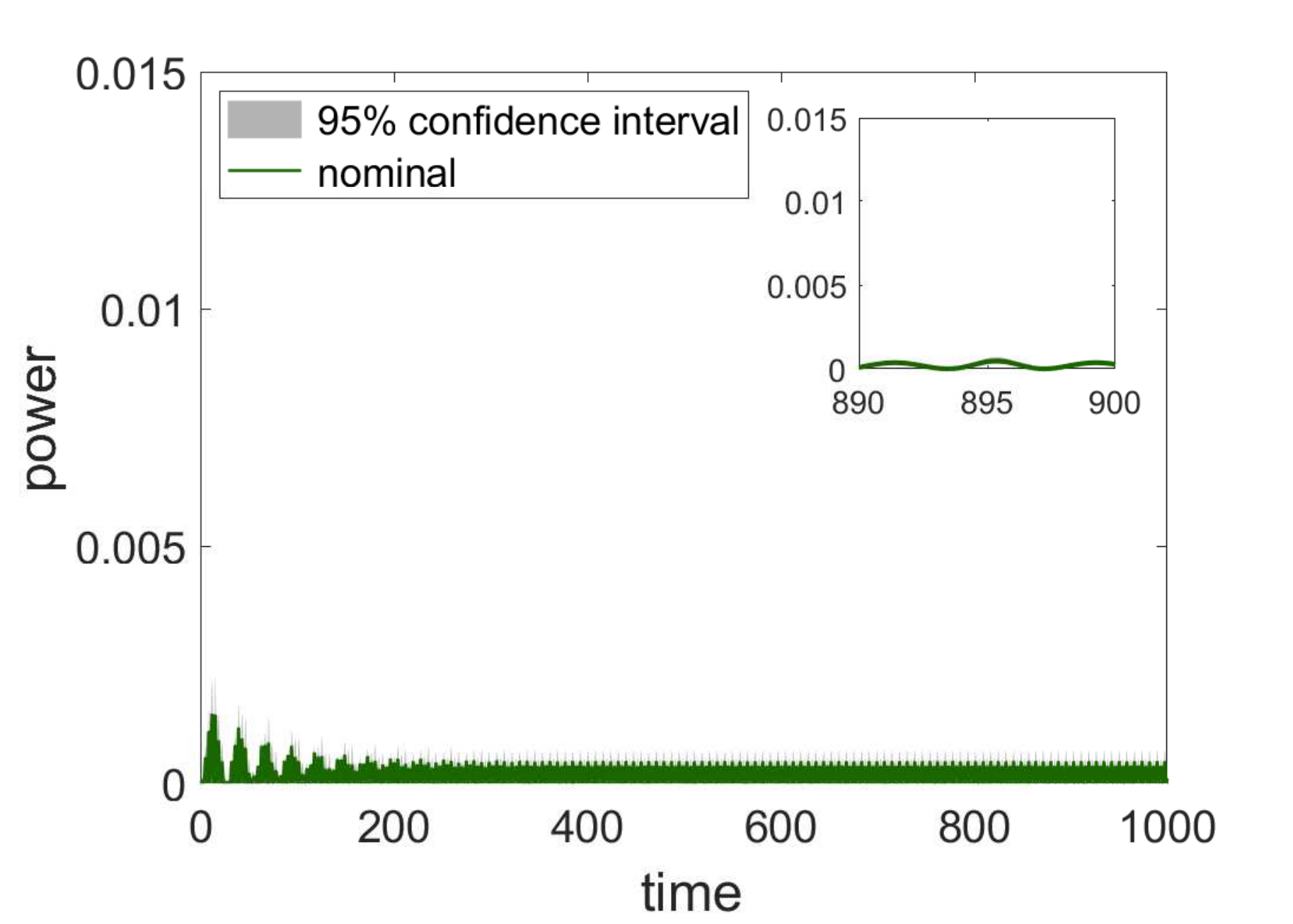}}
    \subfigure[$\mathnormal{f}=0.091$ for $f$]{\includegraphics[width=0.32\textwidth]{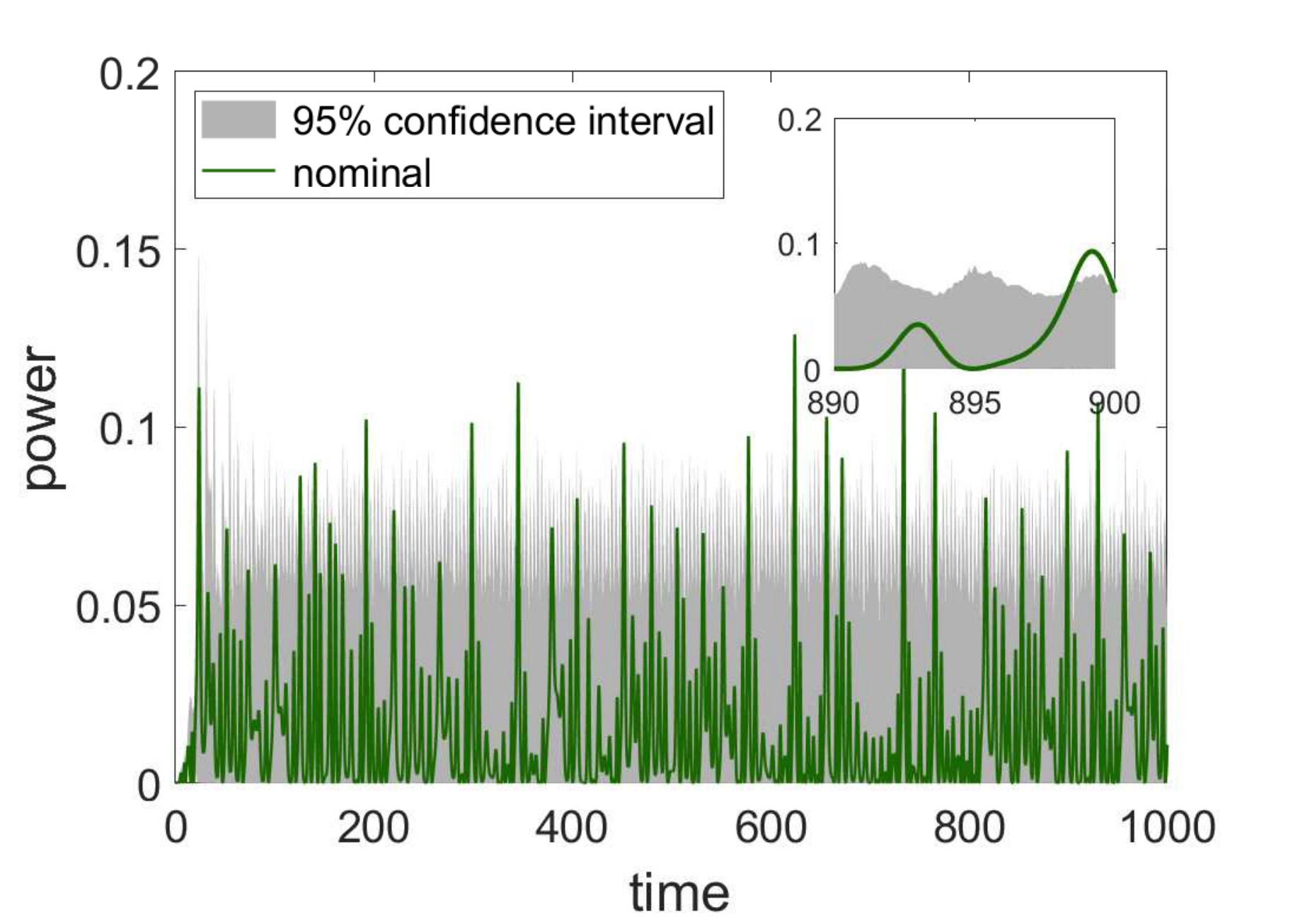}}
    \subfigure[$\mathnormal{f}=0.250$ for $f$]{\includegraphics[width=0.32\textwidth]{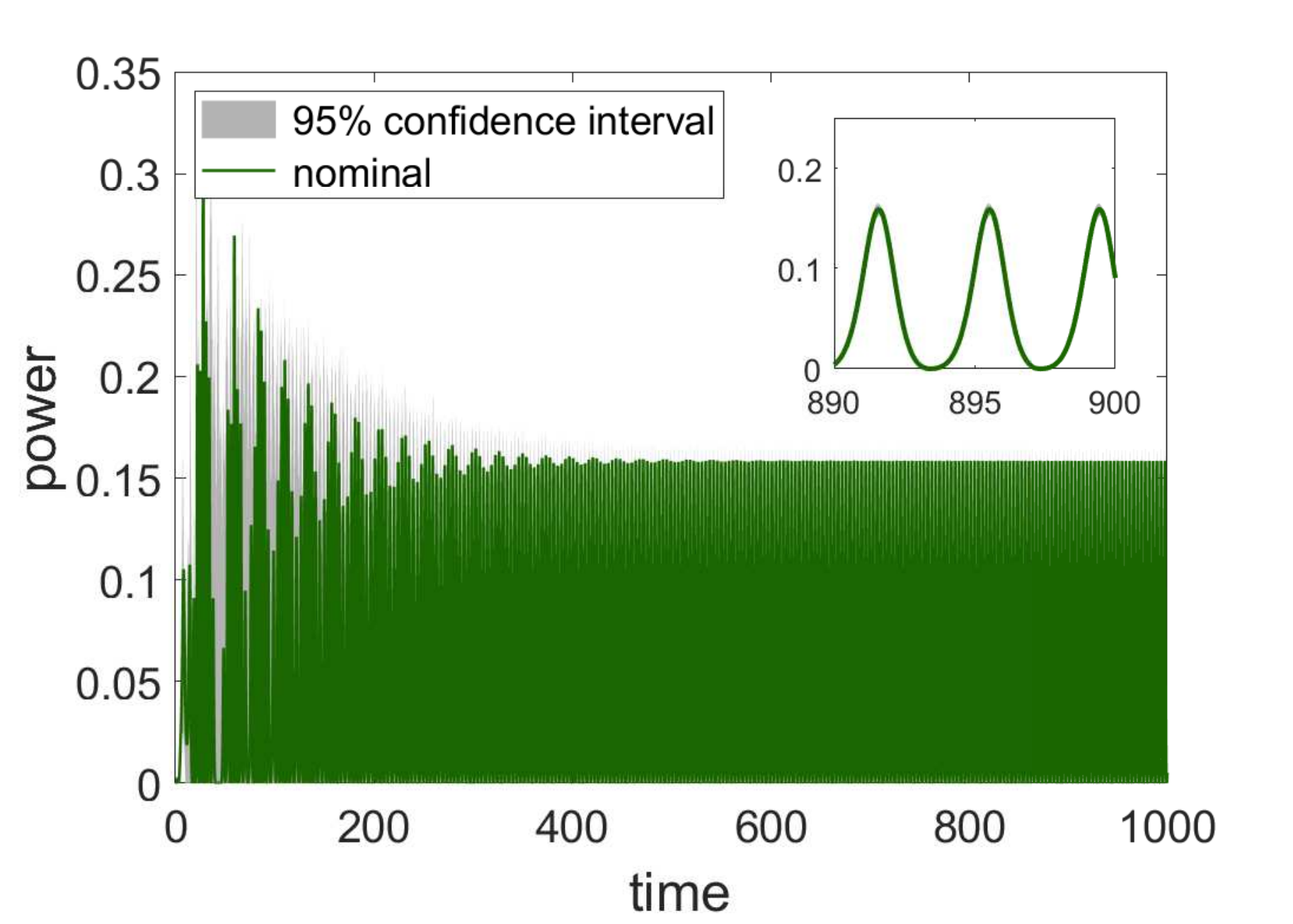}}
    \subfigure[$\mathnormal{f}=0.041$ for $\Omega$]{\includegraphics[width=0.32\textwidth]{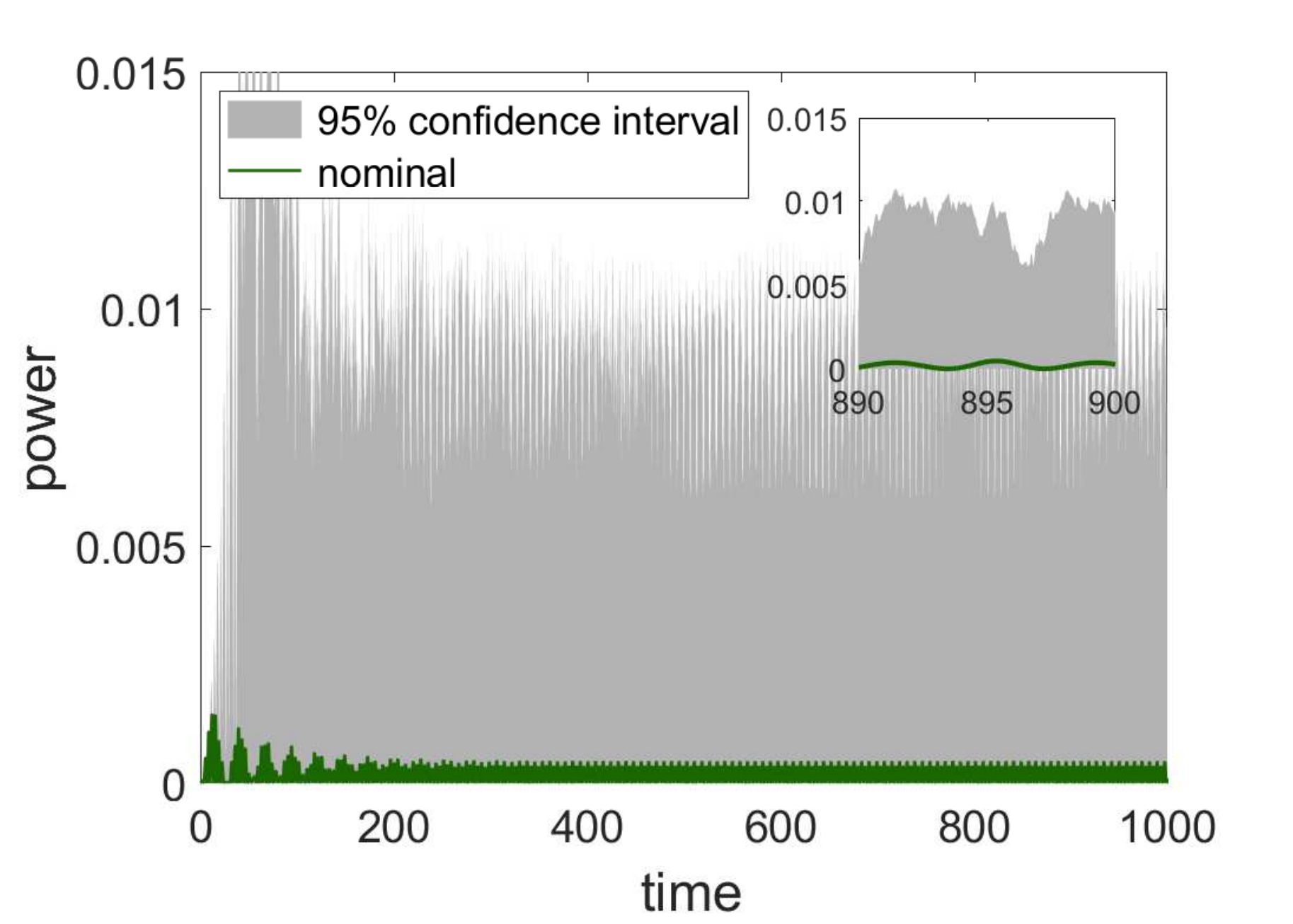}}
    \subfigure[$\mathnormal{f}=0.091$ for $\Omega$]{\includegraphics[width=0.32\textwidth]{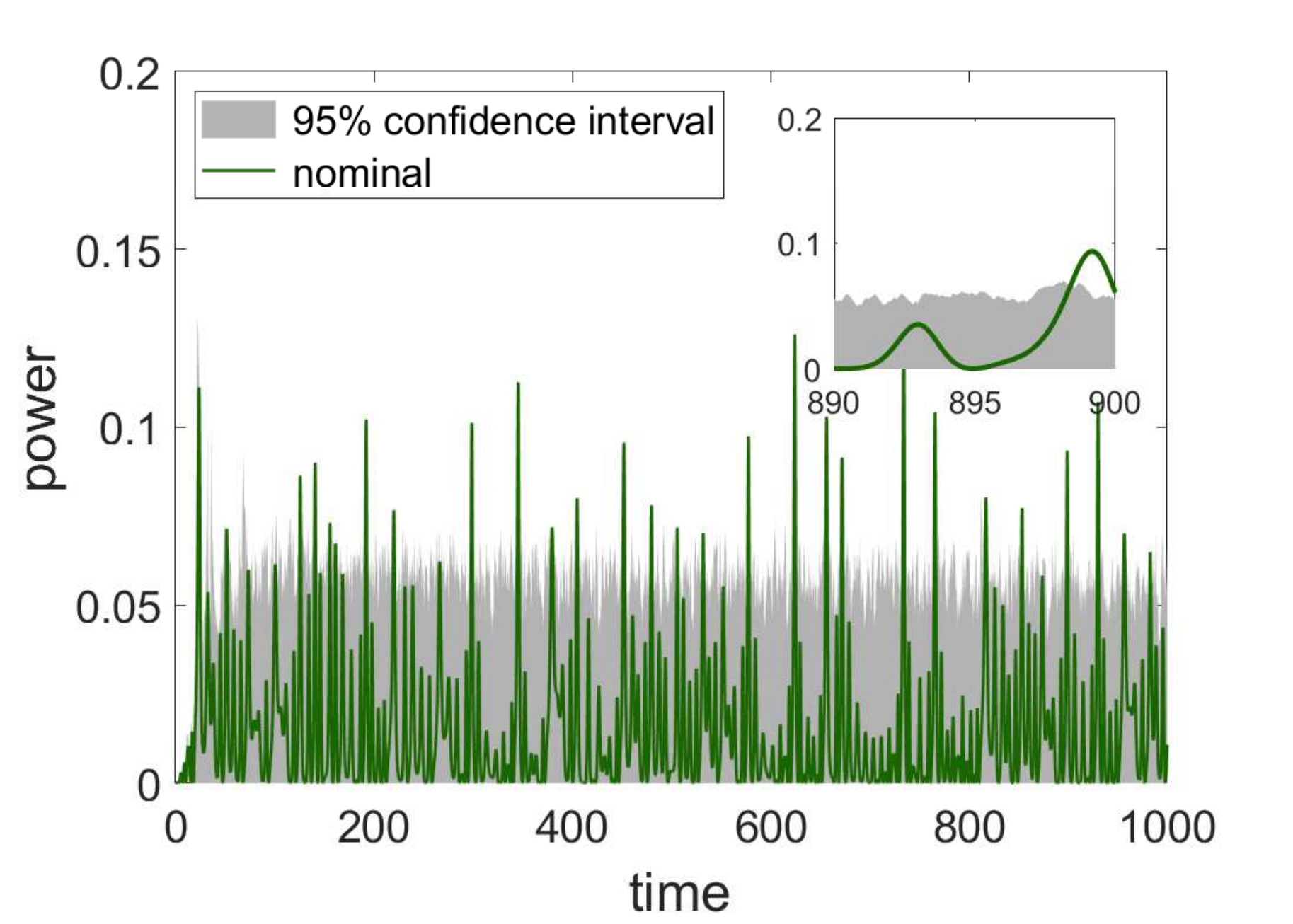}}
    \subfigure[$\mathnormal{f}=0.250$ for $\Omega$]{\includegraphics[width=0.32\textwidth]{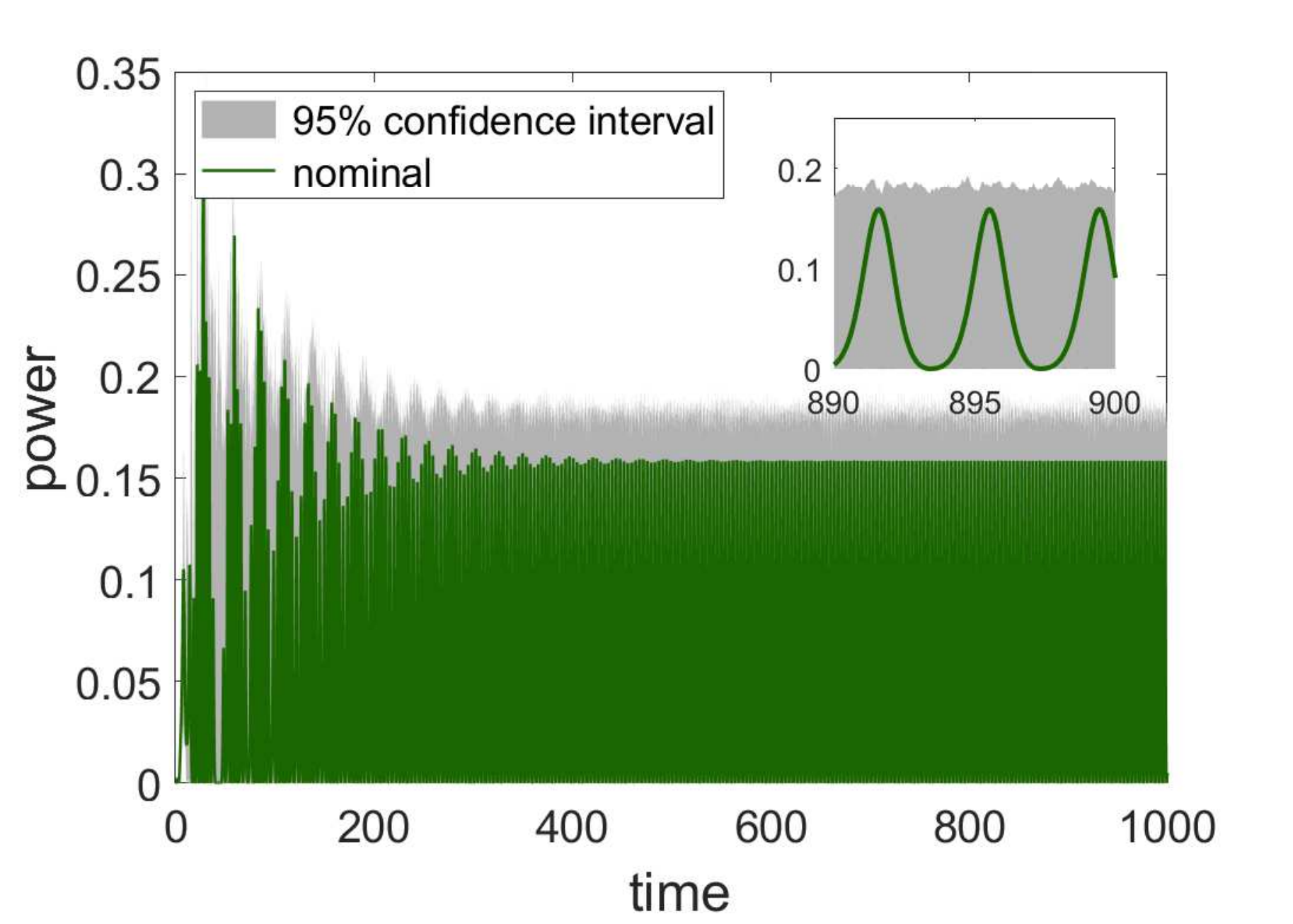}}
    \subfigure[$\mathnormal{f}=0.041$ for $\beta$]{\includegraphics[width=0.32\textwidth]{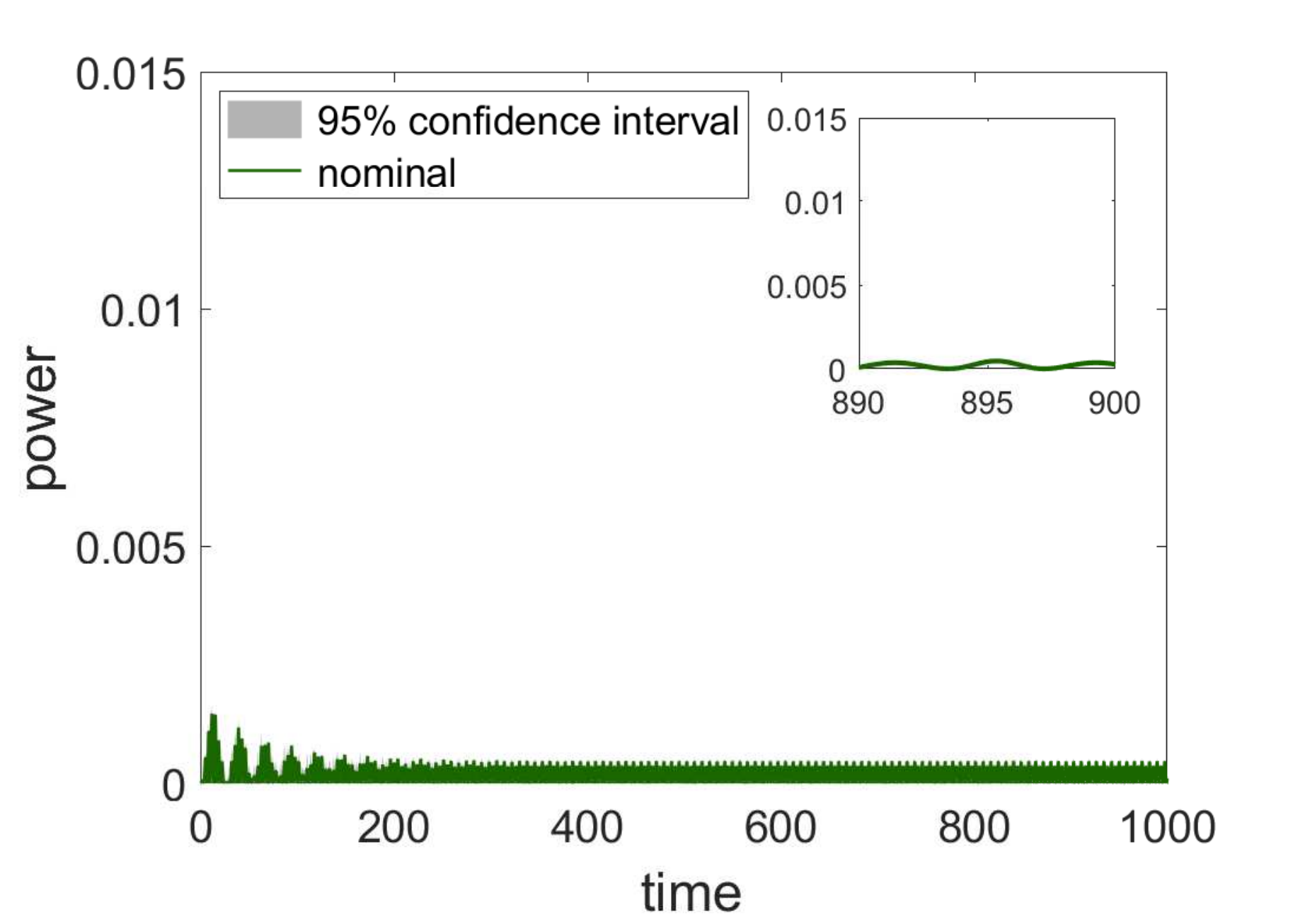}}
    \subfigure[$\mathnormal{f}=0.091$ for $\beta$]{\includegraphics[width=0.32\textwidth]{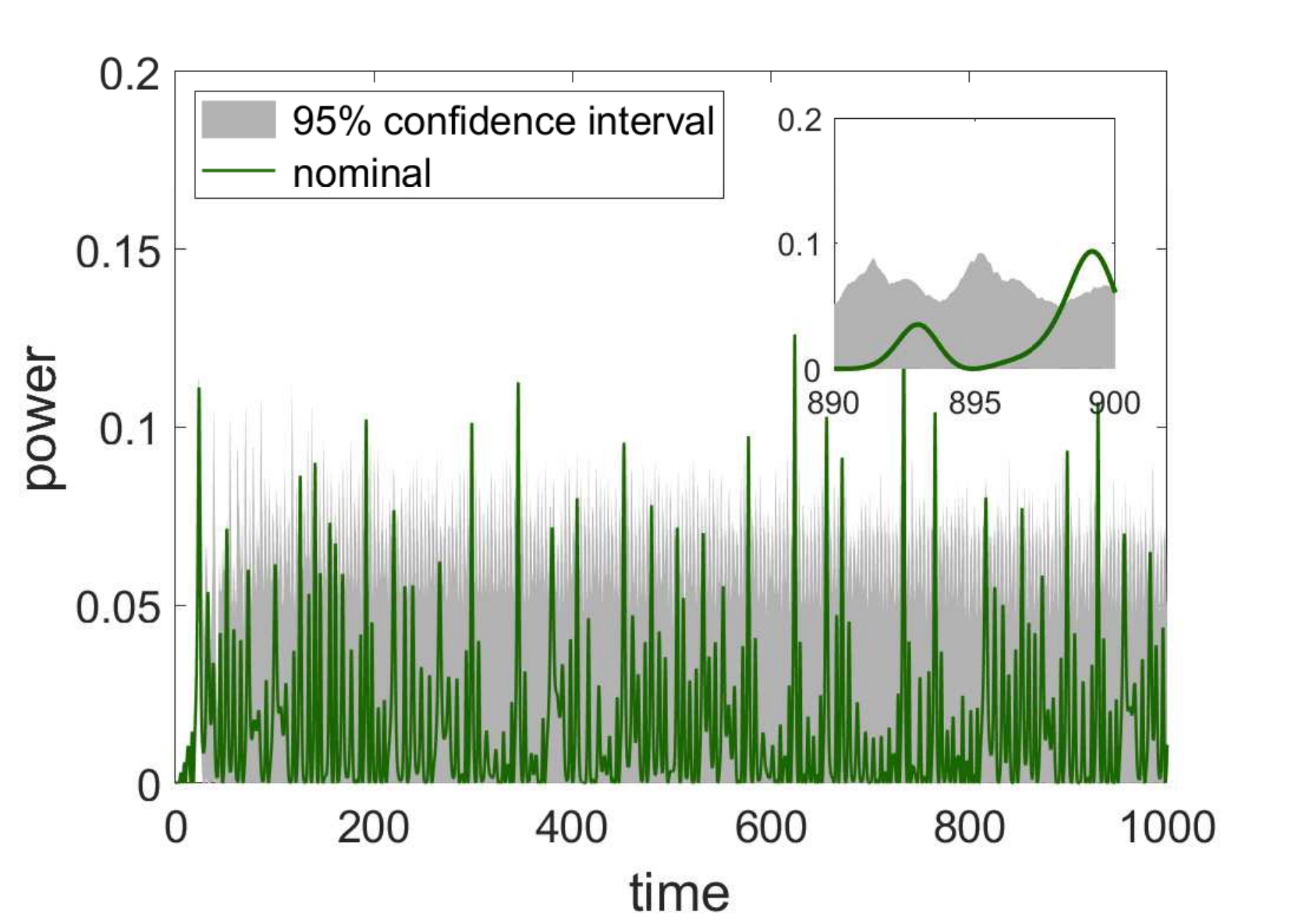}}
    \subfigure[$\mathnormal{f}=0.250$ for $\beta$]{\includegraphics[width=0.32\textwidth]{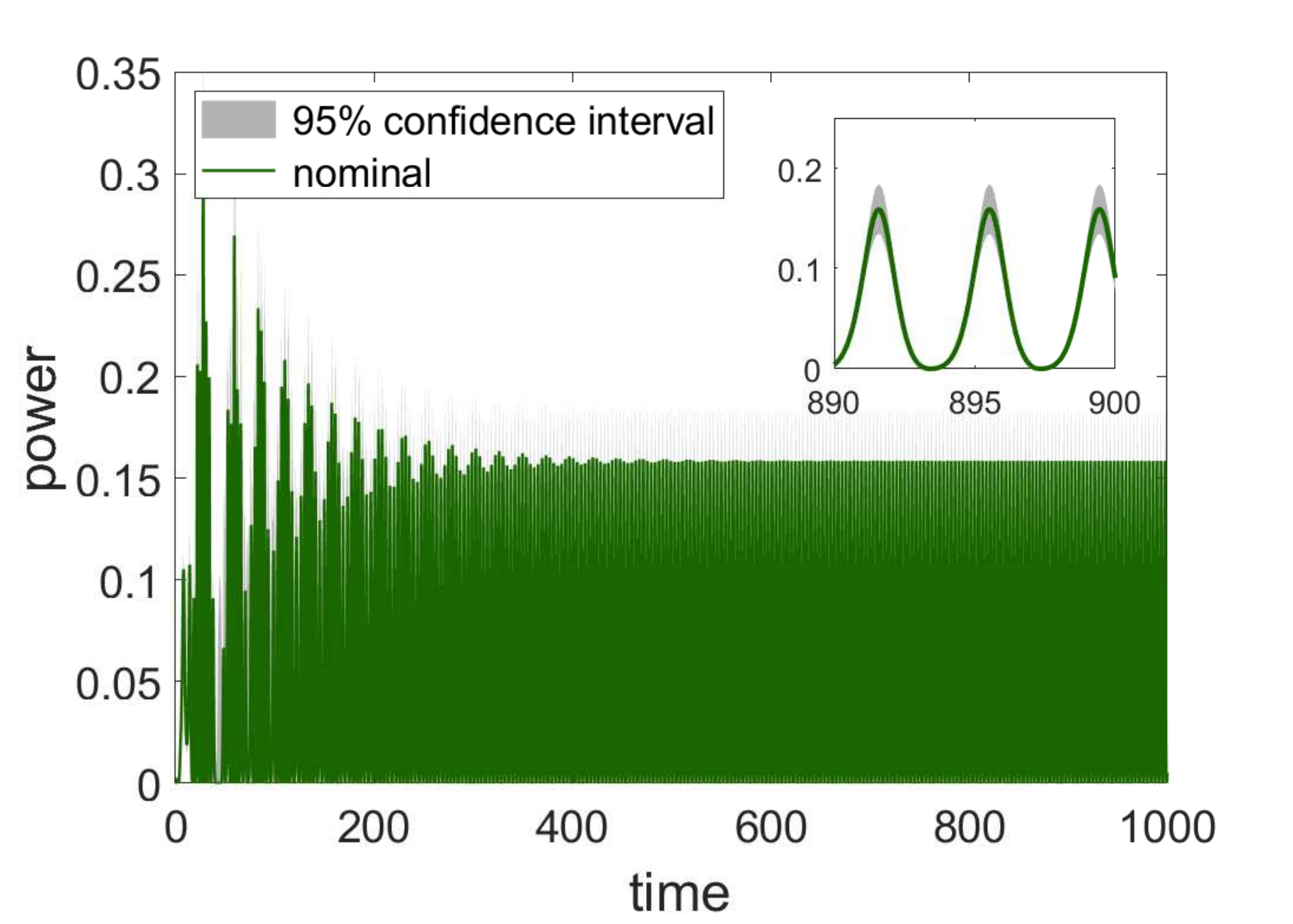}}
    \caption{Propagated uncertainty in the output power time series of the symmetric model with nonlinear piezoelectric coupling is shown under individual parameters: $\lambda$ (first row), $\kappa$ (second row), $\mathnormal{f}$ (third row), $\Omega$ (fourth row), and $\beta$ (fifth row). The columns are divided according to the different motion states of the system: intrawell (left), chaos (middle), and interwell (right).}
    \label{fig:up_BEHn}
\end{figure*}

Figure~\ref{fig:up_BEHn} illustrates the uncertainty propagation on the output power over time for each individual parameter. The plot represents the 95\% confidence interval as well as the nominal series. The left column corresponds to intrawell motion, where $\beta$ generates insignificant variability. Also, $\lambda$, $\kappa$, and $\mathnormal{f}$ exhibit slight variation in the nominal output. In contrast, $\Omega$ remains fundamental to changing the response, leading to an extensive region of confidence. The middle column shows the results for chaotic motion, where all parameters alter the confidence interval. Finally, in the right column, interwell motion, $\kappa$, and $\Omega$ variations affect the mean power similarly to the previous model. However, variations of $\beta$ generate an envelope primarily for the signal amplitude, with little influence on the harvested power's increase or decrease. According to Fig.~\ref{fig:prob_BEHn}, the upper band part is when $\beta$ increases.

\subsection{Asymmetric bistable energy harvester with nonlinear piezoelectric coupling }

\begin{figure*}
     \centering
     \subfigure[$\mathnormal{f}=0.041$]{\includegraphics[width=0.32\textwidth]{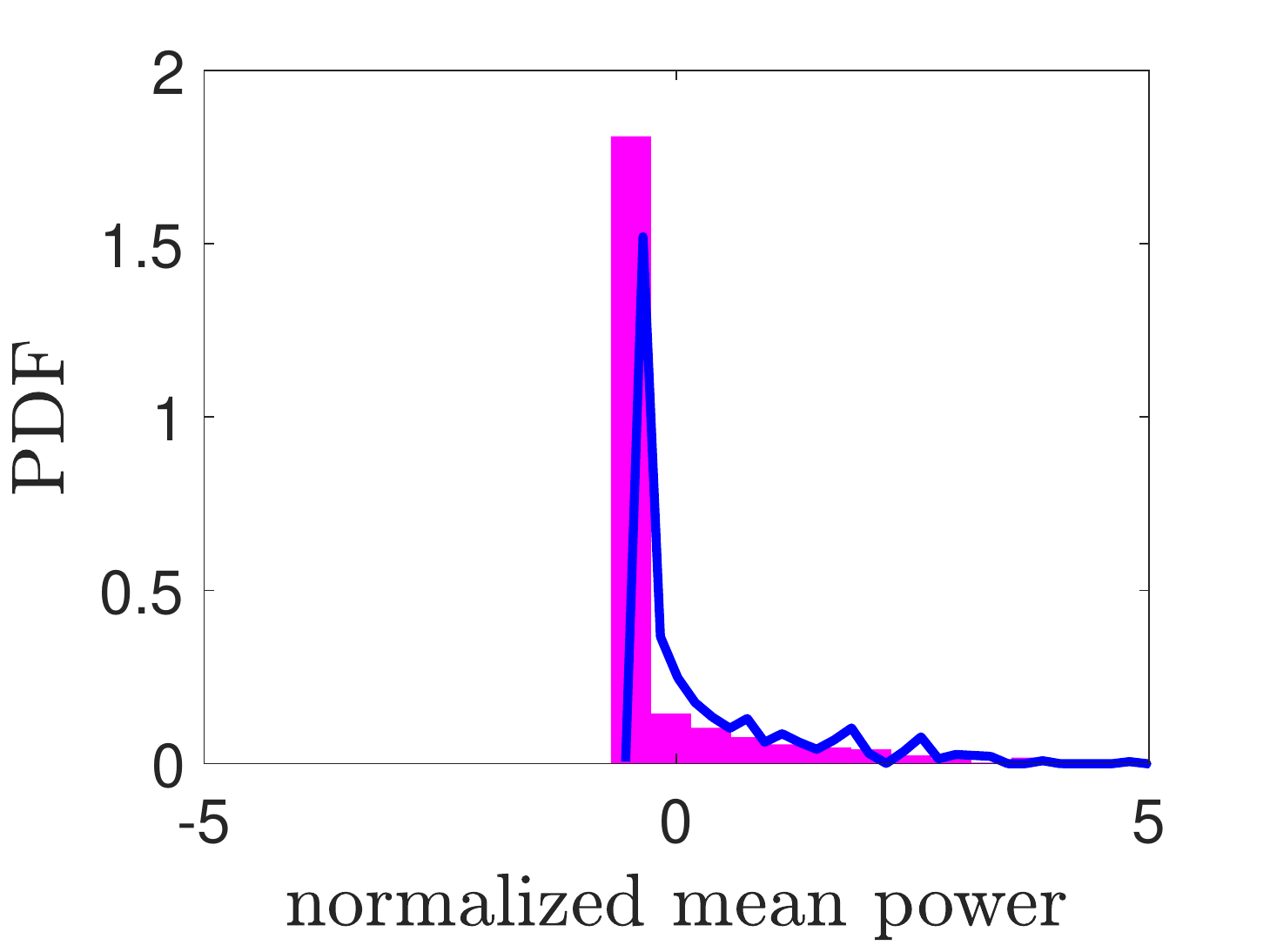}}
     \subfigure[$\mathnormal{f}=0.060$]{\includegraphics[width=0.32\textwidth]{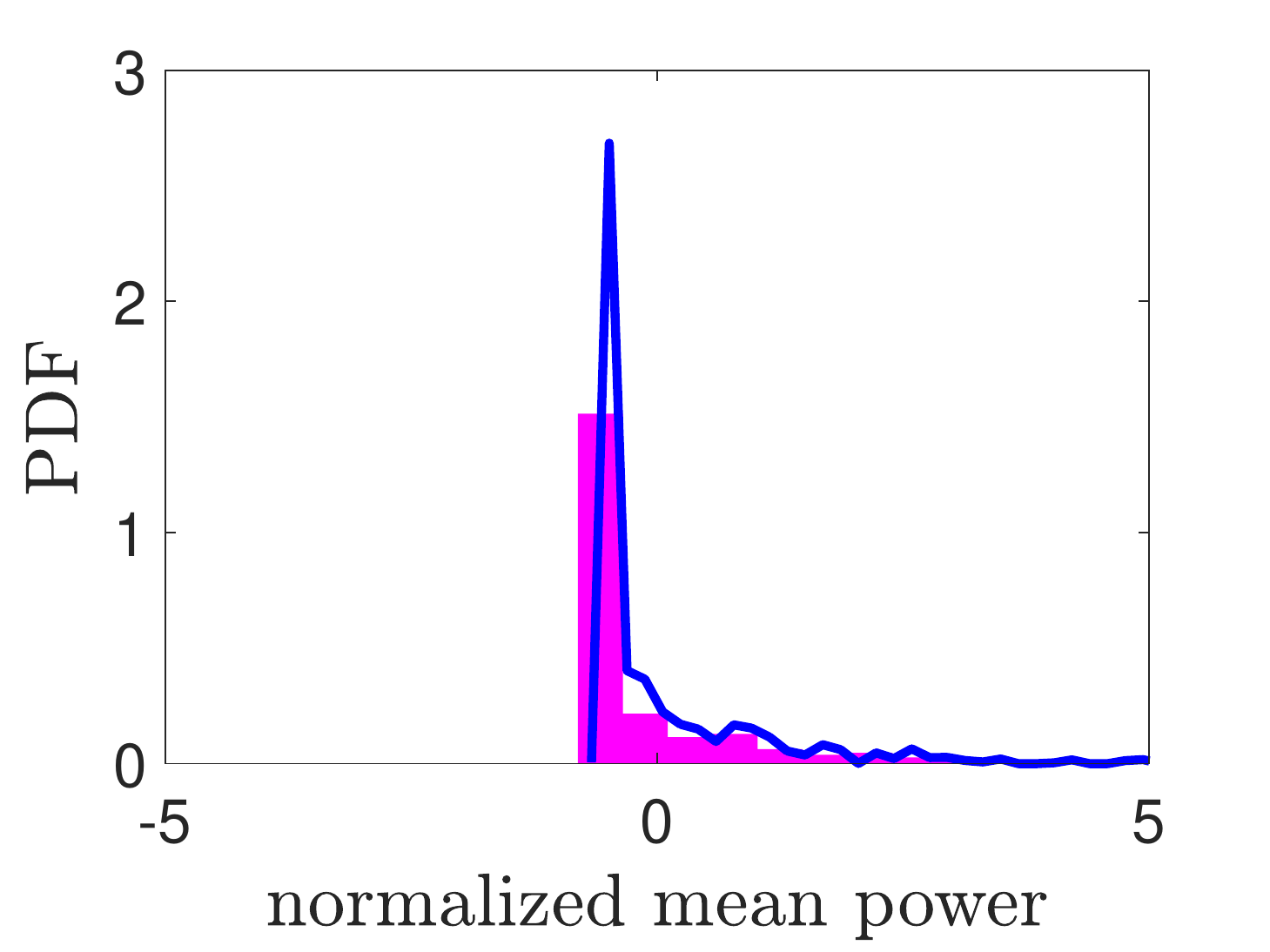}}
     \subfigure[$\mathnormal{f}=0.083$]{\includegraphics[width=0.32\textwidth]{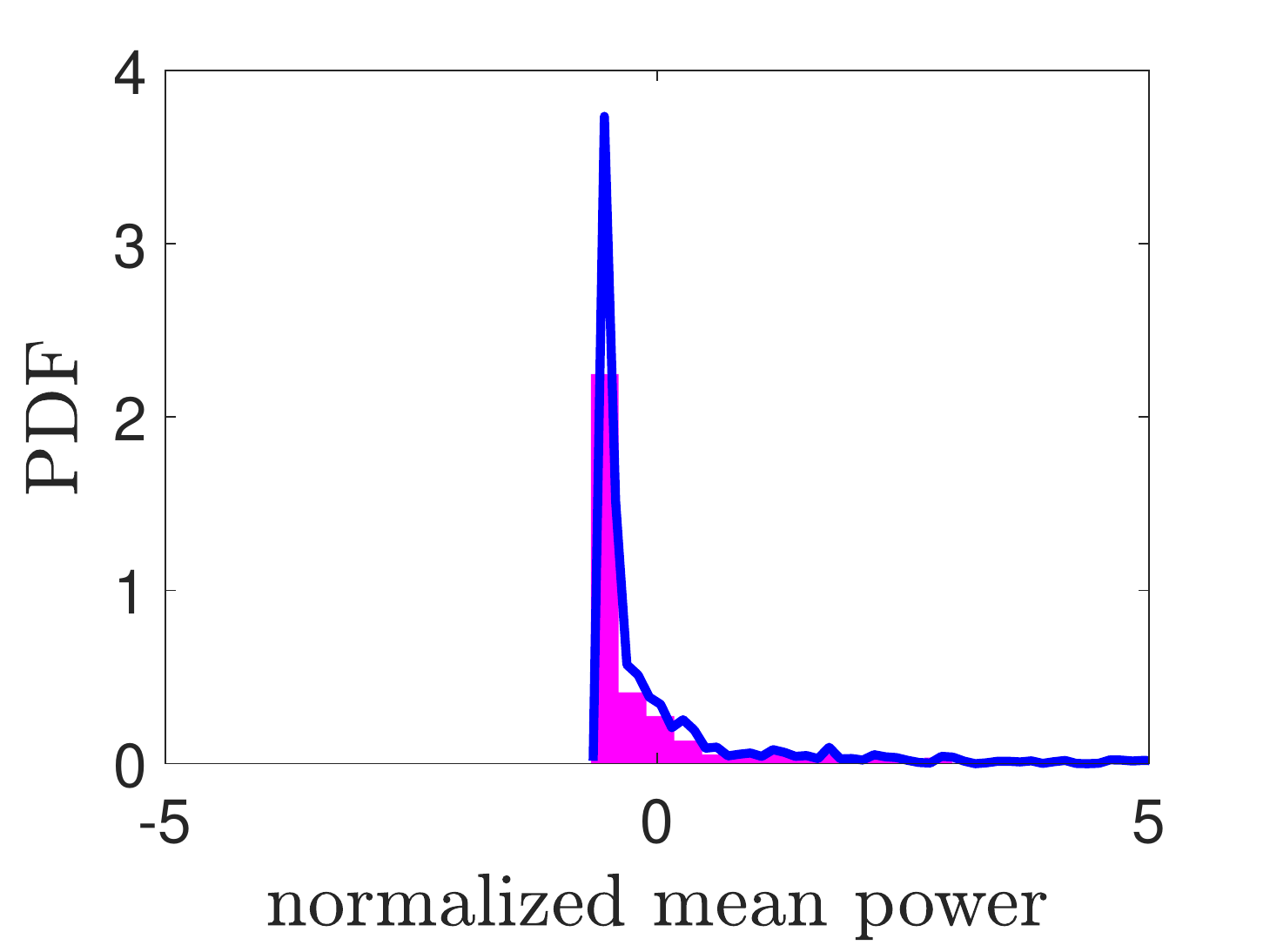}}
     \subfigure[$\mathnormal{f}=0.091$]{\includegraphics[width=0.32\textwidth]{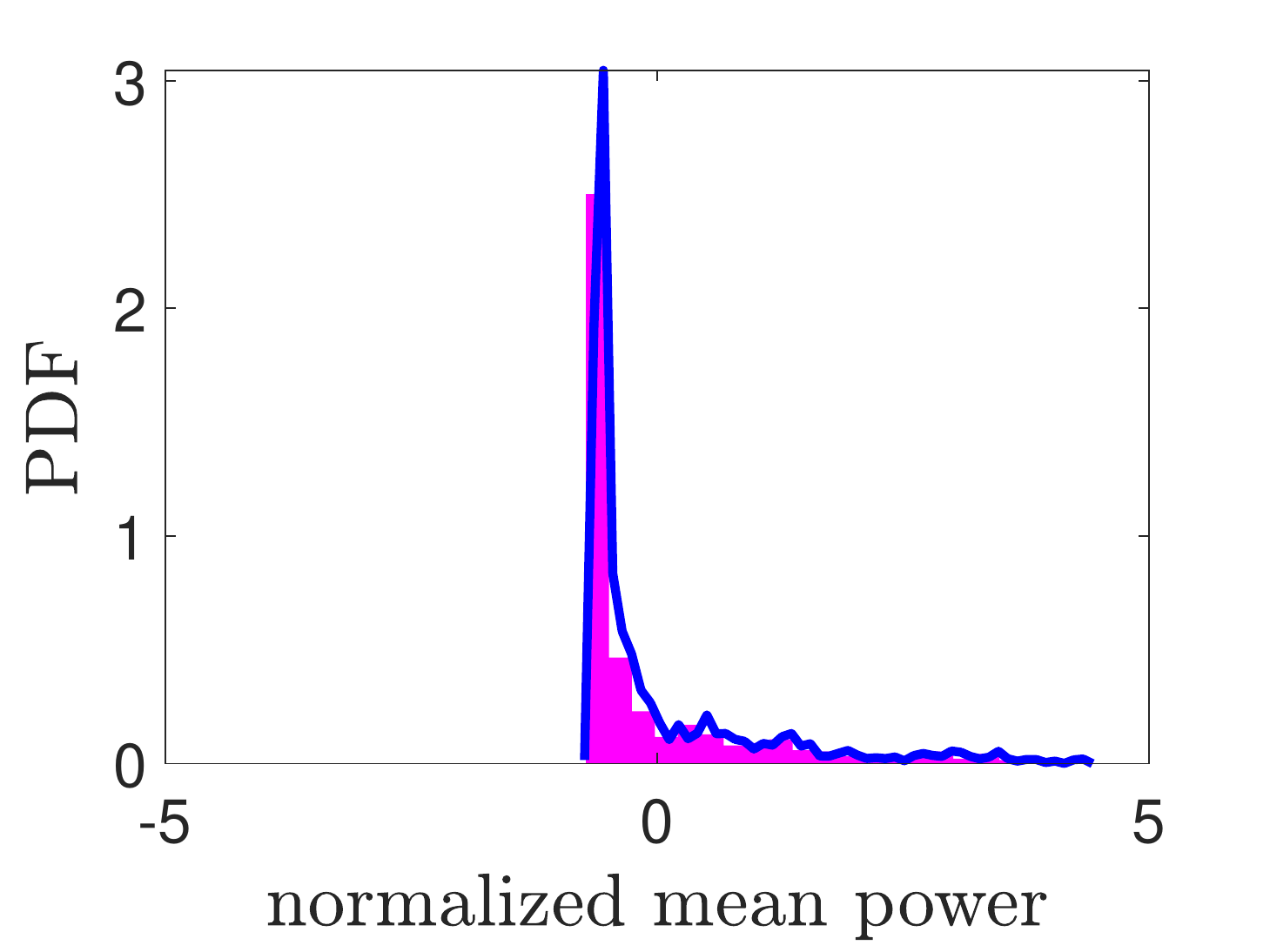}}
     \subfigure[$\mathnormal{f}=0.105$]{\includegraphics[width=0.32\textwidth]{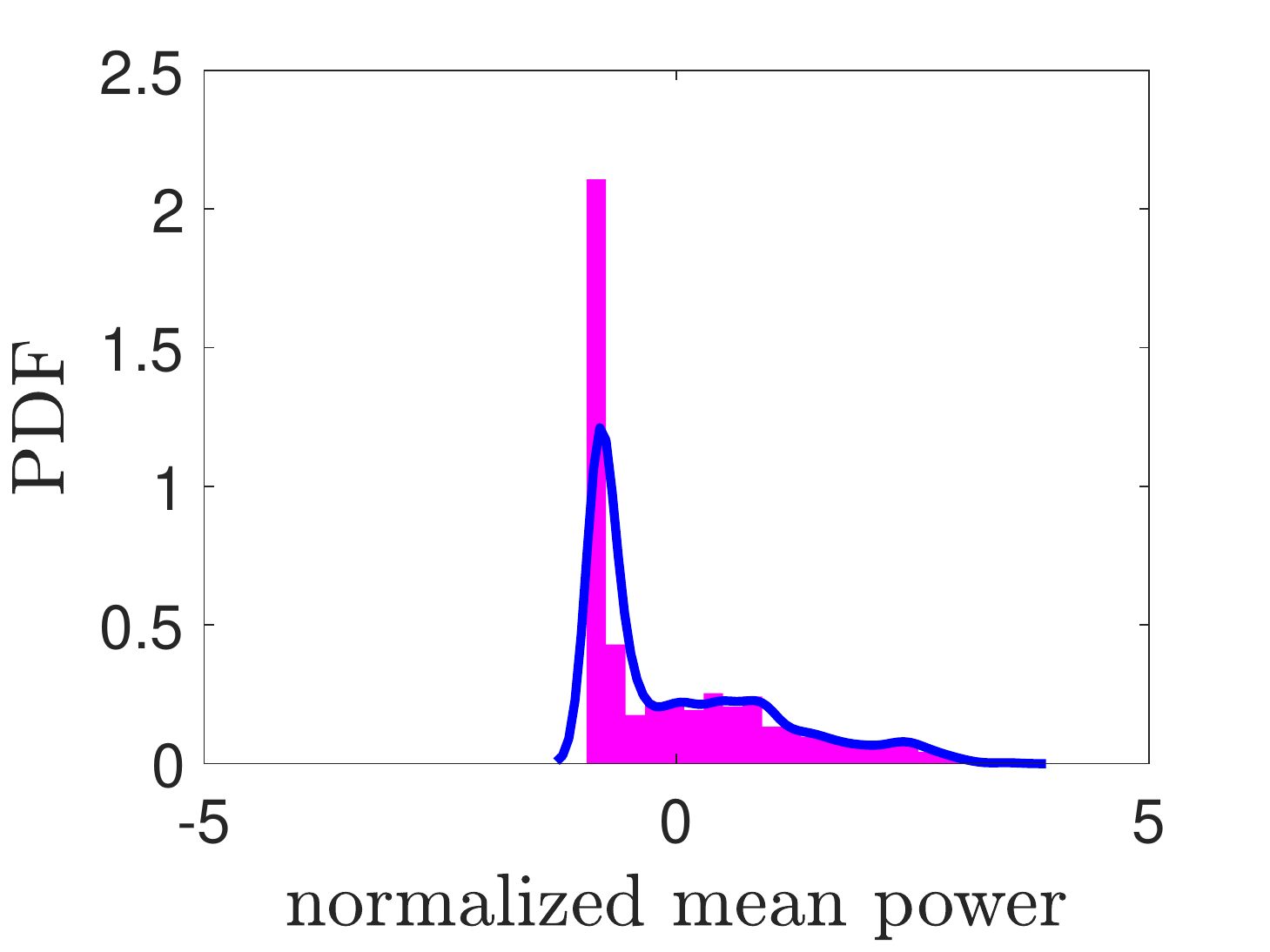}}
     \subfigure[$\mathnormal{f}=0.115$]{\includegraphics[width=0.32\textwidth]{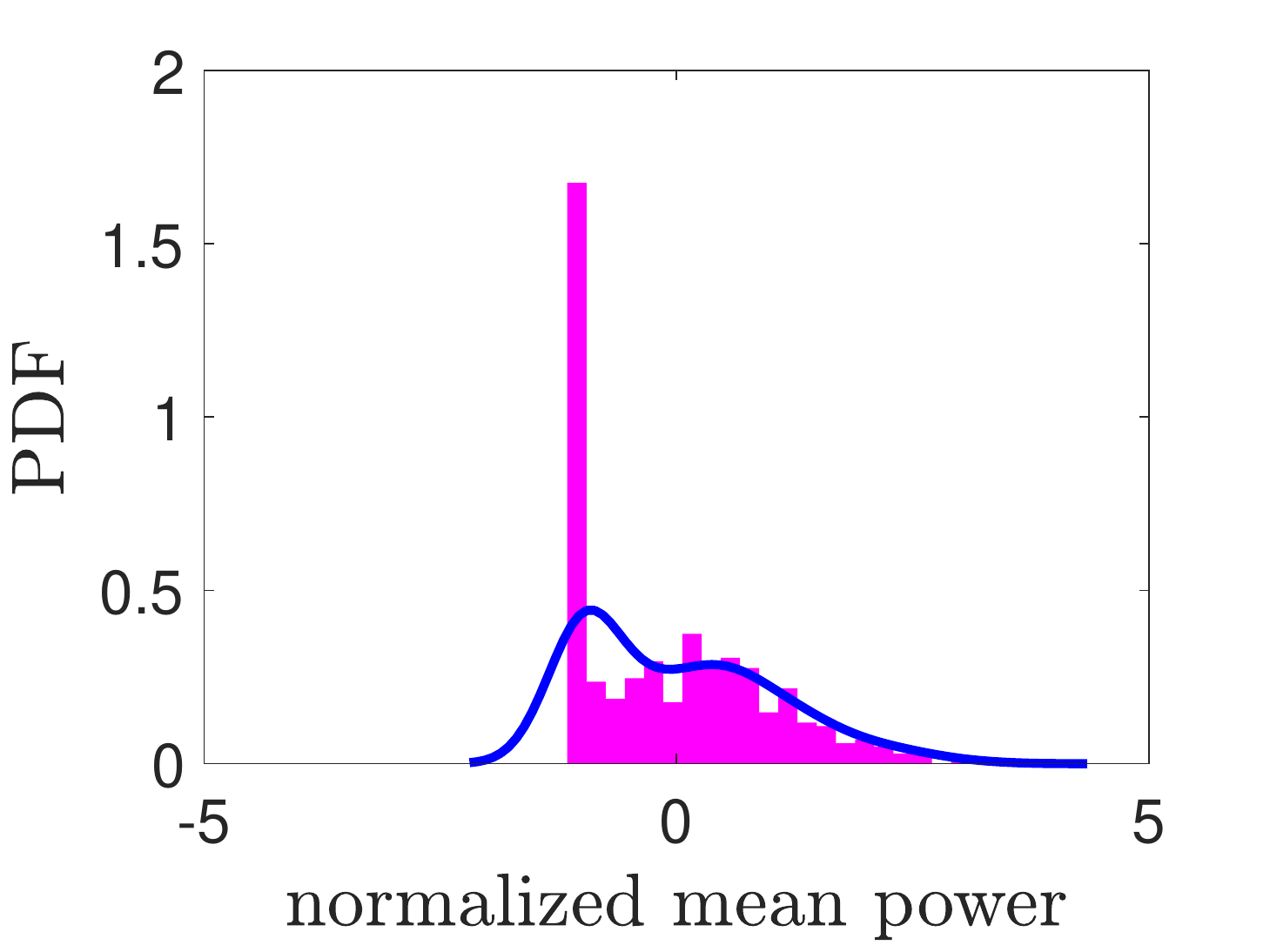}}
     \subfigure[$\mathnormal{f}=0.147$]{\includegraphics[width=0.32\textwidth]{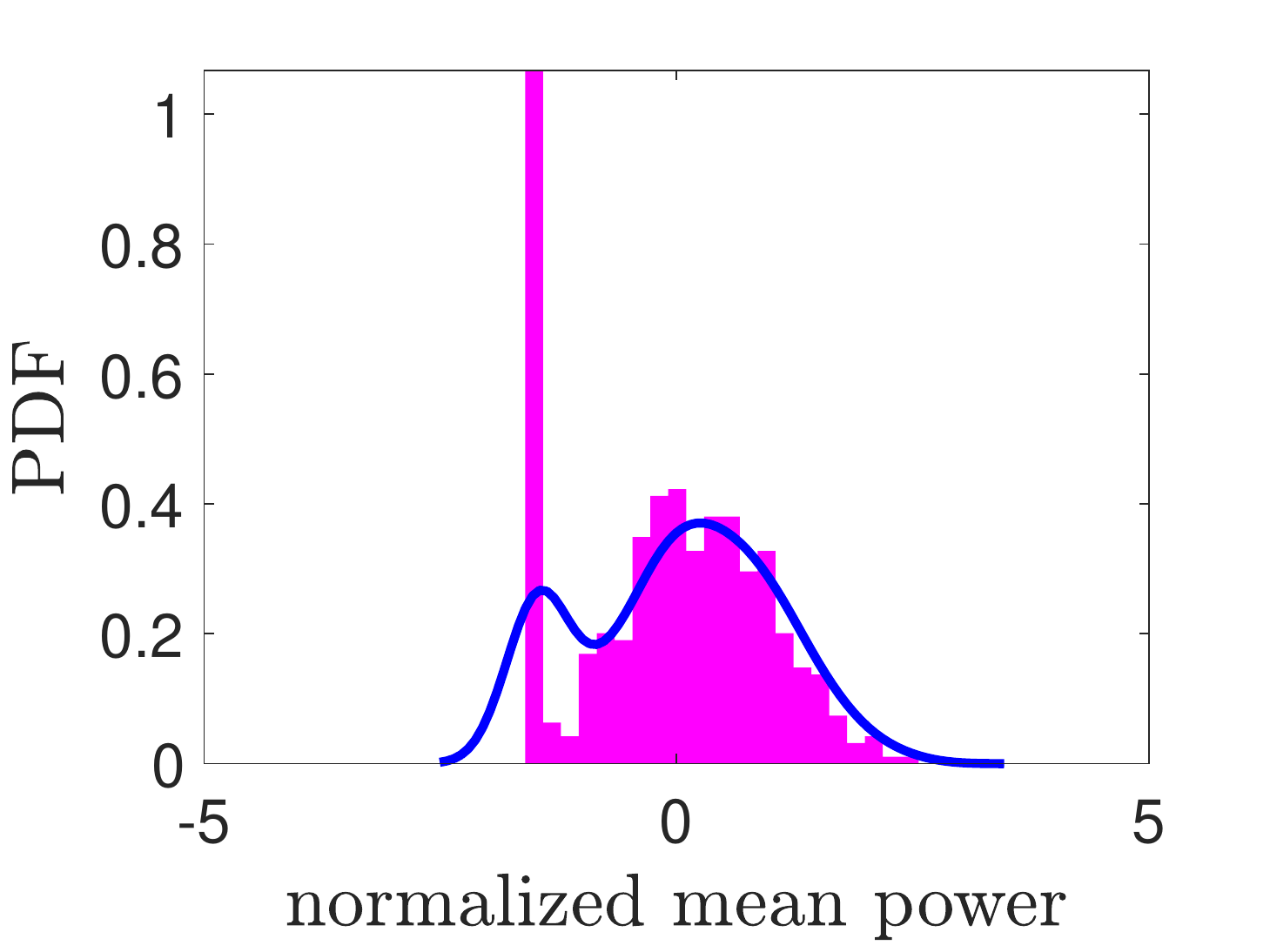}}
     \subfigure[$\mathnormal{f}=0.200$]{\includegraphics[width=0.32\textwidth]{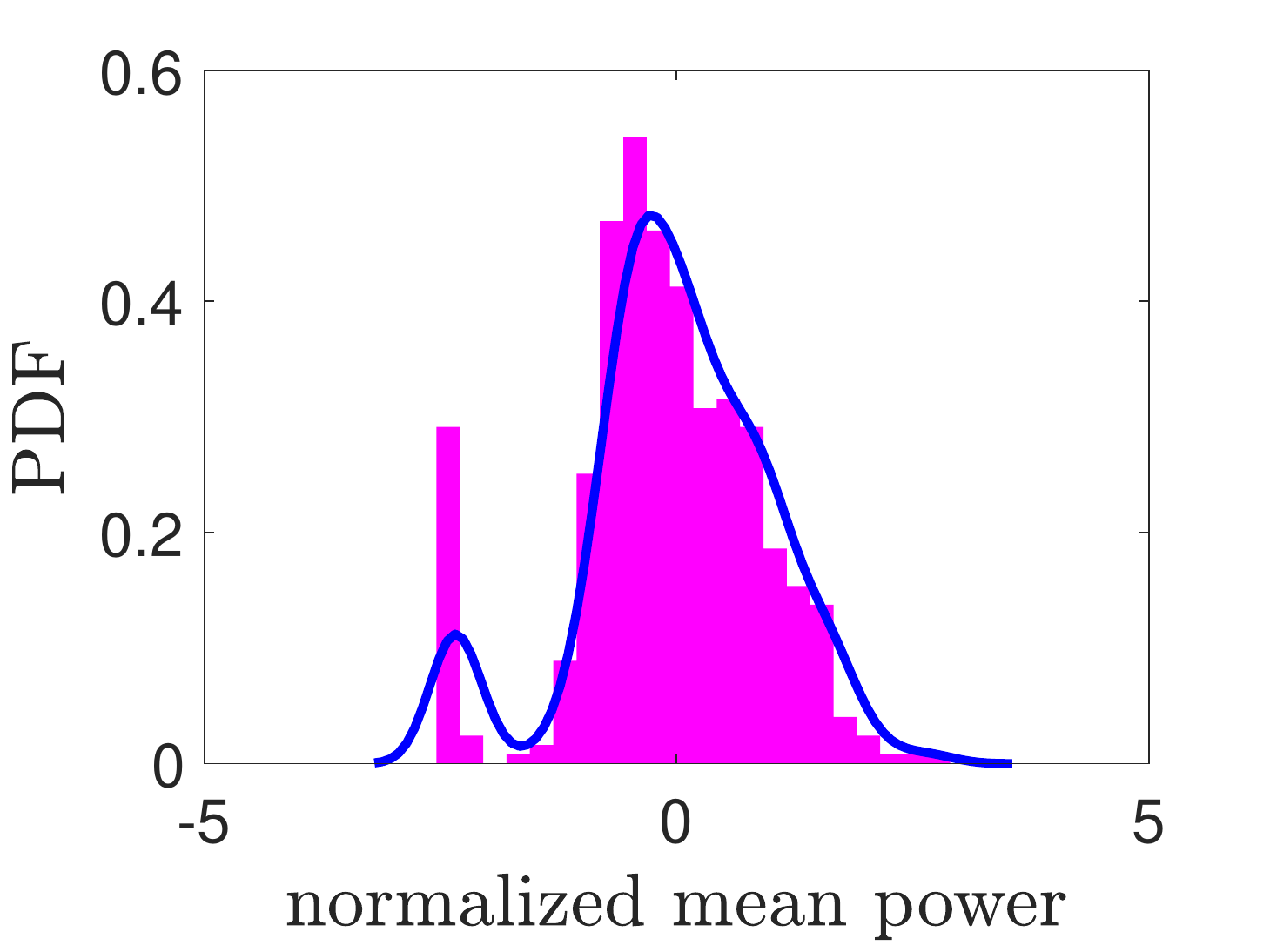}}
     \subfigure[$\mathnormal{f}=0.250$]{\includegraphics[width=0.32\textwidth]{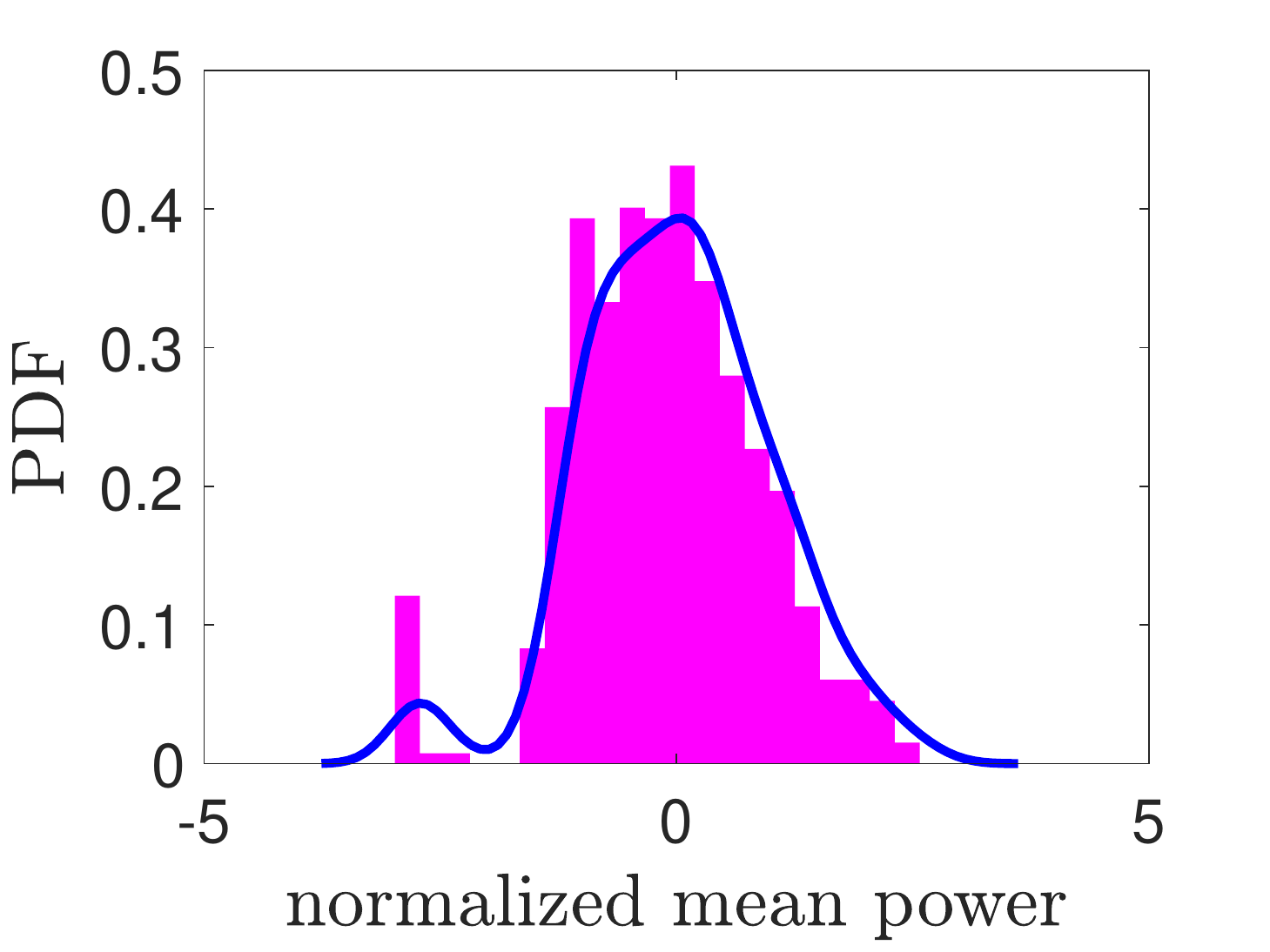}}
     \caption{Probability density function of the normalized mean power for the asymmetric energy harvester model with nonlinear piezoelectric coupling under different excitation amplitudes. The kernel density function is represented by the blue line.}
     \label{fig:pdf_pmeha}
\end{figure*}

Figure~\ref{fig:pdf_pmeha} displays the probability density function of the normalized mean power for an asymmetric bistable energy harvester with nonlinear coupling ($\delta \neq 0$, $\phi \neq 0$ and $\beta \neq 0$) across a range of excitation amplitudes. For low values of $\mathnormal{f}$ ($<0.091$), the distribution takes the form of an exponential function. As the amplitude of excitation increases, a bimodal distribution is observed. Firstly, the negative peak is higher, but for $\mathnormal{f}>0.147$ the positive peak becomes more prominent, indicating that the bistable motion occurs frequently. In contrast to the symmetric model, the distribution also reveals some regions with low probability that negatively impact power harvesting, suggesting monostable motions. If the excitation force were further increased, the negative peak would disappear, and the distribution would become unimodal. To obtain a unimodal distribution with a high mean, it is necessary to increase the excitation force beyond the analyzed range.

\begin{figure*}
    \centering
    \includegraphics[width=1\textwidth]{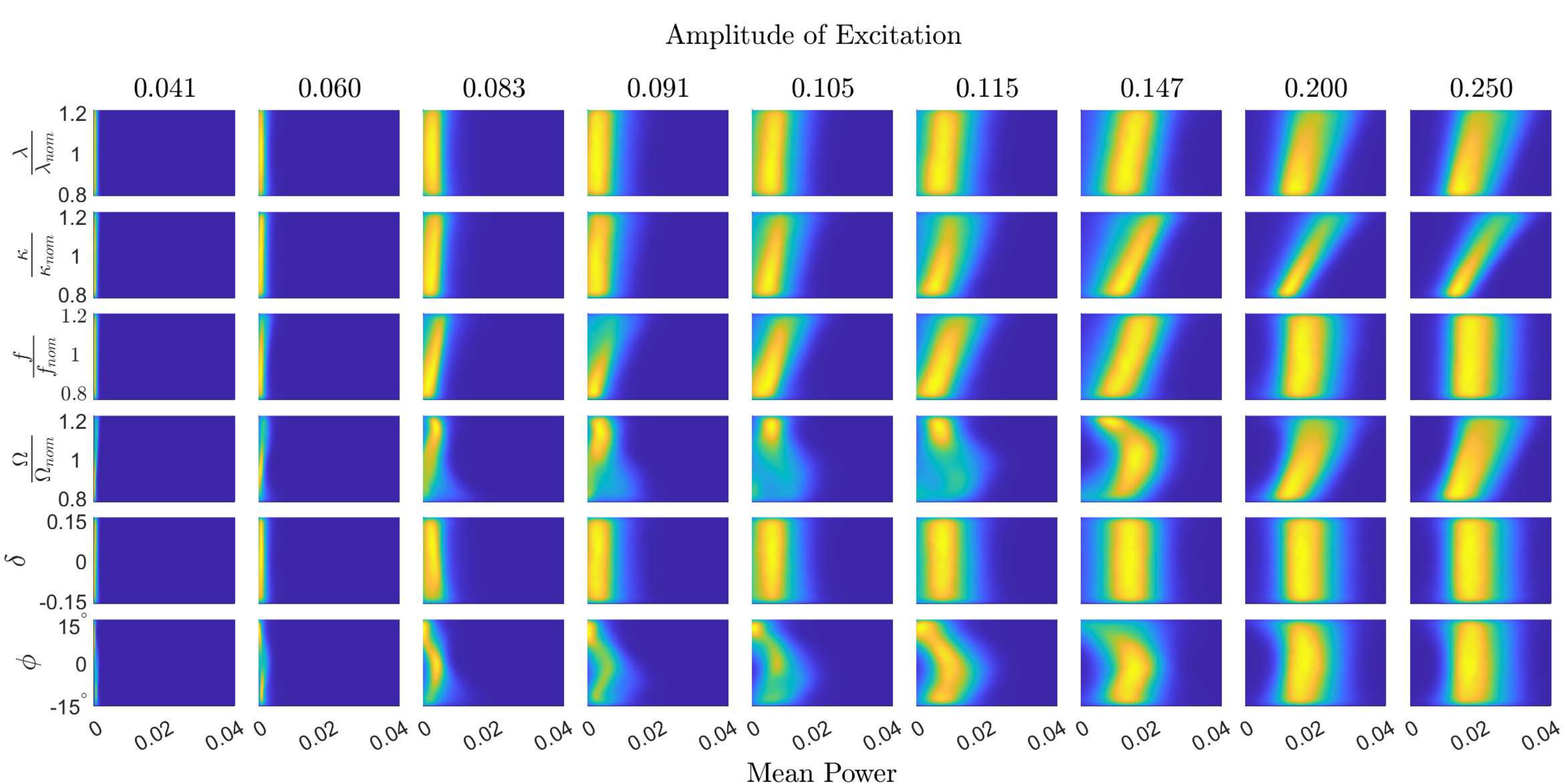}
    \caption{Joint-CDF of mean power conditioned on each parameter of interest ($\lambda$, $\kappa$, $\mathnormal{f}$, $\Omega$, $
    \delta$, $\phi$) under different values of excitation amplitude for the asymmetric model with nonlinear piezoelectric coupling.}
    \label{fig:joint_BEHa}
\end{figure*}

Figure~\ref{fig:joint_BEHa} displays the joint-CDF of the mean power conditioned on each parameter of interest under different nominal excitation conditions for the asymmetric bistable energy harvester with nonlinear coupling. The influence of the asymmetry coefficient ($\delta$) on the mean power is negligible, while the angle ($\phi$) increases the mean power as it approaches zero, suggesting that the asymmetry is not fruitful. However, this effect is not observed in regions with high excitation amplitude. The other parameters exhibit similar effects on the mean power, as seen in the previous cases analyzed.

To examine the conditional probability, we divide the analysis into two scenarios: one where the asymmetry is strong, characterized by high values of $\phi$ and $\delta$, and another where the asymmetry is weak, with $\phi$ and $\delta$ close to zero. These situations are detailed in Tab.~\ref{tab:interval_1}.

\begin{figure*}
    \centering
    \subfigure[Domain $\mathcal{D}_{1}$]{\includegraphics[width=0.7\textwidth]{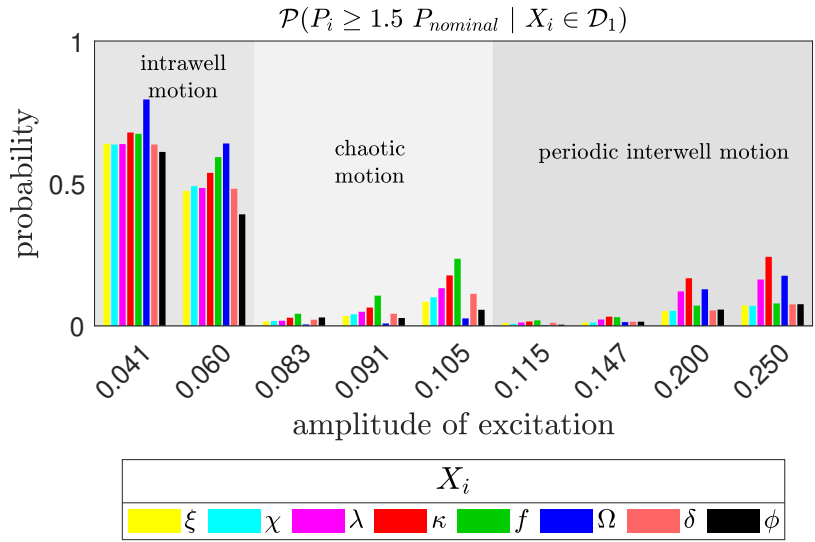}\vspace{-1.5cm} }
    \subfigure[Domain $\mathcal{D}_{2}$]{\includegraphics[width=0.7\textwidth]{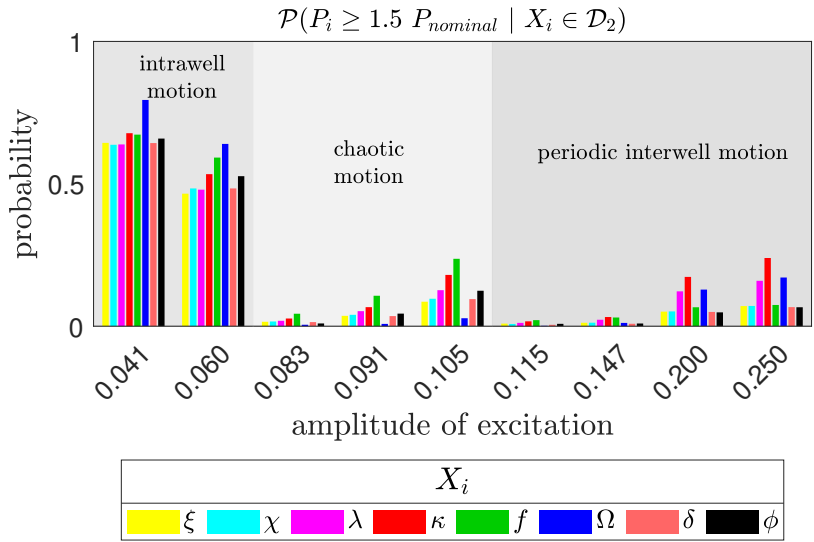}}
    \caption{Probability of increasing the nominal mean power by 50\% as parameter $X_\mathnormal{i}$ is increased by 10\%, plotted against the excitation amplitude, for the asymmetric model with nonlinear piezoelectric coupling.}
    \label{fig:prob_BEHa}
\end{figure*}

Figure~\ref{fig:prob_BEHa}a displays the conditional probability of increasing the mean power by 50\% of its nominal value, given an increase of 10\% in the parameter of interest. The conditional events are displayed in Tab.~\ref{tab:interval_1} as domain $\mathcal{D}_1$. The objective is to assess the influence of asymmetry on power generation, both at the potential coefficient and the angle. In the intrawell motion region, the increase in excitation frequency generates an 70\% probability of increasing the average power, while higher values of the asymmetry coefficient and angle lead to a lower probability of increased power. Especially the angle has a low probability of enhancing the mean power. In the chaos region, the impact of increasing the asymmetry remains low probability, and increasing the amplitude is the best alternative for generating more energy. Finally, for the interwell motion region, the increase in $\kappa$ provides the highest chance of increased power, while asymmetry terms remain low.

In Figure~\ref{fig:prob_BEHa}b, domain $\mathcal{D}_2$ from Tab.~\ref{tab:interval_1} is used to evaluate when the asymmetric system approximates the symmetrical conditions. The result indicates that the probability of increasing the mean power is higher for $\phi$ and $\delta$ when the system has a low asymmetry condition. However, the effect is not superior to the ones brought by frequency, excitation amplitude, and piezoelectric coupling in the intrawell, chaos, and interwell regions.

\begin{table}[!htb]
    \centering
    \caption{Conditional event for each parameter when the system has a strong level of asymmetry ($\mathcal{D}_1$) and a weak level of asymmetry ($\mathcal{D}_2$).}
        \begin{tabular}{ccc}
        \specialrule{.07em}{.05em}{.05em}\noalign{\smallskip}
        \textbf{} & \textbf{Domain $\mathcal{D}_1$}& \textbf{Domain $\mathcal{D}_2$} \\
        \noalign{\smallskip}\specialrule{.07em}{.05em}{.05em}\noalign{\smallskip}
        $\bar{\mathcal{D}}_{\lambda}$  & $X_{\lambda} \geq 1.1 ~ \bar{X}_{\lambda}$ & $X_{\lambda} \geq 1.1 ~ \bar{X}_{\lambda} $  \\ \\
        \noalign{\smallskip}\specialrule{.04em}{.05em}{.05em}\noalign{\smallskip}
        $\bar{\mathcal{D}}_{\kappa}$   & $X_{\kappa} \geq 1.1 ~ \bar{X}_{\kappa}  $ & $X_{\kappa} \geq 1.1 ~ \bar{X}_{\kappa}    $ \\\\
        \noalign{\smallskip}\specialrule{.04em}{.05em}{.05em}\noalign{\smallskip} 
        $\bar{\mathcal{D}}_{\mathnormal{f}}$  & $X_{\mathnormal{f}} \geq 1.1 ~ \bar{X}_{\mathnormal{f}}$ & $X_{\mathnormal{f}} \geq 1.1 ~ \bar{X}_{\mathnormal{f}}$ \\
        \noalign{\smallskip}\specialrule{.04em}{.05em}{.05em}\noalign{\smallskip}
        $\bar{\mathcal{D}}_{\Omega}$   & $X_{\Omega} \geq 1.1 ~ \bar{X}_{\Omega} $  & $X_{\Omega} \geq 1.1 ~ \bar{X}_{\Omega} $ \\
        \noalign{\smallskip}\specialrule{.04em}{.05em}{.05em}\noalign{\smallskip}
        $\bar{\mathcal{D}}_{\delta}$   & $ \left| X_{\delta}\right| \geq 0.1 $ & $ \left| X_{\delta}\right| \leq 0.1 $ \\
        \noalign{\smallskip}\specialrule{.04em}{.05em}{.05em}\noalign{\smallskip}
        $\bar{\mathcal{D}}_{\phi}$   & $\left|X_{\phi}\right| \geq  10^{\circ} $  & $\left|X_{\phi}\right| \leq  10^{\circ}$ \\
        \noalign{\smallskip}\specialrule{.07em}{.05em}{.05em}
        \end{tabular}
    \label{tab:interval_1} 
\end{table}

\begin{figure*}
    \centering
    \subfigure[$\mathnormal{f}=0.041$ for $\lambda$]{\includegraphics[width=0.25\textwidth]{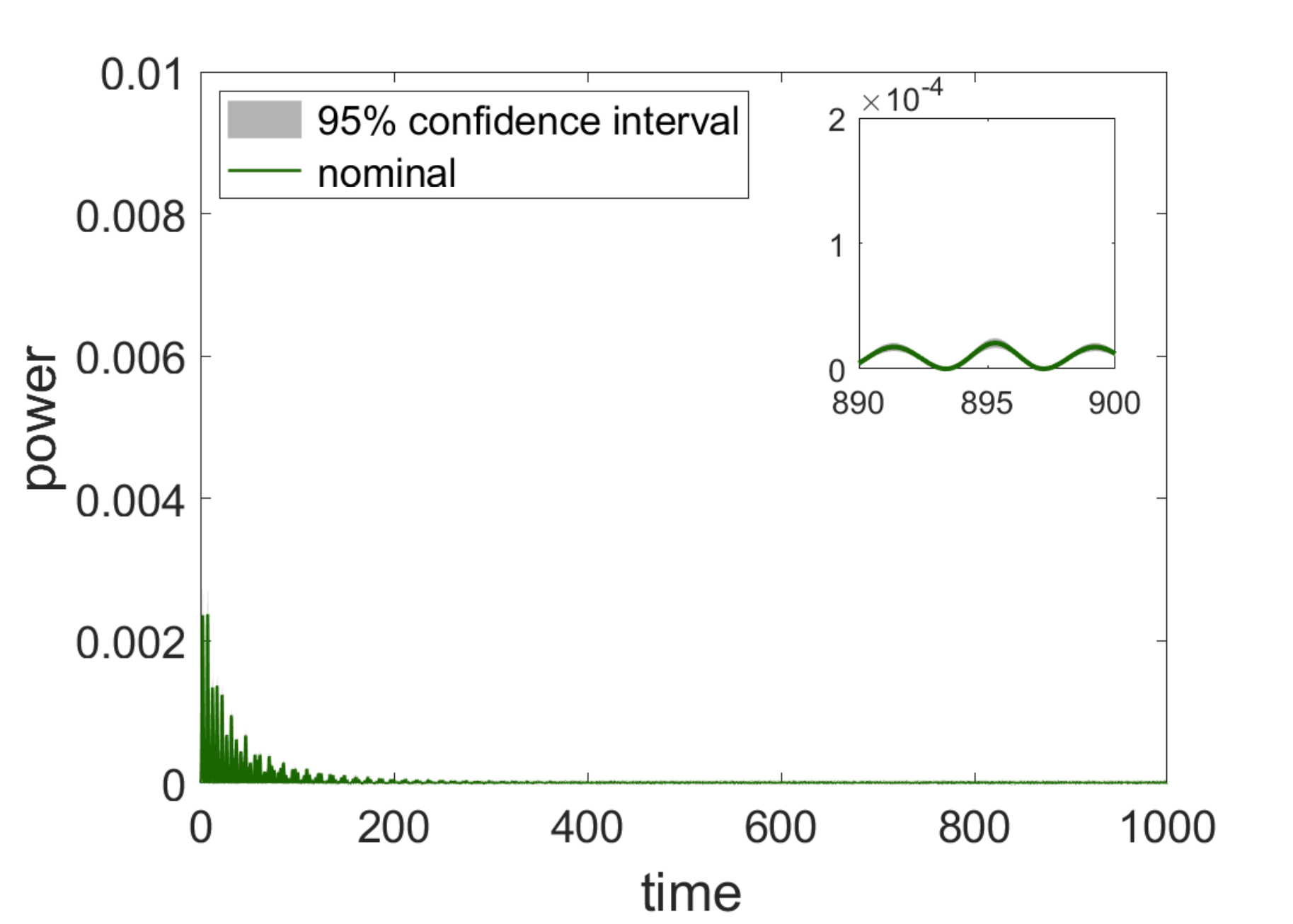}}
    \subfigure[$\mathnormal{f}=0.091$ for $\lambda$]{\includegraphics[width=0.25\textwidth]{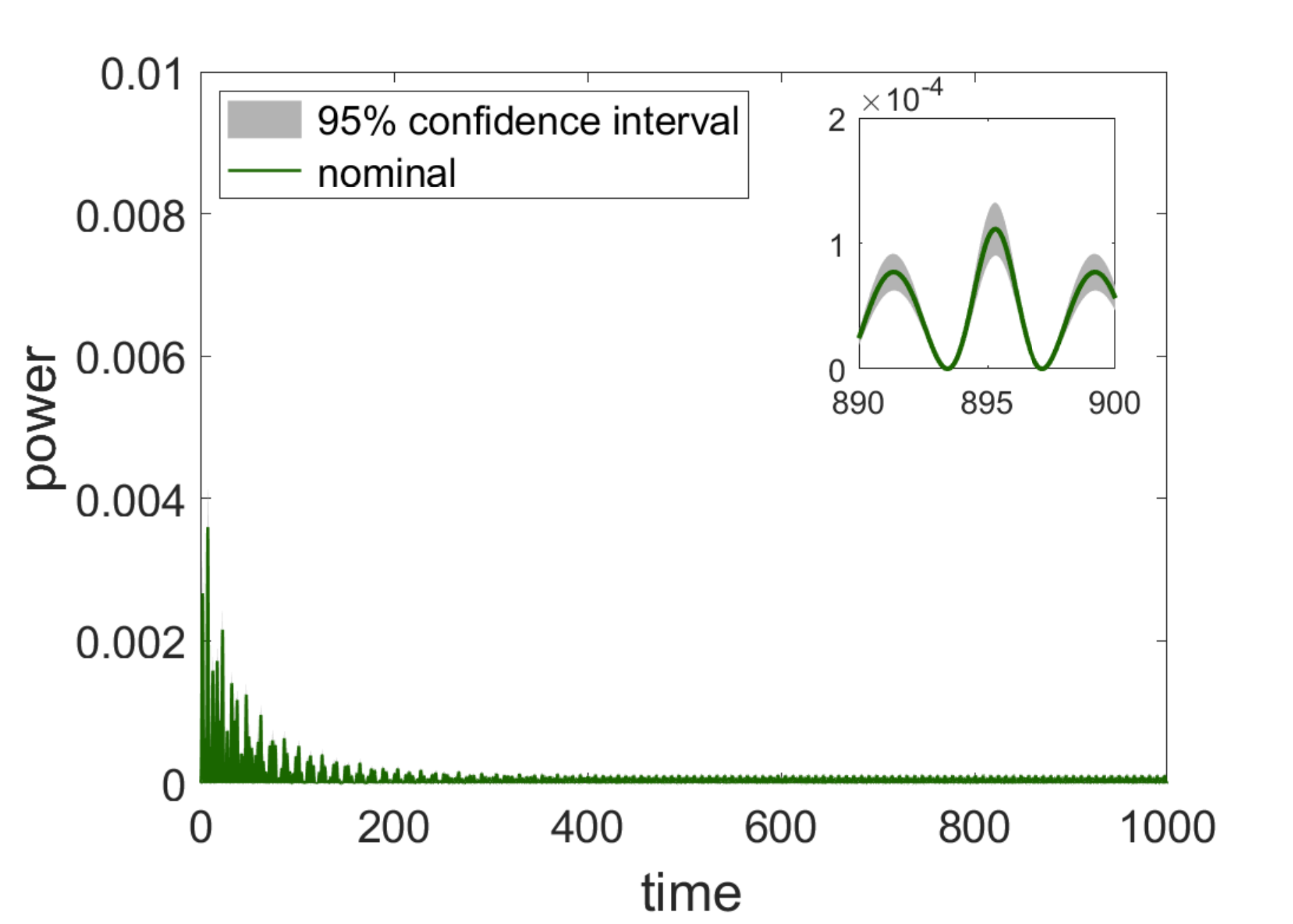}}
    \subfigure[$\mathnormal{f}=0.250$ for $\lambda$]{\includegraphics[width=0.25\textwidth]{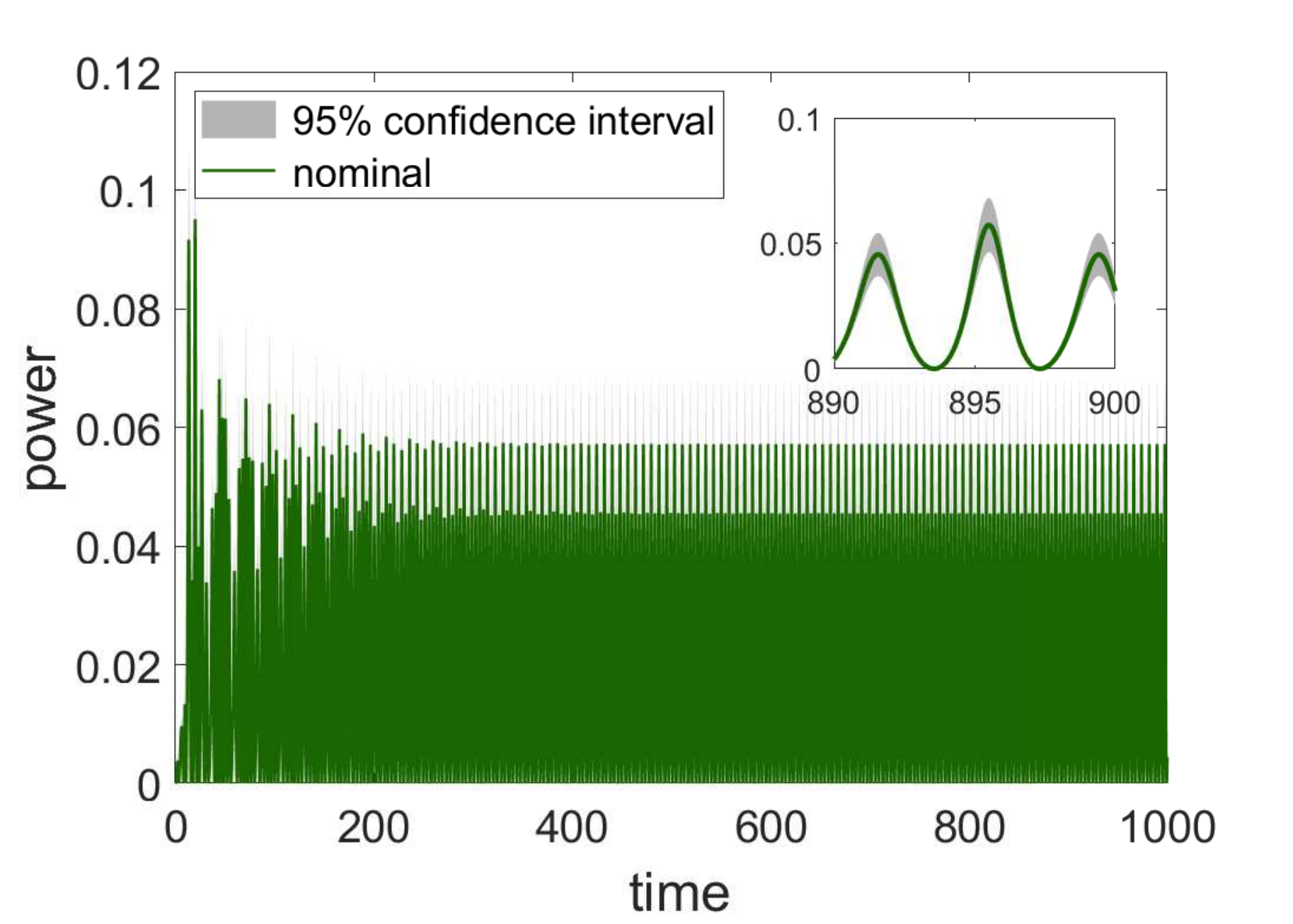}}
    \subfigure[$\mathnormal{f}=0.041$ for $\kappa$]{\includegraphics[width=0.25\textwidth]{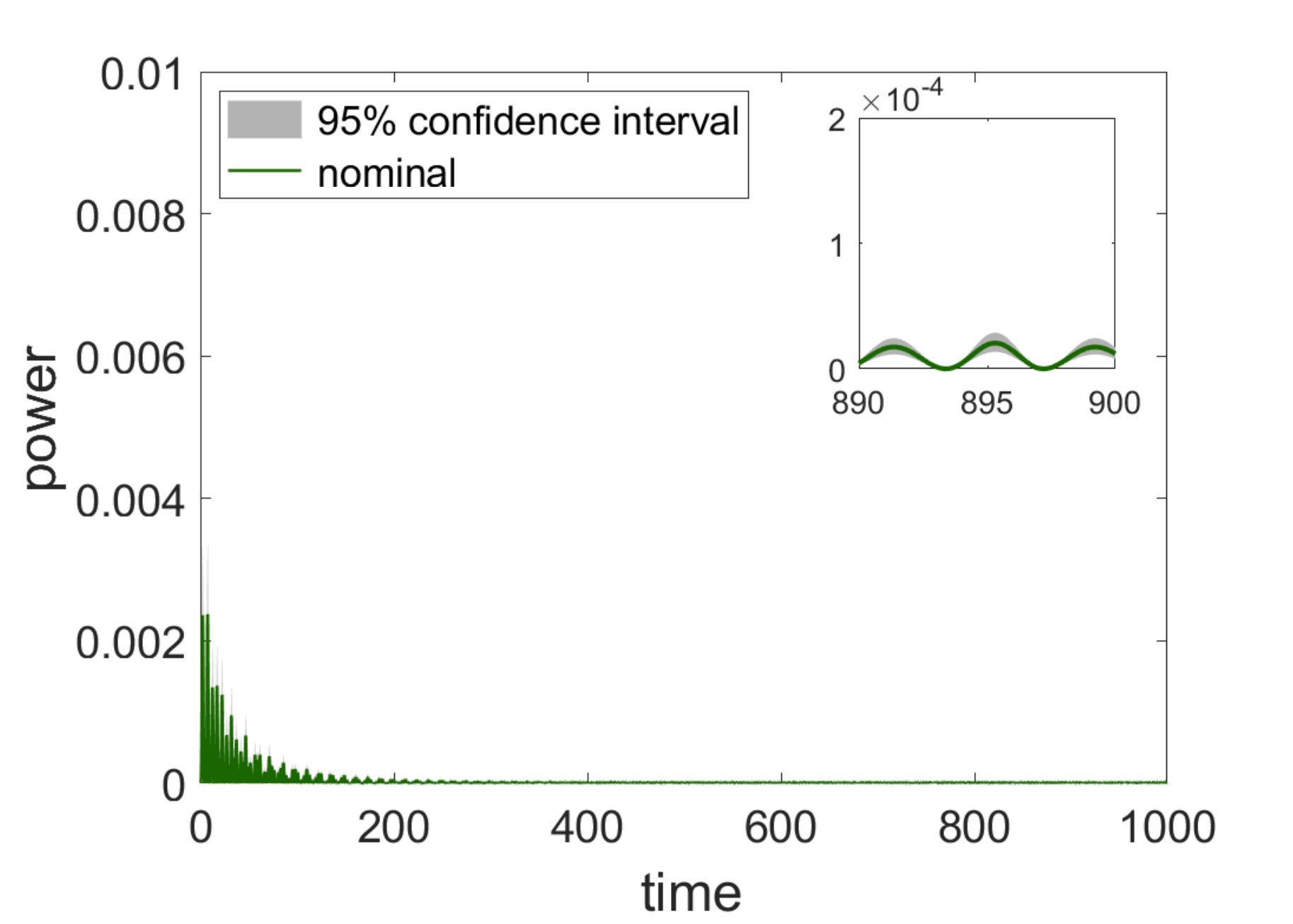}}
    \subfigure[$\mathnormal{f}=0.091$ for $\kappa$]{\includegraphics[width=0.25\textwidth]{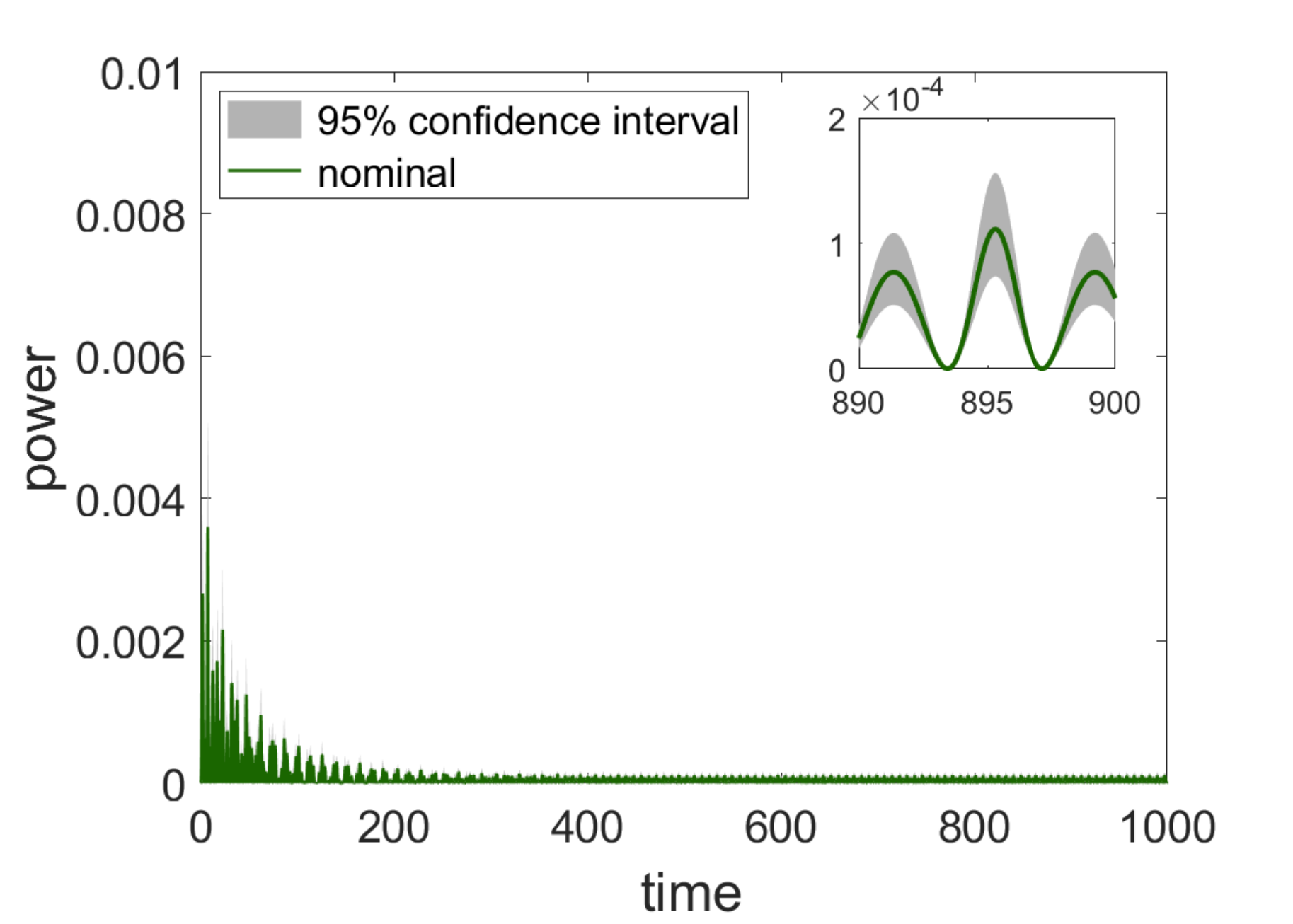}}
    \subfigure[$\mathnormal{f}=0.250$ for $\kappa$]{\includegraphics[width=0.25\textwidth]{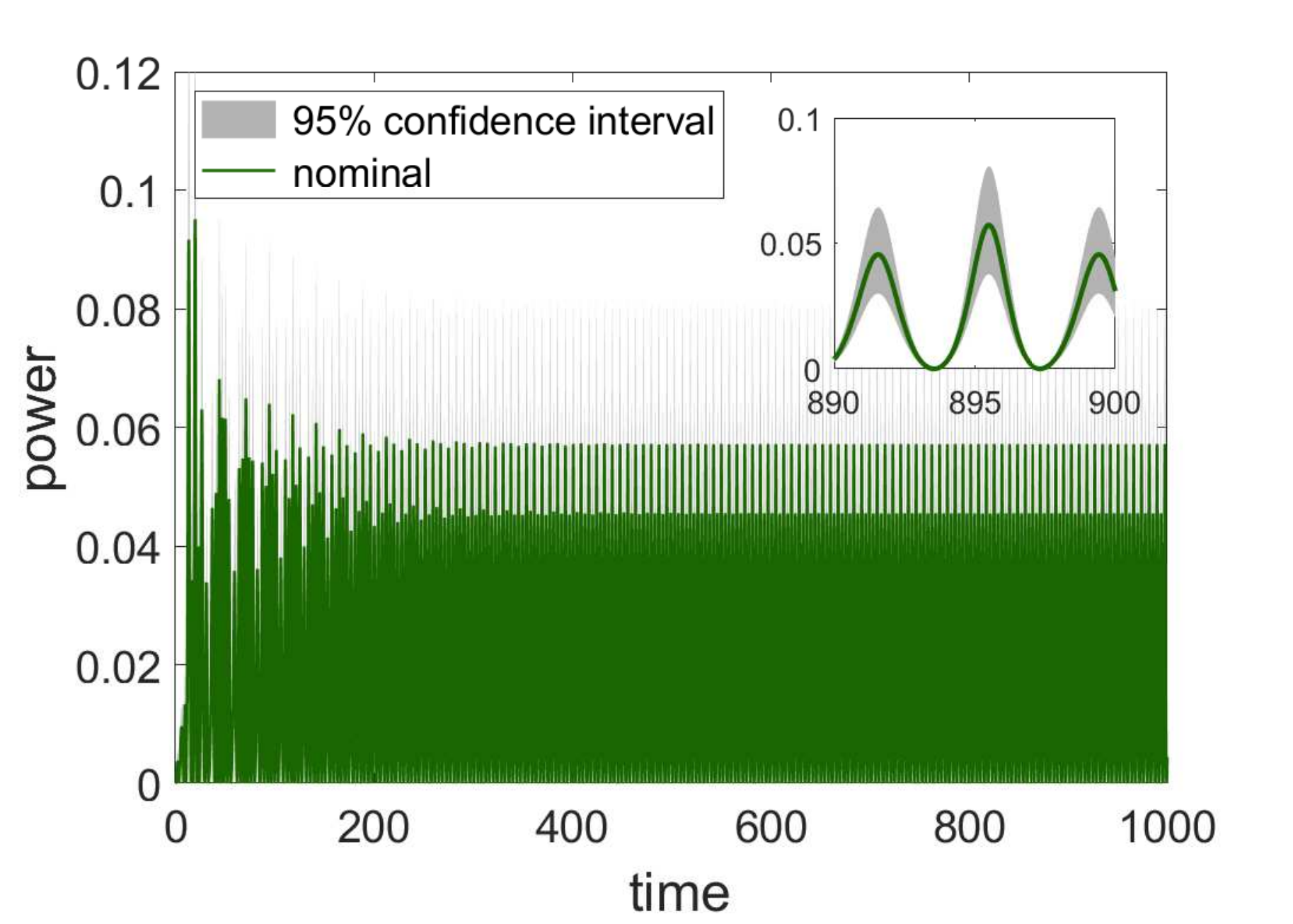}}
    \subfigure[$\mathnormal{f}=0.041$ for $f$]{\includegraphics[width=0.25\textwidth]{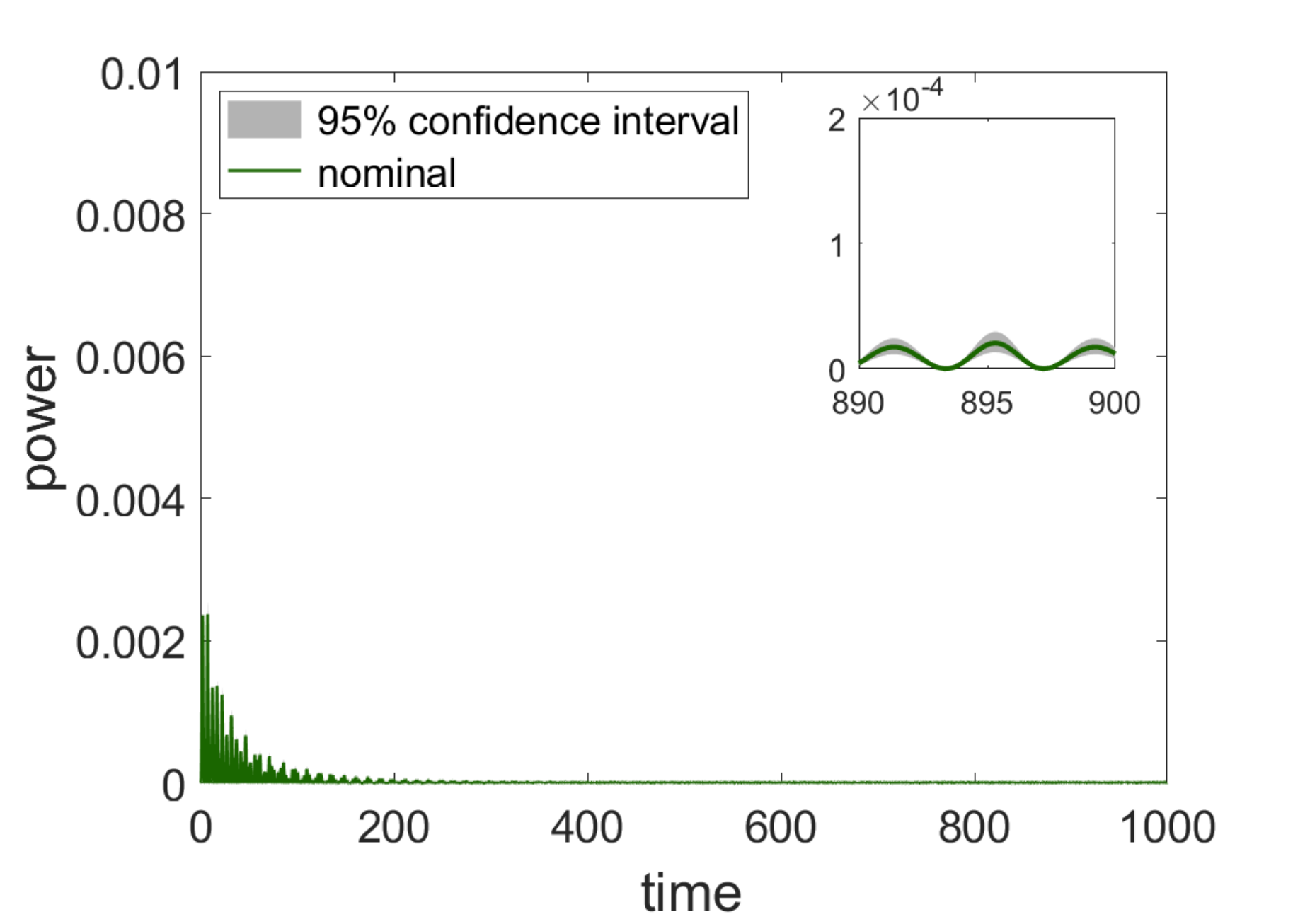}}
    \subfigure[$\mathnormal{f}=0.091$ for $f$]{\includegraphics[width=0.25\textwidth]{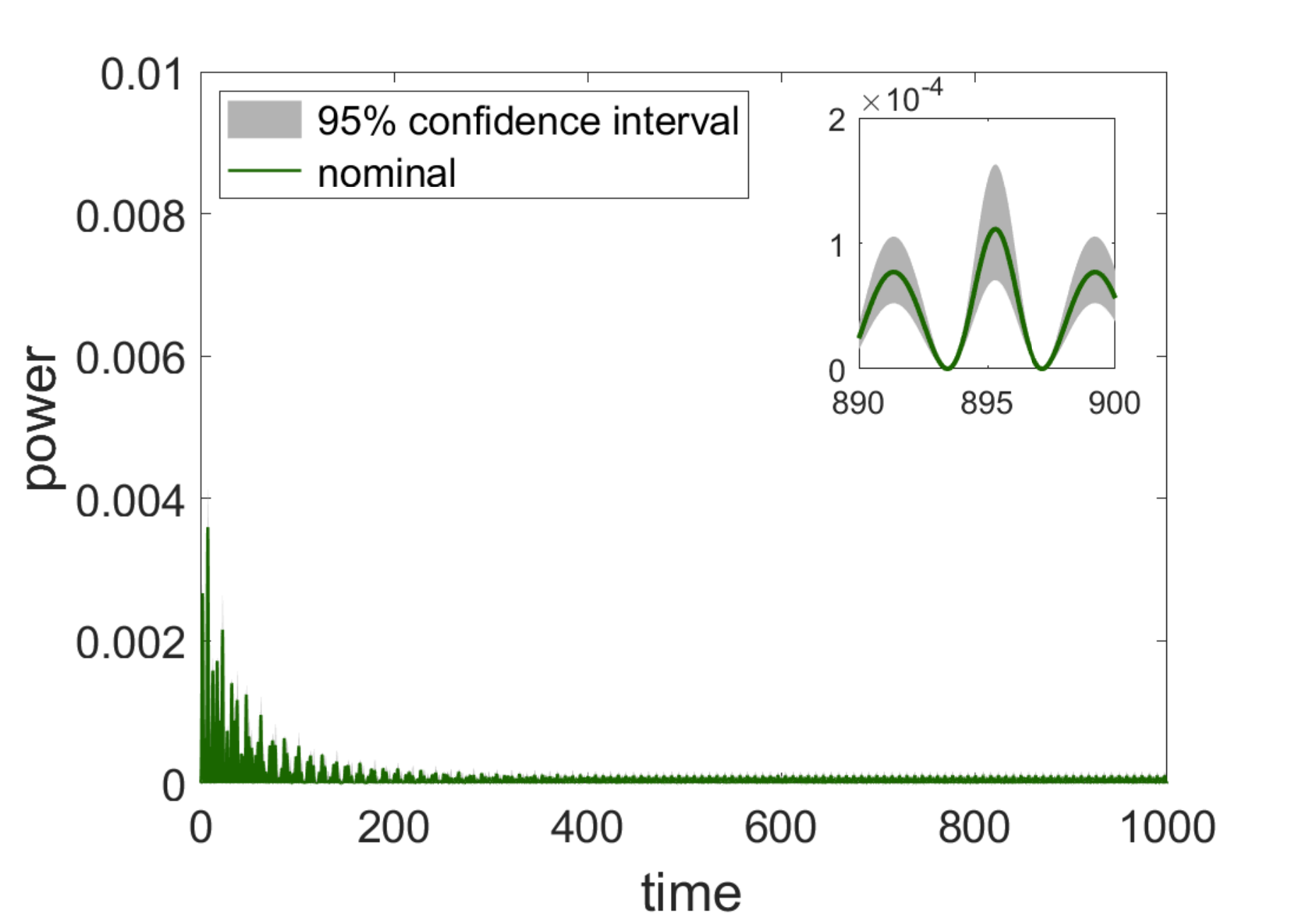}}
    \subfigure[$\mathnormal{f}=0.250$ for $f$]{\includegraphics[width=0.25\textwidth]{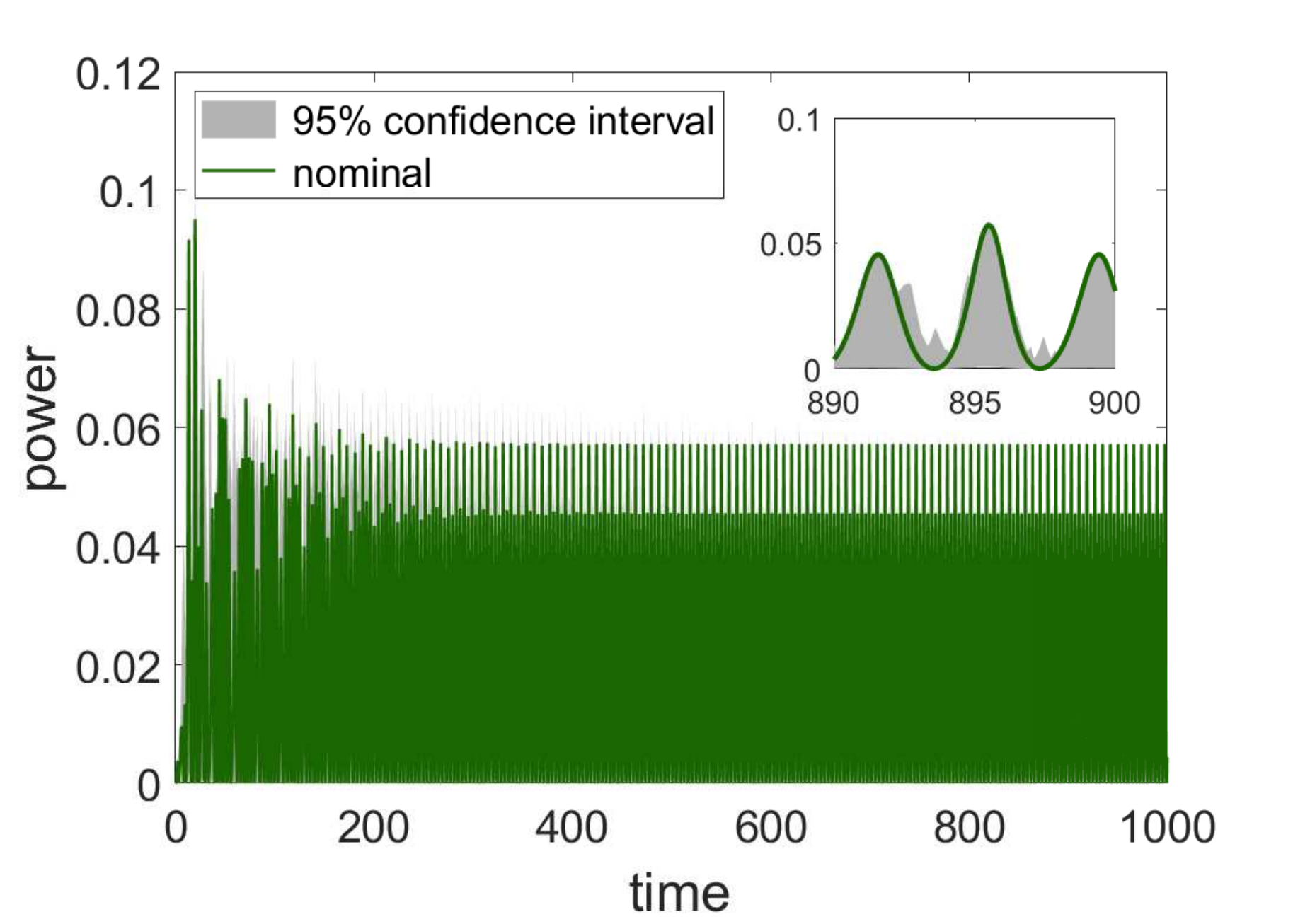}}
    \subfigure[$\mathnormal{f}=0.041$ for $\Omega$]{\includegraphics[width=0.25\textwidth]{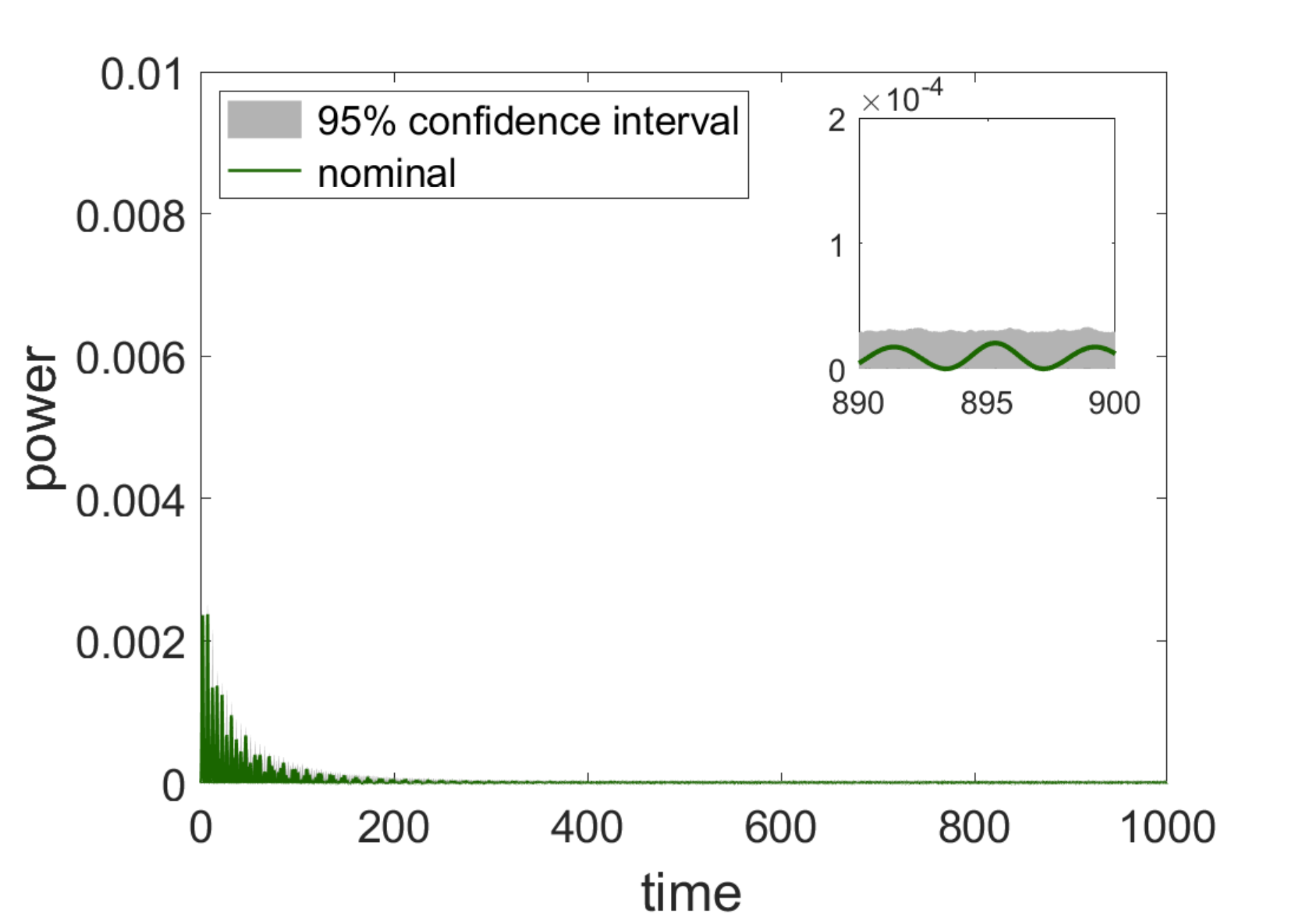}}
    \subfigure[$\mathnormal{f}=0.091$ for $\Omega$]{\includegraphics[width=0.25\textwidth]{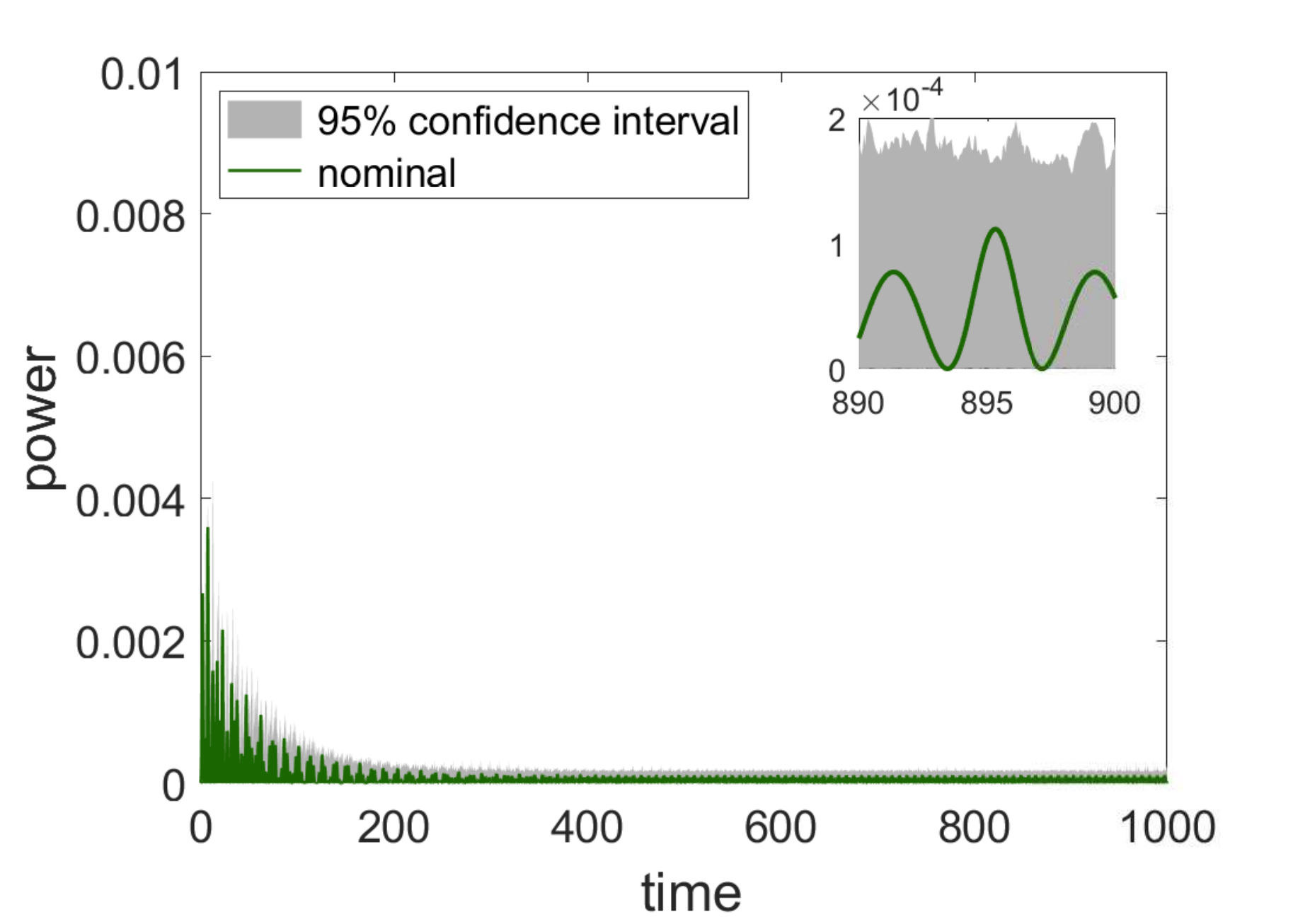}}
    \subfigure[$\mathnormal{f}=0.250$ for $\Omega$]{\includegraphics[width=0.25\textwidth]{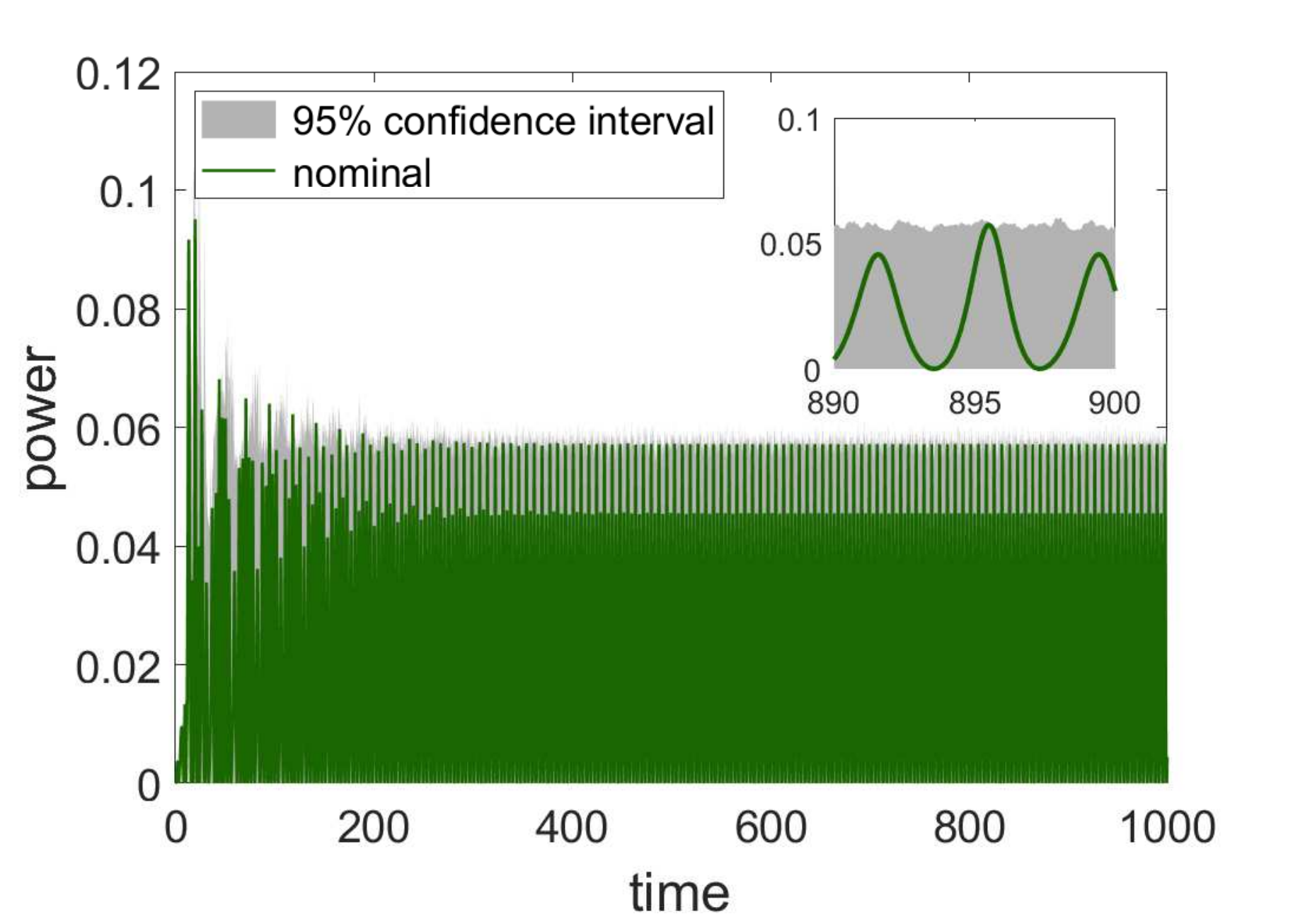}}
    \subfigure[$\mathnormal{f}=0.041$ for $\delta$]{\includegraphics[width=0.25\textwidth]{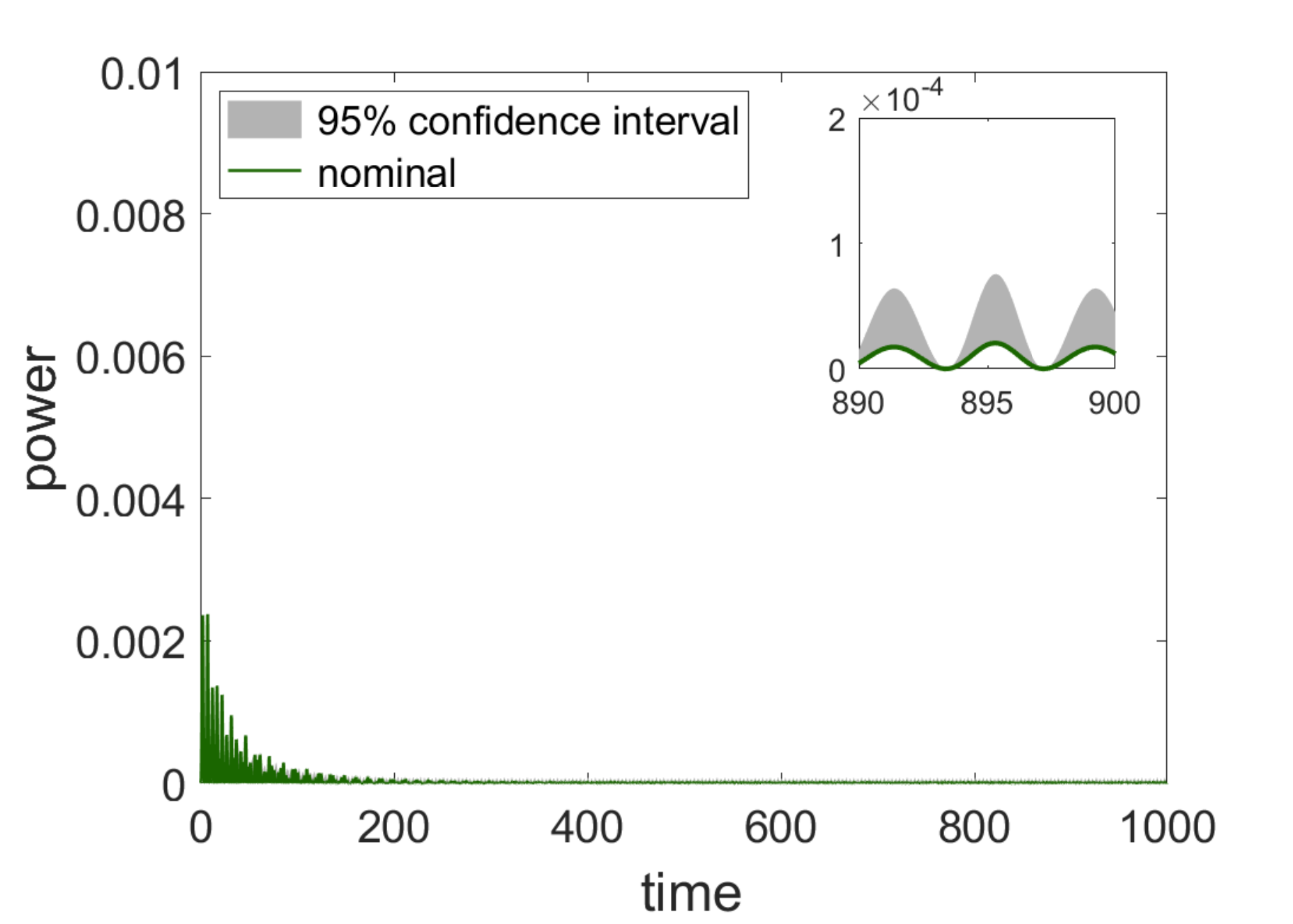}}
    \subfigure[$\mathnormal{f}=0.091$ for $\delta$]{\includegraphics[width=0.25\textwidth]{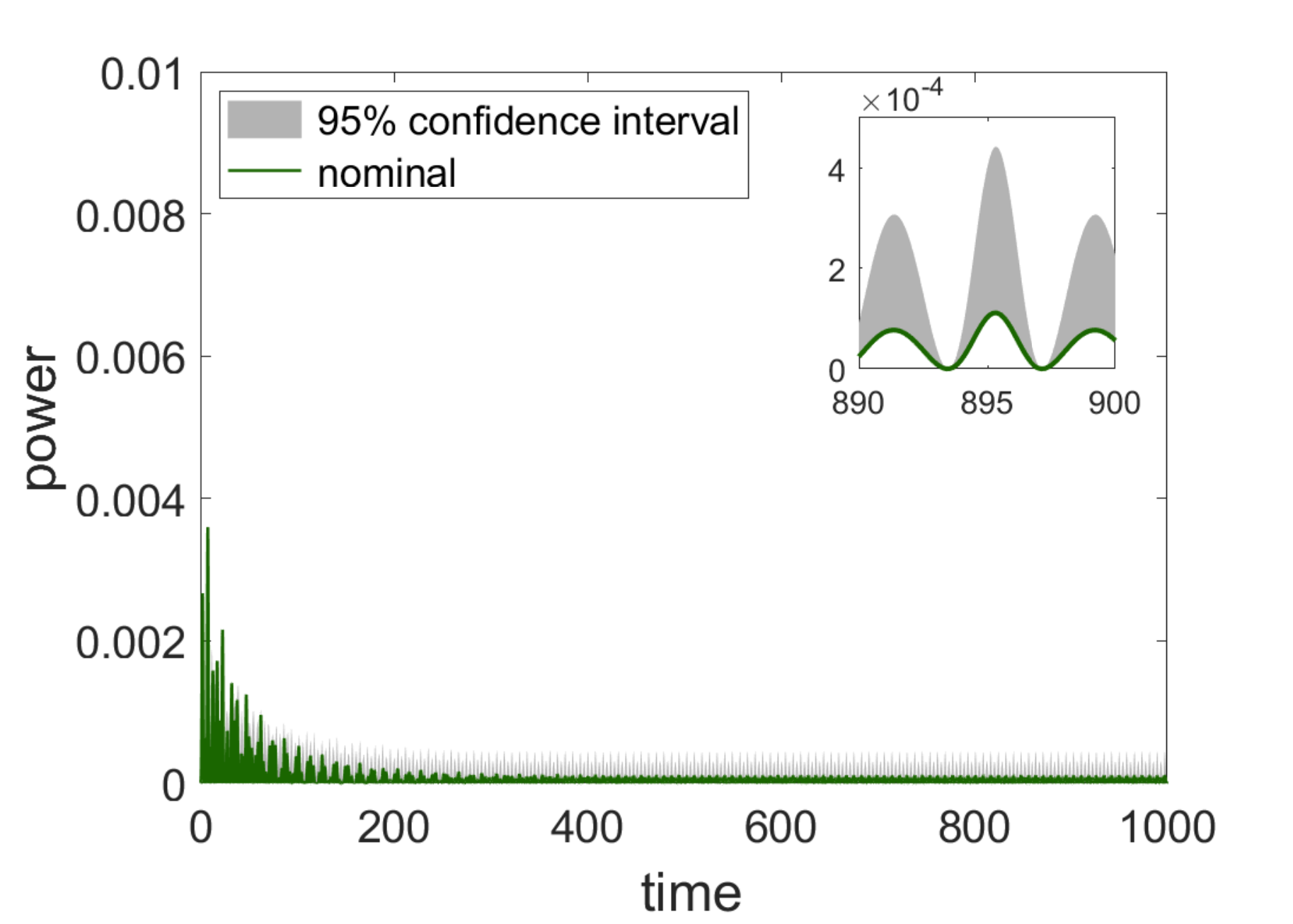}}
    \subfigure[$\mathnormal{f}=0.250$ for $\delta$]{\includegraphics[width=0.25\textwidth]{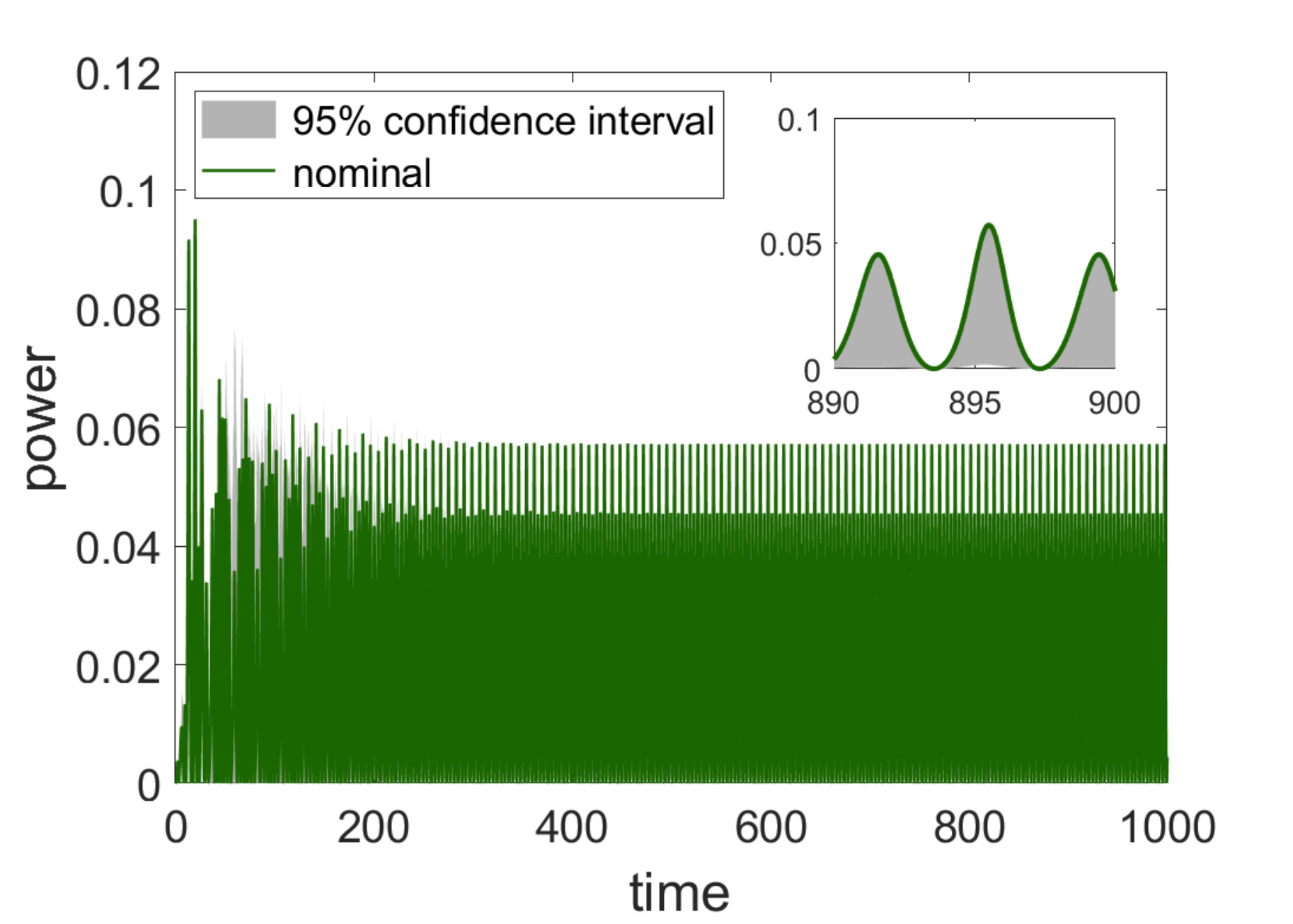}}
    \subfigure[$\mathnormal{f}=0.041$ for $\phi$]{\includegraphics[width=0.25\textwidth]{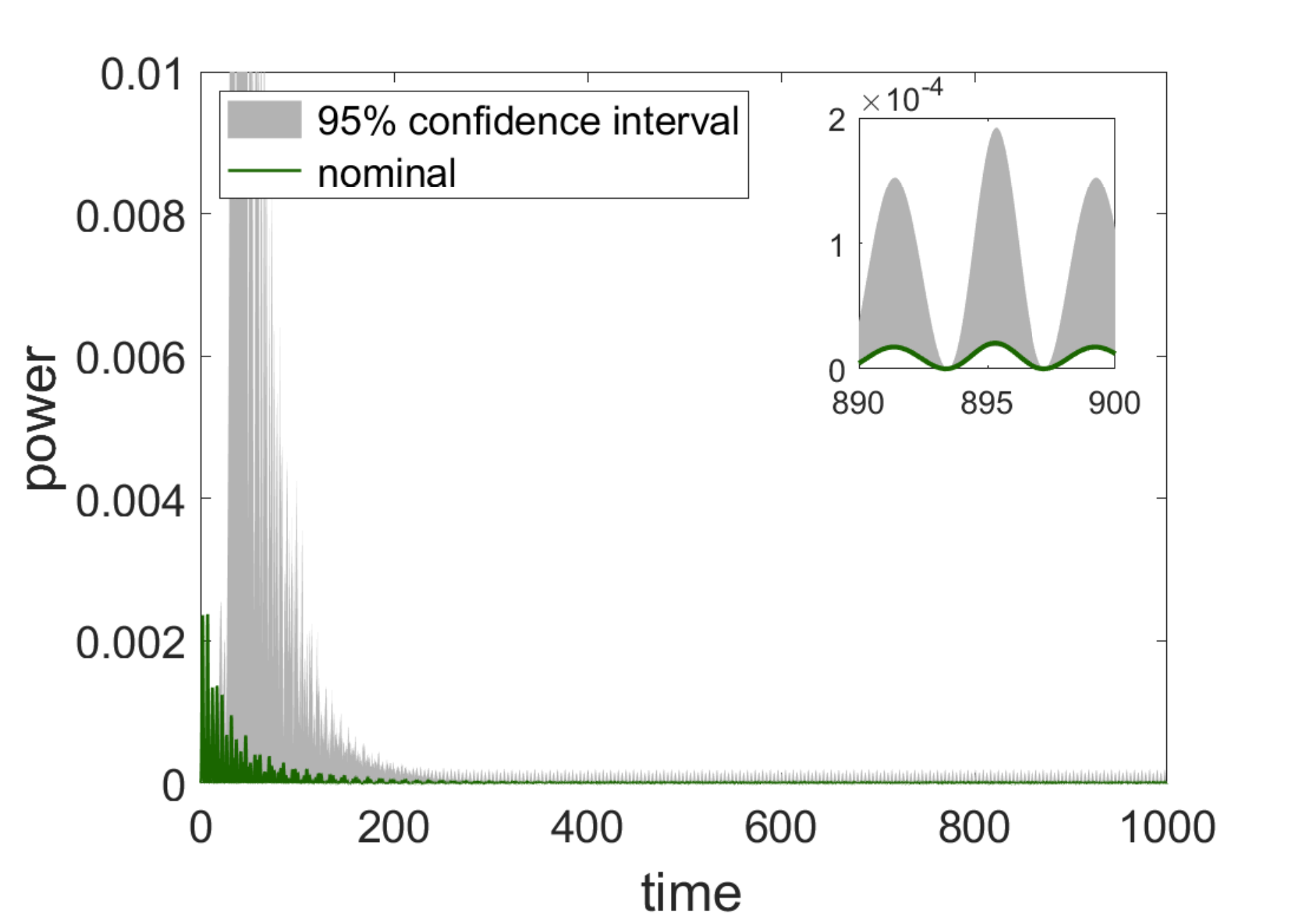}}
    \subfigure[$\mathnormal{f}=0.091$ for $\phi$]{\includegraphics[width=0.25\textwidth]{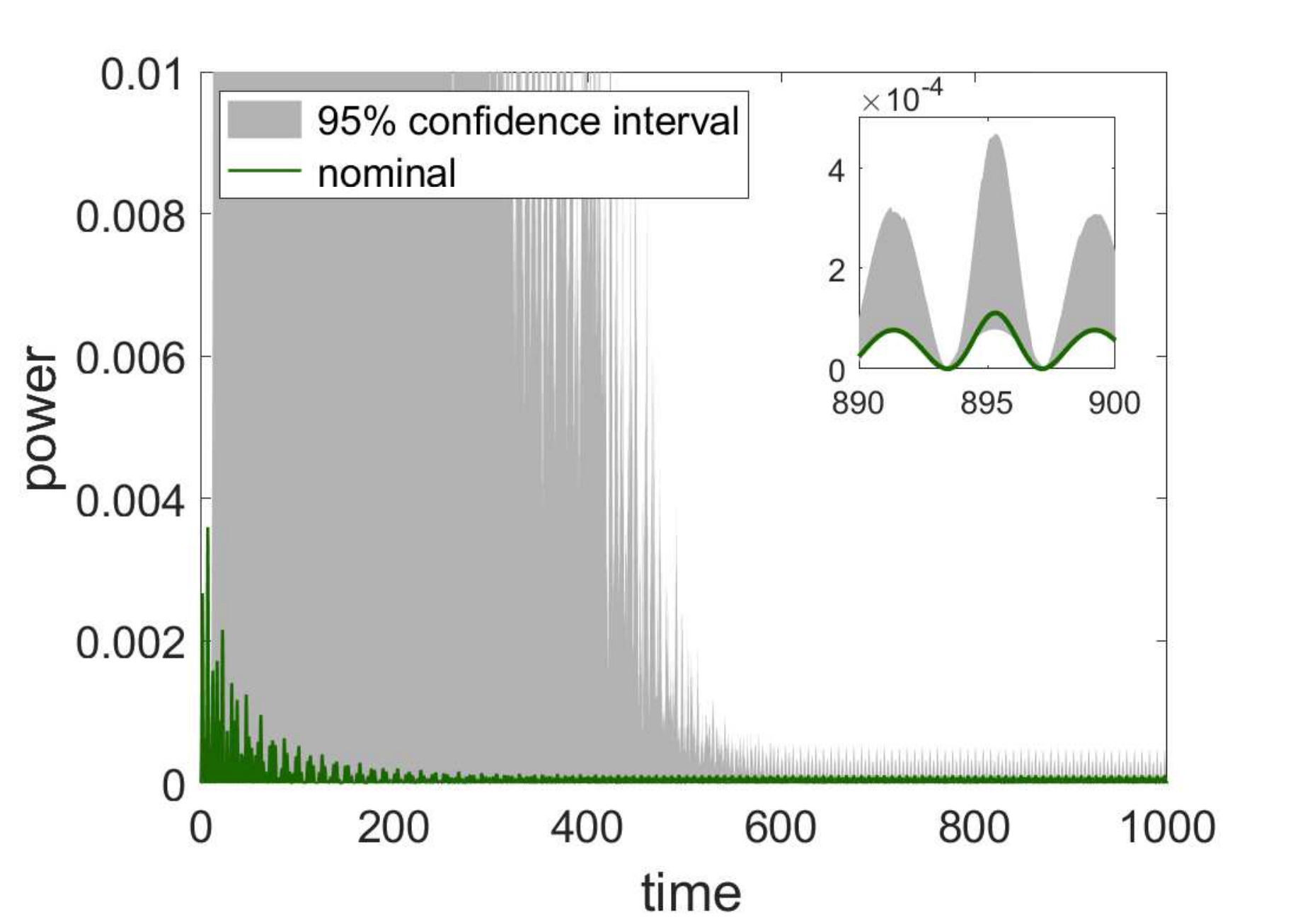}}
    \subfigure[$\mathnormal{f}=0.250$ for $\phi$]{\includegraphics[width=0.25\textwidth]{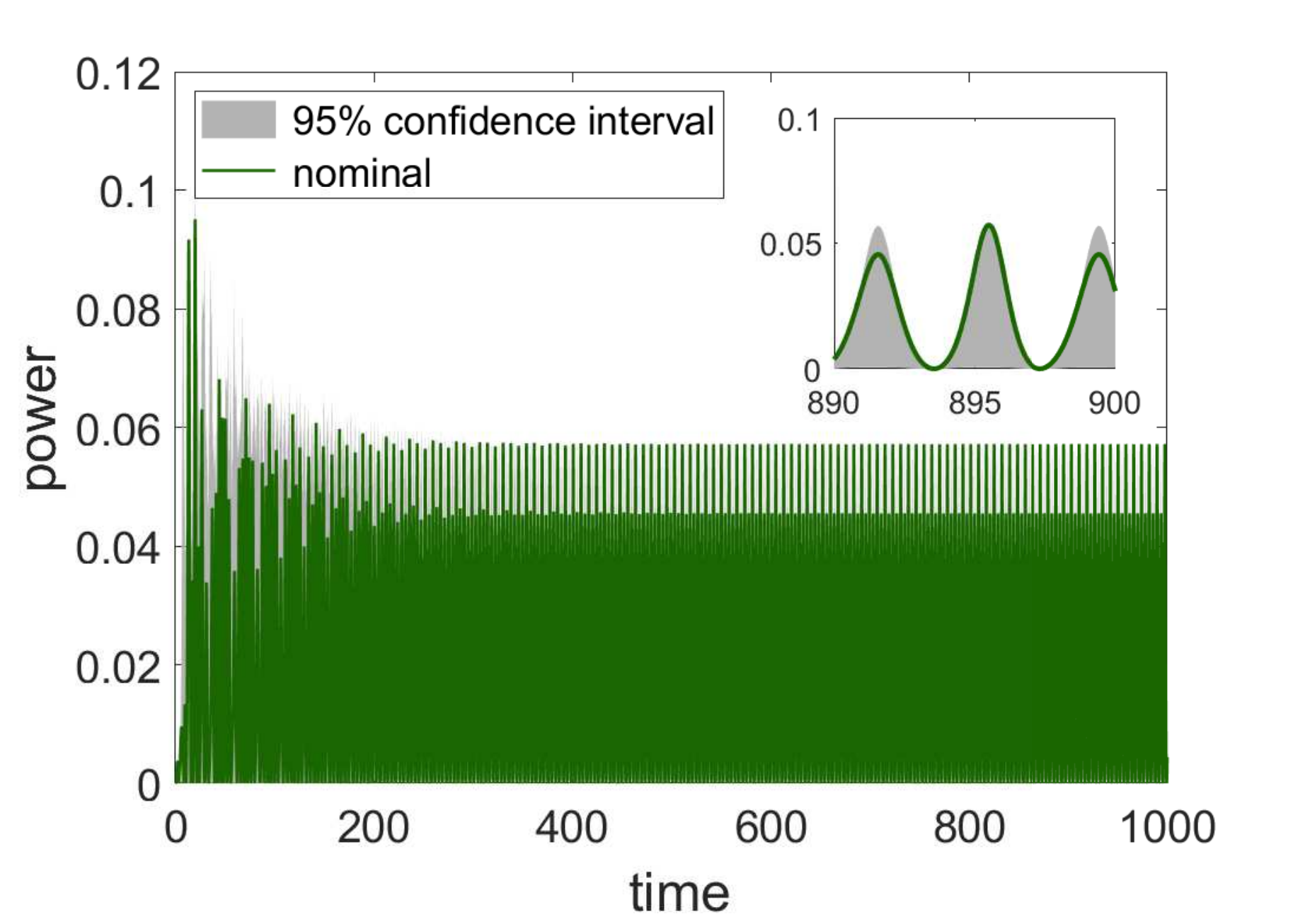}}
    \caption{Propagated uncertainty in the output power time series of the asymmetric model with nonlinear piezoelectric coupling is shown under individual parameters: $\lambda$ (first row), $\kappa$ (second row), $\mathnormal{f}$ (third row), $\Omega$ (fourth row), $\Omega$ (fifth row), and $\phi$ (sixth row). The columns are divided according to the different motion states of the system: intrawell (left), chaos (middle), and interwell (right).}
    \label{fig:up_BEHa}
\end{figure*}

Figure~\ref{fig:up_BEHa} shows the uncertainty propagation for the output power over time for each parameter individually, considering the 95\% confidence interval and the nominal series. In the left column for intrawell motion, the parameters $\lambda$, $\kappa$, $\mathnormal{f}$, and $\Omega$ behave similarly to the previous model. The asymmetric parameters, $\delta$, and $\phi$, generate a significant confidence interval, especially the latter. They also suggest a high amplitude of harvested power, as expected for conditions with small values of $\phi$, as visualized in Fig.~\ref{fig:prob_BEHa}. At the chaotic motion, all parameters drastically alter the confidence interval, with the transient regime for $\phi$ variations exhibiting the most significant effect. The asymmetries again indicate a substantial interval of confidence. In the right column for interwell motion, $\lambda$, $\kappa$, $\mathnormal{f}$, and $\Omega$ affect the mean power in the same way as the previous models at the same condition. However, $\delta$ and $\phi$ generate lower energy harvesting values, demonstrating undesirable behavior, as shown in Figure~\ref{fig:prob_BEHa}.

% --------------------------------------------------------------

\subsection{Importance of probabilistic maps insights}

To highlight the significance of studying conditional probabilities and emphasize the findings of this research, this section provides a concise comparison of results. Specifically, it examines the scenario where all variables vary simultaneously versus the case where only one variable varies while the others remain constant. For this analysis, we employ the nonlinear harvesting system investigated in \cite{PANYAM2014153}. While that study offered valuable insights into the influence of system parameters using perturbation techniques, it primarily focused on individual variations of the parameters. Consequently, to illustrate the complexity of this situation and demonstrate the potential divergence in conclusions, we calculate the resulting variations in harvested power for both joint and individual parameter changes.

\begin{figure*}
    \centering
    \includegraphics[width=0.7\textwidth]{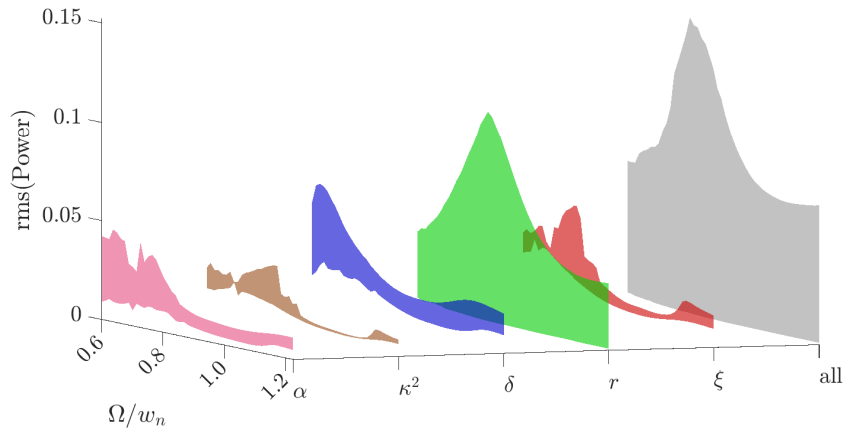}
    \caption{Confidence band of 95\% of the root mean square of the recovered power as a function of the excitation frequency ratio when the excitation force is 0.1. The model used to derive these numerical results is described in detail in the work \cite{PANYAM2014153}. The magenta band is the response when the parameter $\alpha$ (ratio of time constant) is varied uniformly in the range $[0.025~~0.075]$ while the other is fixed in their nominal values. The brown band corresponds to parameter $\kappa^2$ (linear electromechanical coupling) varying between $[0.01~~0.2]$. The blue band represents parameter $\delta$ (coefficient of cubic nonlinearity) in the range $[0.25~~0.75]$. The green band represents parameter $r$ (potential shape coefficient) ranging from $ [1.1~~2]$. The red band depicts parameter $\xi$ (damping ratio) varying between $[0.025~~0.075]$. Lastly, the gray band shows the combined effect of all parameters simultaneously varying.}
    \label{fig:comparison}
\end{figure*}

Figure~\ref{fig:comparison} illustrates the uncertainty confidence band of the root mean square of the power as a function of the excitation frequency ratio. The gray band represents the power response of the system when all variables are considered stochastic, considering joint variations. The other bands depict the power response when only the indicated parameter is varying, representing the individual variations.
The confidence bands for the parameters $\alpha$, $\kappa^2$, $\delta$, and $\xi$ exhibit a similar pattern. They indicate a comprehensive band in low-frequency regions and narrower at higher frequencies. On the other hand, the confidence band for the parameter $r$ is wide across all frequencies, with a noticeable peak around $\Omega/\omega_n$ equal to 1.2. This behavior resembles the gray confidence band, where all parameters vary simultaneously.
This similarity suggests a high sensitivity of power to the parameter $r$. But comparing the solution region, the joint effect results in a larger domain than the green band.  This suggests that the variation of all parameters, when combined, also exerts a significant influence. It further implies that the sensitivity of joint parameters can be pronounced. In other words, when all parameters are varying, the system becomes more sensitive to their combined effects.
Considering this effect, the calculation of conditional probabilities and the formulation of probability maps presented in this work prove fruitful in demonstrating the influence of each parameter when the other parameters are also varying.

% rev by Americo
% --------------------------------------------------------------
\section{Conclusions}
\label{concl_sect}

This paper presented a comprehensive uncertainty quantification analysis of the most sensitive parameters on several bistable energy harvesting systems. The analysis considered asymmetries and nonlinear piezoelectric coupling, and determined the probability density function of mean power for several excitation situations. The most sensitive parameters were treated as uniformly distributed random variables to obtain the joint CDF of mean power conditioned for each parameter of interest, which was then used to check the correlation between them and the mean power. The PDF was estimated using a kernel density estimator, and the conditional probability for increasing the mean power was also calculated.

The analysis results suggest that there is a specific parameter that has a consistently positive impact on the energy harvested, regardless of the model being considered. When it comes to intrawell motion, increasing the excitation frequency increases the probability of generating more power. On the other hand, in chaotic regimes, increasing the excitation amplitude is the preferred option for enhancing energy harvesting. For interwell motion, increasing the piezoelectric coupling offers a higher chance of increasing the power. However, in asymmetric systems, there are exceptions where changes in the plane angle can be favorable or unfavorable to the harvesting performance. Notably, weaker levels of asymmetry are advantageous, while stronger asymmetry values can be harmful. The configuration of asymmetry and excitation conditions plays a crucial role in determining performance, highlighting the complexity of nonlinearity. Additionally, some parameters exhibit heightened influence during the transient regime, despite having low sensitivity. Lastly, the uncertainties associated with the most sensitive parameters result in a confidence interval solely in the magnitude of output power, except for the excitation frequency.

This work presents important information about the uncertainty influence in the nonlinear energy harvesting processes and describes when and how each parameter's uncertainty can impact the harvesting process. The findings of this work are consistent with the sensitivity study reported in \cite{norenbergNoDy_2022}, demonstrating that the most sensitive parameters can increase or decrease the electrical power. In summary, the contributions of this work provide valuable insight into the uncertainties of nonlinear energy harvesting systems, which can guide the design and optimization of these systems.
% --------------------------------------------------------------

% rev by Americo
% --------------------------------------------------------------
\section*{Acknowledgements}

The authors gratefully acknowledge the insightful discussions on the results presented in this paper with Professors Grzegorz Litak (Lublin University of Technology) and Marcelo Savi (Federal University of Rio de Janeiro).
% --------------------------------------------------------------

% rev by Americo
% --------------------------------------------------------------
\section*{Funding}

This research was financially supported by the Brazilian agencies Coordena\c{c}\~{a}o de Aperfei\c{c}oamento de Pessoal de N\'{\i}vel Superior (CAPES) under Finance Code 001, Conselho Nacional de Desenvolvimento Cient\'{i}fico e Tecnol\'{o}gico under the grants 306526 /2019-0 and 305476 /2022-0, and the Carlos Chagas Filho Research Foundation of Rio de Janeiro State (FAPERJ) under grants 210.167/2019, 211.037/2019, and 201.294/2021.
% --------------------------------------------------------------

% rev by Americo
% --------------------------------------------------------------
\section*{Code availability}

The simulations presented in this paper were performed using the computational code {\bf STONEHENGE - Suite for Nonlinear Analysis of Energy Harvesting Systems} \cite{STONEHENGE_paper}, which is available for free on GitHub \cite{STONEHENGE}.
% --------------------------------------------------------------

% --------------------------------------------------------------
\section*{Declarations}
% --------------------------------------------------------------

% rev by Americo
% --------------------------------------------------------------
\section*{Conflict of Interest }
The authors declare they have no conflict of interest.
% --------------------------------------------------------------

% --------------------------------------------------------------
\section*{Disclaimer}

This manuscript has undergone a comprehensive revision utilizing artificial intelligence-powered tools, such as Grammarly and ChatGPT, to enhance its grammatical accuracy and clarity. However, the authors assume full responsibility for the original language and phrasing used in the manuscript.
% --------------------------------------------------------------

% rev by Americo
\bibliographystyle{spmpsci}      % mathematics and physical sciences

\bibliography{references}

\begin{thebibliography}{10}
\providecommand{\url}[1]{{#1}}
\providecommand{\urlprefix}{URL }
\expandafter\ifx\csname urlstyle\endcsname\relax
  \providecommand{\doi}[1]{DOI~\discretionary{}{}{}#1}\else
  \providecommand{\doi}{DOI~\discretionary{}{}{}\begingroup
  \urlstyle{rm}\Url}\fi

\bibitem{Ali_2010}
Ali, S.F., Friswell, M.I., Adhikari, S.: Piezoelectric energy harvesting with
  parametric uncertainty.
\newblock Smart Materials and Structures \textbf{19}(10), 105010 (2010).
\newblock \doi{10.1088/0964-1726/19/10/105010}

\bibitem{Chatterjee_2022}
Chatterjee, T., Karlicic, D., Adhikari, S., Friswell, M.: Parametric
  amplification in a stochastic nonlinear piezoelectric energy harvester via
  machine learning.
\newblock In: Proceedings of the Society for Experimental Mechanics Series
  (2022)

\bibitem{Lua_2020}
Costa, L.G., da~Silva~Monteiro, L.L., Pacheco, P.M.C.L., Savi, M.A.: A
  parametric analysis of the nonlinear dynamics of bistable vibration-based
  piezoelectric energy harvesters.
\newblock Journal of Intelligent Material Systems and Structures
  \textbf{32}(7), 699--723 (2021).
\newblock \doi{10.1177/1045389X20963188}

\bibitem{cottone2009p080601}
Cottone, F., Vocca, H., Gammaitoni, L.: Nonlinear energy harvesting.
\newblock Phys. Rev. Lett. \textbf{102}, 080601 (2009).
\newblock \doi{10.1103/PhysRevLett.102.080601}

\bibitem{37_Crestaux}
Crestaux, T., {Le Ma\^{\i}tre}, O., Martinez, J.M.: Polynomial chaos expansion
  for sensitivity analysis.
\newblock Reliability Engineering \& System Safety \textbf{94}(7), 1161--1172
  (2009).
\newblock \doi{10.1016/j.ress.2008.10.008}.
\newblock Special Issue on Sensitivity Analysis

\bibitem{Cunha_uq2017}
Cunha, A.: Modeling and Quantification of Physical Systems Uncertainties in a
  Probabilistic Framework, pp. 127--156.
\newblock Springer International Publishing, Cham (2017).
\newblock \doi{10.1007/978-3-319-55852-3_8}

\bibitem{cunhajr2021p137}
{Cunha~Jr}, A.: Enhancing the performance of a bistable energy harvesting
  device via the cross-entropy method.
\newblock Nonlinear Dynamics \textbf{103}, 137--155 (2021).
\newblock \doi{10.1007/s11071-020-06109-0}

\bibitem{cunhajr2014p1355}
{Cunha~Jr}, A., Nasser, R., Sampaio, R., Lopes, H., Breitman, K.: Uncertainty
  quantification through {M}onte {C}arlo method in a cloud computing setting.
\newblock Computer Physics Communications \textbf{185}, 1355 – 1363 (2014).
\newblock \doi{10.1016/j.cpc.2014.01.006}

\bibitem{19_dAQAQ}
Daqaq, M., Masana, R., Erturk, A., Quinn, D.: On the role of nonlinearities in
  vibratory energy harvesting: A critical review and discussion.
\newblock Applied Mechanics Reviews \textbf{66}, 040801 (2014).
\newblock \doi{10.1115/1.4026278}

\bibitem{12_duToit}
duToit, N.E., Wardle, B.L.: Experimental verification of models for
  microfabricated piezoelectric vibration energy harvesters.
\newblock AIAA Journal \textbf{45}(5), 1126--1137 (2007).
\newblock \doi{10.2514/1.25047}

\bibitem{erturk2009p254102}
Erturk, A., Hoffmann, J., Inman, D.J.: A piezomagnetoelastic structure for
  broadband vibration energy harvesting.
\newblock Applied Physics Letters \textbf{94}(25), 254102 (2009).
\newblock \doi{10.1063/1.3159815}

\bibitem{ghanem2003}
Ghanem, R.G., Spanos, P.D.: Stochastic Finite Elements: A Spectral Approach,
  2nd edn.
\newblock Dover~Publications, New~York (2003)

\bibitem{Godoy_2012}
Godoy, T., Trindade, M.: Effect of parametric uncertainties on the performance
  of a piezoelectric energy harvesting device.
\newblock Journal of the Brazilian Society of Mechanical Sciences and
  Engineering \textbf{34}, 552--560 (2012).
\newblock \doi{10.1590/S1678-58782012000600003}

\bibitem{Huang_2020}
Huang, D., Zhou, S., Han, Q., Litak, G.: Response analysis of the nonlinear
  vibration energy harvester with an uncertain parameter.
\newblock Proceedings of the Institution of Mechanical Engineers, Part K:
  Journal of Multi-body Dynamics \textbf{234}(2), 393--407 (2020).
\newblock \doi{10.1177/1464419319893211}

\bibitem{Jia2020}
Jia, Y.: Review of nonlinear vibration energy harvesting: Duffing, bistability,
  parametric, stochastic and others.
\newblock Journal of Intelligent Material Systems and Structures
  \textbf{31}(7), 921--944 (2020).
\newblock \doi{10.1177/1045389X20905989}

\bibitem{Kapur_1992}
Kapur, J., Kesavan, H.: Entropy optimization principles with applications.
\newblock Academic Press, Cambridge (1992)

\bibitem{Khovanov_2021}
Khovanov, I.A.: The response of a bistable energy harvester to different
  excitations: the harvesting efficiency and links with stochastic and
  vibrational resonances.
\newblock Philosophical Transactions of the Royal Society A: Mathematical,
  Physical and Engineering Sciences \textbf{379}(2198), 20200245 (2021).
\newblock \doi{10.1098/rsta.2020.0245}

\bibitem{Koka_2014}
Koka, A., Zhou, Z., Tang, H., Sodano, H.A.: Controlled synthesis of ultra-long
  vertically aligned batio3 nanowire arrays for sensing and energy harvesting
  applications.
\newblock Nanotechnology \textbf{25}(37), 375603 (2014).
\newblock \doi{10.1088/0957-4484/25/37/375603}

\bibitem{kroese2011}
Kroese, D.P., Taimre, T., Botev, Z.I.: Handbook of Monte Carlo Methods.
\newblock Wiley, New~Jersey (2011)

\bibitem{Li_2019}
Li, Y., Zhou, S., Litak, G.: Uncertainty analysis of excitation conditions on
  performance of nonlinear monostable energy harvesters.
\newblock International Journal of Structural Stability and Dynamics
  \textbf{19}(06), 1950052 (2019).
\newblock \doi{10.1142/S0219455419500524}

\bibitem{Li_2020b}
Li, Y., Zhou, S., Litak, G.: Robust design optimization of a nonlinear
  monostable energy harvester with uncertainties.
\newblock Meccanica \textbf{55}, { 1753–1762} (2020).
\newblock \doi{10.1007/s11012-020-01216-z}

\bibitem{Li_2020}
Li, Y., Zhou, S., Litak, G.: Uncertainty analysis of bistable vibration energy
  harvesters based on the improved interval extension.
\newblock Journal of Vibration Engineering \& Technologies \textbf{8},
  {297--306} (2020).
\newblock \doi{10.1007/s42417-019-00134-z}

\bibitem{Lopes_2019}
Lopes, V.G., Peterson, J.V.L.L., Cunha~Jr., A.: Nonlinear characterization of a
  bistable energy harvester dynamical system.
\newblock In: M.~Belhaq (ed.) Topics in Nonlinear Mechanics and Physics, pp.
  71--88. Springer Singapore, Singapore (2019)

\bibitem{Mahmud2022}
Mahmud, M.A.P., Bazaz, S.R., Dabiri, S., Mehrizi, A.A., Asadnia, M., Warkiani,
  M.E., Wang, Z.L.: Advances in mems and microfluidics-based energy harvesting
  technologies.
\newblock Advanced Materials Technologies \textbf{7}(7), 2101347 (2022).
\newblock \doi{10.1002/admt.202101347}

\bibitem{Mallick2017}
Mallick, D., Constantinou, P., Podder, P., Roy, S.: Multi-frequency mems
  electromagnetic energy harvesting.
\newblock Sensors and Actuators A: Physical \textbf{264}, 247--259 (2017).
\newblock \doi{10.1016/j.sna.2017.08.002}

\bibitem{Mann_2010}
Mann, B., Owens, B.: Investigations of a nonlinear energy harvester with a
  bistable potential well.
\newblock Journal of Sound and Vibration \textbf{329}(9), 1215--1226 (2010).
\newblock \doi{10.1016/j.jsv.2009.11.034}

\bibitem{Martins_2022}
Martins, P., Trindade, M., Varoto, P.: Simplified robust and multiobjective
  optimization of piezoelectric energy harvesters with uncertain parameters.
\newblock International Journal of Mechanics and Materials in Design
  \textbf{18}, 63–85 (2022).
\newblock \doi{10.1007/s10999-021-09586-2}

\bibitem{Nagel2020}
Nagel, J.B., Rieckermann, J., Sudret, B.: Principal component analysis and
  sparse polynomial chaos expansions for global sensitivity analysis and model
  calibration: Application to urban drainage simulation.
\newblock Reliability Engineering \& System Safety \textbf{195}, 106737 (2020).
\newblock \doi{10.1016/j.ress.2019.106737}

\bibitem{Nanda_2015}
Nanda, A., Karami, M., Singla, P.: Uncertainty quantification of energy
  harvesting systems using method of quadratures and maximum entropy principle.
\newblock In: Proceedings of the ASME 2015 Conference on Smart Materials,
  Adaptive Structures and Intelligent Systems (2015)

\bibitem{STONEHENGE}
Norenberg, J., Peterson, J., Lopes, V., Luo, R., de~la Roca, L., Pereira, M.,
  Ribeiro, J., {Cunha~Jr}, A.: {STONEHENGE} --- {S}uite for nonlinear analysis
  of energy harvesting systems (2021).
\newblock \urlprefix\url{https://americocunhajr.github.io/STONEHENGE}

\bibitem{norenbergNoDy_2022}
Norenberg, J.P., Cunha, A., da~Silva, S., Varoto, P.S.: Global sensitivity
  analysis of asymmetric energy harvesters.
\newblock Nonlinear Dynamics \textbf{109}(2), 443--458 (2022).
\newblock \doi{10.1007/s11071-022-07563-8}

\bibitem{NORENBERG2023108542}
Norenberg, J.P., Luo, R., Lopes, V.G., Peterson, J.V.L., Cunha, A.: Nonlinear
  dynamics of asymmetric bistable energy harvesters.
\newblock International Journal of Mechanical Sciences p. 108542 (2023).
\newblock \doi{10.1016/j.ijmecsci.2023.108542}

\bibitem{STONEHENGE_paper}
Norenberg, J.P., Peterson, J.V., Lopes, V.G., Luo, R., {de la Roca}, L.,
  Pereira, M., {Telles Ribeiro}, J.G., Cunha, A.: {STONEHENGE} — {S}uite for
  nonlinear analysis of energy harvesting systems.
\newblock Software Impacts \textbf{10}, 100161 (2021).
\newblock \doi{10.1016/j.simpa.2021.100161}

\bibitem{34_Oladyshkin}
Oladyshkin, S., Nowak, W.: Data-driven uncertainty quantification using the
  arbitrary polynomial chaos expansion.
\newblock Reliability Engineering \& System Safety \textbf{106}, 179--190
  (2012).
\newblock \doi{10.1016/j.ress.2012.05.002}

\bibitem{38_Palar}
Palar, P.S., Zuhal, L.R., Shimoyama, K., Tsuchiya, T.: Global sensitivity
  analysis via multi-fidelity polynomial chaos expansion.
\newblock Reliability Engineering \& System Safety \textbf{170}, 175--190
  (2018).
\newblock \doi{10.1016/j.ress.2017.10.013}

\bibitem{PANYAM2014153}
Panyam, M., Masana, R., Daqaq, M.F.: On approximating the effective bandwidth
  of bi-stable energy harvesters.
\newblock International Journal of Non-Linear Mechanics \textbf{67}, 153--163
  (2014).
\newblock \doi{https://doi.org/10.1016/j.ijnonlinmec.2014.09.002}.
\newblock
  \urlprefix\url{https://www.sciencedirect.com/science/article/pii/S0020746214001796}

\bibitem{Ruiz2017_uq_gsa}
Ruiz, R.O., Meruane, V.: Uncertainties propagation and global sensitivity
  analysis of the frequency response function of piezoelectric energy
  harvesters.
\newblock Smart Materials and Structures \textbf{26}(6), 065003 (2017).
\newblock \doi{10.1088/1361-665X/aa6cf3}

\bibitem{SEOL_2013}
Seol, M.L., Choi, J.M., Kim, J.Y., Ahn, J.H., Moon, D.I., Choi, Y.K.:
  Piezoelectric nanogenerator with a nanoforest structure.
\newblock Nano Energy \textbf{2}(6), 1142--1148 (2013).
\newblock \doi{10.1016/j.nanoen.2013.04.006}

\bibitem{35_SEPAHVAND}
Sepahvand, K., Marburg, S., Hardtke, H.J.: Uncertainty quantification in
  stochastic systems using polynomial chaos expansion.
\newblock International Journal of Applied Mechanics \textbf{02}(02), 305--353
  (2010).
\newblock \doi{10.1142/S1758825110000524}

\bibitem{Soize_2017}
Soize, C.: Uncertainty Quantification: an accelerated course with advanced
  applications in computational engineering.
\newblock Springer, New York (2017)

\bibitem{31_sudret}
Sudret, B.: Global sensitivity analysis using polynomial chaos expansions.
\newblock Reliability Engineering \& System Safety \textbf{93}(7), 964--979
  (2008).
\newblock \doi{10.1016/j.ress.2007.04.002}.
\newblock Bayesian Networks in Dependability

\bibitem{13_triplett}
Triplett, A., Quinn, D.D.: The effect of non-linear piezoelectric coupling on
  vibration-based energy harvesting.
\newblock Journal of Intelligent Material Systems and Structures
  \textbf{20}(16), 1959--1967 (2009).
\newblock \doi{10.1177/1045389X09343218}

\bibitem{Varoto_2019}
Varoto, P.: Dynamic behavior and performance analysis of piezoelastic energy
  harvesters under model and parameter uncertainties.
\newblock In: Proceedings of the Society for Experimental Mechanics Series
  (2019)

\bibitem{Wasserman_2007}
Wasserman, L.: All of Nonparametric Statistics.
\newblock Springer, New York (2007)

\bibitem{xiu2002p619}
Xiu, D., Karniadakis, G.E.: The {W}iener-{A}skey polynomial chaos for
  stochastic differential equations.
\newblock {SIAM} Journal on Scientific Computing \textbf{24}, 619--644 (2002).
\newblock \doi{10.1137/S1064827501387826}

\end{thebibliography}
% --------------------------------------------------------------

\end{document}